\documentclass[10pt]{article}
\usepackage{latexsym}
\usepackage{tikz}
\usepackage{tkz-graph}
\usepackage{amsmath}
\usepackage{longtable}
\usepackage{array}
\usepackage{graphicx}
\graphicspath{ {./images/} }
\usepackage{caption}
\usepackage{subcaption}
\usepackage[toc,page]{appendix}
\newcommand*\circled[1]{\tikz[baseline=(char.base)]{
		\node[shape=circle,draw,inner sep=2pt] (char) {#1};}}
\usetikzlibrary{shapes}

\newtheorem{Theorem}{Theorem}[section]

\newtheorem{Lemma}[Theorem]{Lemma}
\newtheorem{Observation}[Theorem]{Observation}

\newtheorem{Conjecture}[Theorem]{Conjecture}

\newtheorem{Property}{P \hspace*{-0.7em}}

\def\inst#1{$^{#1}$}
\def\CCD{($claw$, $4K_1$, $bridge$, $C_4$-twin)}
\newcommand{\qed}{ $\Box$}

\begin{document}
	
\title{On graphs without four-vertex induced subgraphs \\}

\author{
	Kathie Cameron\inst{1}\thanks{The  authors acknowledge the support of the Natural Sciences and Engineering Research Council of Canada (NSERC), [funding reference number RGPIN-2016-06517 for the first and third authors and DDG-2024-00015 for the second author]. Cette recherche a \'et\'e financ\'ee par le Conseil de recherches en sciences naturelles et en g\'enie du Canada (CRSNG), [num\'ero de r\'ef\'erence RGPIN-2016-06517 pour le premier et le troisi\`eme auteurs, et DDG-2024-00015 pour the deuxi\`eme auteur].}
	\and Ch\'inh T. Ho\`ang\inst{2}\footnotemark[1]
	\and Taite LaGrange\inst{3}\footnotemark[1]
}
\maketitle

\begin{center}
{\footnotesize

\inst{1} Department of Mathematics, Wilfrid Laurier
University, \\Waterloo, Ontario, Canada \\
\inst{2} Department of Physics and Computer Science, Wilfrid Laurier
University, \\Waterloo, Ontario, Canada \\
\inst{3} David R. Cheriton School of Computer Science, University of Waterloo, \\Waterloo, Ontario, Canada }\\

\end{center}

\begin{abstract}
Given a family $F$ of graphs, a graph $G$ is $F$-free if it does not
contain any graph in $F$ as an induced subgraph. The problem of determining the complexity of colouring ($claw$, $4K_1$)-free graphs is a well-known open problem. In this paper we solve the colouring problem for a subclass of ($claw$, $4K_1$)-free graphs.
We design a polynomial-time algorithm to colour \CCD-free graphs. 
This algorithm is derived from a structural theorem on \CCD-free graphs.


\end{abstract}

\section{Introduction}\label{sec:intro}
{\it Vertex colouring} is the problem of assigning a minimum number of colours to the vertices of a graph $G$ so that no two adjacent vertices receive the same colour.   Vertex colouring is NP-hard for general graphs. The problem is even NP-hard for some well-studied restricted classes of graphs such as $claw$-free graphs (definitions not given here will be given later). 
However vertex colouring can be solved in polynomial time for many well-known classes of graph such as perfect graphs. 

Let $L$ be a set of graphs. A graph $G$ is $L$-free if $G$ does not contain any graph of $L$ as an induced subgraph. There has been keen interest in graphs whose forbidden list $L$ contains graphs with four vertices, and a recent paper of Lozin and Malyshev \cite{Lozin} discusses the state of the art on this problem, identifying three outstanding classes: $L = (claw, 4K_1), L = (claw, 4K1, co-diamond)$, and $L = (4K_1,C_4)$. To be more formal, if $L$ is a set of four-vertex graphs, then it is known that vertex colouring $L$-free graphs is polynomial-time solvable or NP-hard except when $L$ is one of the above three classes. We are interested in the class $L=(claw, 4K_1)$. ($4K_1$ is the stable set on four vertices.)

\begin{figure}
\begin{center}
	\begin{tikzpicture} [scale = 1]
	\tikzstyle{every node}=[font=\small]
	
	\newcommand{\size}{1}
	
	\newcommand{\claw}{1}{
		\path (0,0) coordinate (g1);
		\path (g1) +(\size / 2, \size / 2) node (g1_1){}; \path (g1) +(\size / 2, 1.5 * \size) node
		(g1_2){}; \path (g1) +(-\size / 2, 0) node (g1_3){}; \path (g1)
		+(1.5 * \size, 0) node (g1_4){};
		\foreach \Point in {(g1_1),(g1_2),(g1_3),(g1_4)}{
			\node at \Point {\textbullet};
		}
		\draw   (g1_1) -- (g1_2)
		(g1_1) -- (g1_3)
		(g1_1) -- (g1_4);
		\path (g1) ++(\size  / 2,-\size / 2) node[draw=none,fill=none] { {\large Claw}};
	}

	\newcommand{\g2}{2}{
		\path(3 * \size, 0) coordinate(g2);
		\path (g2) +(\size / 2, 0) node(g2_1){};
		\path (g2) +(- \size / 2, \size / 2) node(g2_2){};
		\path (g2) +(1.5 * \size, \size / 2) node(g2_3){};
		\path (g2) +(\size / 2, \size) node(g2_4){};
		\path (g2) +(- \size / 2, 1.5 * \size ) node(g2_5){};
		\path (g2) +(1.5 * \size, 1.5 * \size) node(g2_6){};

		\foreach \Point in {(g2_1),(g2_2),(g2_3),(g2_4), (g2_5), (g2_6)}{
			\node at \Point {\textbullet};
		}
		\draw   (g2_1) -- (g2_2)
		(g2_1) -- (g2_3)
		(g2_1) -- (g2_4)
		(g2_4) -- (g2_2)
		(g2_4) -- (g2_3)
		(g2_2) -- (g2_5)
		(g2_3) -- (g2_6);
		\path (g2) ++(\size  / 2,-\size / 2) node[draw=none,fill=none] { {\large {\it $P_5$-Twin}}};
	}
	
	\newcommand{\g3}{3}{
		\path( 6.5 * \size, 0) coordinate(g3);
		\path(g3) +(\size / 2, 0) node(g3_1){};
		\path(g3) +(0, \size) node(g3_2){};
		\path(g3) +(\size, \size) node(g3_3){};
		\path(g3) +(\size / 2, 2 * \size) node(g3_4){};
		\path(g3) +(-1 * \size, \size) node(g3_5){};
		\path(g3) +(2 * \size, \size) node(g3_6){};
		
		\foreach \Point in {(g3_1),(g3_2),(g3_3),(g3_4), (g3_5), (g3_6)}{
			\node at \Point {\textbullet};
		}
		\draw   (g3_1) -- (g3_2)
		(g3_1) -- (g3_3)
		(g3_1) -- (g3_4)
		(g3_4) -- (g3_2)
		(g3_4) -- (g3_3)
		(g3_1) -- (g3_5)
		(g3_2) -- (g3_5)
		(g3_4) -- (g3_5)
		(g3_3) -- (g3_6);
		\path (g3) ++(\size  / 2,-\size / 2) node[draw=none,fill=none] { {\large co-$R $}};
		
	}

	\newcommand{\g4}{4}{
		\path(0, -3 * \size) coordinate(g4);
		\path (g4) +(\size / 2, 0) node(g4_1){};
		\path (g4) +(- \size / 2, \size / 2) node(g4_2){};
		\path (g4) +(1.5 * \size, \size / 2) node(g4_3){};
		\path (g4) +(\size / 2, \size) node(g4_4){};
		\path (g4) +(\size / 2, 2 * \size ) node(g4_5){};

		\foreach \Point in {(g4_1),(g4_2),(g4_3),(g4_4), (g4_5)}{
			\node at \Point {\textbullet};
		}
		\draw   (g4_1) -- (g4_2)
		(g4_1) -- (g4_3)
		(g4_4) -- (g4_2)
		(g4_4) -- (g4_3)
		(g4_2) -- (g4_5)
		(g4_3) -- (g4_5)
		(g4_1) -- (g4_4);
		\path (g4) ++(\size  / 2,-\size / 2) node[draw=none,fill=none] { {\large $C_4$-twin}};
	}
	
	\newcommand{\g5}{5}{
		\path(3 * \size, -3 * \size) coordinate(g5);
		\path (g5) +(\size / 2, 0) node(g5_1){};
		\path (g5) +(- \size / 2, \size / 2) node(g5_2){};
		\path (g5) +(1.5 * \size, \size / 2) node(g5_3){};
		\path (g5) +(\size / 2, \size) node(g5_4){};
		\path (g5) +(- \size / 2, 1.5 * \size ) node(g5_5){};
		\path (g5) +(1.5 * \size, 1.5 * \size) node(g5_6){};

		\foreach \Point in {(g5_1),(g5_2),(g5_3),(g5_4), (g5_5), (g5_6)}{
			\node at \Point {\textbullet};
		}
		\draw   (g5_1) -- (g5_2)
		(g5_1) -- (g5_3)
		(g5_4) -- (g5_2)
		(g5_4) -- (g5_3)
		(g5_2) -- (g5_5)
		(g5_3) -- (g5_6)
		(g5_5) -- (g5_6)
		(g5_1) -- (g5_4);
		\path (g5) ++(\size  / 2,-\size / 2) node[draw=none,fill=none] { {\large $C_5$-twin}};
	}
	
	\newcommand{\g6}{6}{
		\path( 6.5 * \size, -3 * \size) coordinate(g6);
		\path(g6) +(\size / 2, 0) node(g6_1){};
		\path(g6) +(0, \size) node(g6_2){};
		\path(g6) +(\size, \size) node(g6_3){};
		\path(g6) +(\size / 2, 2 * \size) node(g6_4){};
		\path(g6) +(-1 * \size, \size) node(g6_5){};
		\path(g6) +(2 * \size, \size) node(g6_6){};
		
		\foreach \Point in {(g6_1),(g6_2),(g6_3),(g6_4), (g6_5), (g6_6)}{
			\node at \Point {\textbullet};
		}
		\draw   (g6_1) -- (g6_2)
		(g6_1) -- (g6_3)
		(g6_1) -- (g6_4)
		(g6_4) -- (g6_2)
		(g6_4) -- (g6_3)
		(g6_1) -- (g6_5)
		(g6_2) -- (g6_5)
		(g6_4) -- (g6_5)
		(g6_1) -- (g6_6)
		(g6_3) -- (g6_6)
		(g6_4) -- (g6_6);
		\path (g6) ++(\size  / 2,-\size / 2) node[draw=none,fill=none] { {\large bridge}};
		
	}
	
	\newcommand{\g7}{7}{
		\path(\size / 2, -6 * \size) coordinate(g7);
		\path(g7) +(0, \size) node(g7_1){};
		\path(g7) +(\size / 2, 0) node(g7_2){};
		\path(g7) +(-\size / 2 , 0) node(g7_3){};
		\path(g7) +(-6 / 5 * \size, 6/5 * \size) node(g7_4){};
		\path(g7) +(0, 2 * \size) node(g7_5){};
		\path(g7) +(6 / 5 * \size, 6/5 * \size) node(g7_6){};
		\foreach \Point in {(g7_1),(g7_2),(g7_3),(g7_4), (g7_5), (g7_6)}{
			\node at \Point {\textbullet};
		}
		\draw   (g7_1) -- (g7_2)
		(g7_1) -- (g7_3)
		(g7_1) -- (g7_4)
		(g7_1) -- (g7_5)
		(g7_1) -- (g7_6)
		(g7_2) -- (g7_3)
		(g7_2) -- (g7_6)
		(g7_3) -- (g7_4)
		(g7_4) -- (g7_5)
		(g7_5) -- (g7_6);
		\path (g7) ++(0,-\size / 2) node[draw=none,fill=none] { {\large $5-wheel$}};
		
	}
	
	\newcommand{\g8}{8}{
		\path(\size * 3.5, -6 * \size) coordinate(g8);
		\path(g8) +(-\size, 1.5 * \size) node(g8_1){};
		\path(g8) +(0, 1.5 *  \size) node(g8_2){};
		\path(g8) +(\size , 1.5 * \size) node(g8_3){};
		\path(g8) +(-\size / 2, 0) node(g8_4){};
		\path(g8) +(\size / 2, 0) node(g8_5){};
		\path(g8) +(1.5 * \size, 0) node(g8_6){};
		\foreach \Point in {(g8_1),(g8_2),(g8_3),(g8_4), (g8_5), (g8_6)}{
			\node at \Point {\textbullet};
		}
		\draw   (g8_1) -- (g8_2)
		(g8_1) -- (g8_4)
		(g8_2) -- (g8_3)
		(g8_2) -- (g8_4)
		(g8_2) -- (g8_5)
		(g8_3) -- (g8_5)
		(g8_3) -- (g8_6)
		(g8_4) -- (g8_5)
		(g8_5) -- (g8_6);
		\path (g8) ++(\size  / 3.5,-\size / 2) node[draw=none,fill=none] { {\large co-$A$}};
	}
	
	\newcommand{\g9}{9}{
		\path(\size * 7.0, -6 * \size) coordinate(g9);
		\path(g9) +(0, 0) node(g9_1){};
		\path(g9) +(-\size, \size/2) node(g9_2){};
		\path(g9) +(\size, \size/2) node(g9_3){};
		\path(g9) +(-\size, 1.5 * \size) node(g9_4){};
		\path(g9) +(\size, 1.5 * \size) node(g9_5){};
		\foreach \Point in {(g9_1),(g9_2),(g9_3),(g9_4), (g9_5)}{
			\node at \Point {\textbullet};
		}
		\draw   (g9_1) -- (g9_2)
		(g9_1) -- (g9_3)
		(g9_1) -- (g9_4)
		(g9_1) -- (g9_5)
		(g9_2) -- (g9_3)
		(g9_2) -- (g9_4)
		(g9_2) -- (g9_5)
		(g9_3) -- (g9_4)
		(g9_3) -- (g9_5);
		\path (g9) ++(0,-\size / 2) node[draw=none,fill=none] { {\large $K_5 - e$}};
	}

	\end{tikzpicture}
\end{center}
	\caption{The nine forbidden sugraphs for line-graphs}\label{fig:line-graphs}
\end{figure}
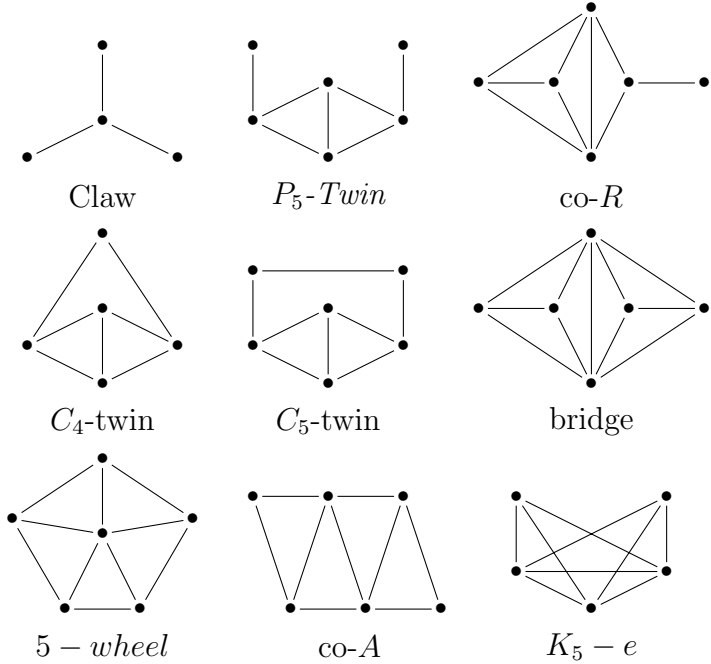

Given a graph $G$, the line-graph $L(G)$ of $G$ is defined to be the graph whose vertices are the edges of $G$, and two vertices of $L(G)$ are adjacent if their corresponding edges in $G$ are incident. Beineke \cite{Bei} derived a characterization of line-graphs in terms of a set of nine forbidden induced subgraphs (see Figure~\ref{fig:line-graphs}), two of which are the $claw$ and $K_5 - e$. Fraser, Hamel, Ho\`ang and Maffray \cite{FraHam2018} showed that vertex colouring can be solved in polynomial time for $(claw, 4K_1, K_5 - e)$-free graphs. Given the connection between $(claw, 4K_1)$-free graphs and $(4K_1)$-free line-graphs, it is natural to ask about polynomial-time solvability of vertex colouring of graphs where $K_5 - e$ is replaced by one of the other nine forbidden induced subgraphs. Three of the graphs in Beineke's list contain a $K_4$ (a clique on four vertices). They are $K_5 -e$, co-$R$, and $bridge$. The following problem then might be of interest: Is there a polynomial-time algorithm for vertex colouring $(claw, 4K_1, F)$-free graphs, where $F$ is the graph co-$R$ or the $bridge$? Abuadas and Ho\`ang \cite{AbuHoa2022} positively solved the above problem for $F$ = co-$R$.  In this paper, we investigate the question for the $bridge$. We are not able to fully answer the question for the $bridge$. However, we are able to provide the answer if we add the $C_4$-twin to the list of forbidden subgraphs. Note that $C_4$-twin is one of the nine forbidden induced subgraphs for line-graphs. Furthermore, $C_4$-twin-free graphs do not contain complements of chordless cycles with at least seven vertices. In fact, $C_4$-twin is the only graph of the nine forbbiden induced subgraphs with this property. This fact allows us to use known perfect graph algorithms to establish our main result.

\begin{Theorem}\label{thm:main}
	There is a polynomial-time algorithm to colour a {\CCD}-free graph.
\end{Theorem}

The purpose of this paper is to prove Theorem \ref{thm:main}.
Our results have connections to perfect graphs and clique-width. As we will see later, the theorem below implies Theorem~\ref{thm:main}.
\begin{Theorem}\label{thm:structure}
	Let $G$ be a  \CCD-free graph. Then one of the following holds.
	\begin{description}
		\item[(i)] $G$ is perfect. 
		\item[(ii)] $G$ has bounded clique-width.
	\end{description}
\end{Theorem}

In Section~\ref{sec:background}, we discuss background results related to our work and establish a number of new results on clique-width. In Sections \ref{sec:proofs}-\ref{sec:c5}, we study the structure of {\CCD}-free graphs. In Section~\ref{sec:proof-of-main-results}, we prove Theorems~\ref{thm:main} and~\ref{thm:structure}. Finally, in Section~\ref{sec:conclusions}, we discuss open problems related to our work.

\section{Definitions and Background}\label{sec:background}
In this section, we discuss the background results needed to prove our main theorem. This section has two subsections. In the first subsection, we discuss $claw$-free graphs. In the second subsection, we discuss clique-width of graphs.
\subsection{Claw-free graphs}\label{subsec:claw}
Let $\chi(G)$ denote the chromatic number of $G$ (the minimum number of colours for which a graph can be coloured), and $\omega(G)$ denote the number of vertices in a largest clique of $G$. 
A {\it hole} is an induced cycle on at least four vertices. An {\it antihole} is the complement of a hole. A $k$-hole is a hole with $k$ vertices.
Our results rely on known theorems on perfect $claw$-free graphs, and we discuss these results now.  A graph G is {\it perfect} if for each induced subgraph $H$ of $G$, we have $\chi(H) = \omega(H)$. A graph is {\it Berge} if it does not contain an odd hole or odd antihole as an induced subgraph. Parthasarathy and Ravindra \cite{ParRav1976} proved that $claw$-free Berge graphs are perfect. Chv\'atal and Sbihi \cite{ChvSbi1998} showed that $claw$-free perfect graphs can be recognized in polynomial time. Hsu \cite{Hsu1981} showed that these graphs can be coloured in polynomial time. Chudnovsky, Robertson, Seymour and Thomas \cite{ChuRob2006} proved that a graph is perfect if and only if it is Berge, solving a long standing conjecture of Berge \cite{Ber1961}. Perfect graphs can be optimally coloured in polynomial time (Gr\"otschel, Lov\'asz and Schrijver \cite{GroLov1984}). For more information on perfect graphs, see Berge and Chv\'atal \cite{BerChv1984}, or Ho\`ang and Sritharan \cite{HoaSri2015}.  In the next section, we discuss the notion of the clique-width of graphs.

\subsection{Clique-width of graphs}\label{subsec:cliquewidth}
In this section we discuss clique-width and prove a number of results that we will need later.

The clique-width of a graph G is the smallest number of labels needed to construct G using the following four operations:
\begin{itemize}
	\item
	Create a new vertex $v$ with integer $i$, denoted by $i(v)$,
	\item
	Take the disjoint union of two graphs $G$ and $H$, denoted by $G \oplus H$,
	\item
	Relabel all vertices with label $i$ by label $j$, denoted by $\rho_{i\rightarrow j}$,
	\item
	Join every vertex labeled $i$ to every vertex labeled $j$ by an edge, denoted by $\eta_{i,j}$. 
\end{itemize}
We will refer to the above four operations as ``the four clique-width operations''. Every graph can be defined by an expression using the four clique-width operations. For example, the chordless path on three vertices $v_1, v_2, v_3$ can be built with the following operations:
\newline 
\[\eta_{2,3}(( \eta_{1,2}( 1(v_1) \oplus 2(v_2)) ) \oplus (3(v_3)))\]

Such expressions are called $k$-expressions where $k$ is the number of labels. Thus the clique-width   of a graph $G$, denoted by $cwd(G)$, is the smallest $k$ such that there is a $k$-expression defining $G$. The above expression shows that $cwd(P_4) \leq 2$, where $P_4$ is the induced path on four vertices.

It is well-known \cite{CouOla2000} that if the clique-width of 
a graph is bounded then so is that of its complement. In the theorem below, $\overline{G}$ denotes the complement of graph $G$. 
\begin{Theorem}[Theorem 4.1 in \cite{CouOla2000}]\label{thm:complement}
For a graph $G$, $cwd(\overline{G}) \leq 2 \; cwd(G)$.
\end{Theorem}

Let $A$ and $B$ be two disjoint sets of vertices of a graph $G$, we write $A \;\circled{0}\; B$ to mean there are no edges between $A$ and $B$ (the {\it co-join} of $A$ and $B$). We write $A \;\circled{1}\; B$ to mean there are all edges between $A$ and $B$ (the {\it join} of $A$ and $B$). A set $M$ of vertices of a graph $G$ is a {\it module} if every vertex $x$ in $V(G) - M$ satisfies $\{x\} \;\circled{0}\; M$ or  $\{x\} \;\circled{1}\; M$. Trivially, any vertex is a module, and the set $V(G)$ is also a module of $G$. We are interested in ``non-trivial'' modules. These are called homogeneous sets defined as follows. A subset $X$ of $V(G)$ is {\it homogeneous} if $X$ is a module of $G$ and $2 \leq |X| < |V(G)|$.  A graph is {\it prime} if it does not contain a homogenous set.
The following results are folklore. 

\begin{Observation}\label{obs:constant}
	Let $G$ be a graph. Let $k$ be a constant. Then the graph obtained from $G$ by removing any $k$ vertices has bounded clique-width if and only if $G$ has bounded clique-width.
\end{Observation}

\begin{Observation}[\cite{CouMak2000}]\label{obs:prime}
	Let $G$ be a graph. If every prime induced subgraph of $G$ has bounded clique-width, then $G$ has bounded clique-width.
\end{Observation}

Finally, we will need to following result in the literature.
\begin{Lemma}[\cite{DHH}]\label{lem:cliquewidth-2k}
	Let $G$ be a graph such that $V(G)$ can be partitioned into $k$ (disjoint) cliques $X_1, \ldots, X_k$. For a vertex $x$, let $X_{i_x}$ be the clique containing $x$, and  let $N_O(x)$ be the set of neigbours $y$ of $x$ such that $y \in X_j$ for some $j \not= i_x$. Suppose $G$ satisfies the following conditions: (i) for every vertex $x$ and any set $X_j$ with $j \not= i_x$, $x$ has at most one neighbor in $X_j$, and (ii) for any vertex $x$, $N_O(x)$ is a clique. Then $G$ has clique-width at most $2k$. Furthermore, there is a labelling of $G$ such that the vertices of each $X_i$ have labels that are different
     from those of $X_j$'s with $i \not= j$. \qed
\end{Lemma}

In closing this section, we would like to note the results of \cite{BraEng2006}. Let $L^4$ be a set of four-vertex graphs. Brandst\"adt, Engelfried, Le and Lozin \cite{BraEng2006} obtained a complete description of all classes $L^4$ such that $L^4$-free graphs have bounded (or unbounded) clique-width. In particular, $(claw, 4K_1)$-free graphs have unbounded clique-width.   This result is one of the motivations for our work.

\subsection{Labelling Conventions}

To show that a given graph has bounded clique-width, we will describe a labelling algorithm using the four clique-width operations. We will partition the vertices of a graph into sets with names such as $T_i$, $X_i$, or $Y_i$. We will need to assign labels to these vertices. The labels will be 2-tuples such as ($s$, new) or ($s$, old) where $s$ is a name that indicates the set a given vertex is in. In the beginning of the labelling process, a vertex $v$ receives a ``new'' label such as ($s$, new). The label of $v$ may change during the process, but it will eventually change to a ``final'' label which will be ($s$, old). So a vertex in $X_1$ will have label (x1, old) and a vertex in $Y_2$ will have label (y2, old). When we say a set is labelled, we mean all vertices of the set have this label (set, old).

Once all sets $T_i$, $X_i$, and $Y_i$ are labelled, the induced cycle $H$ defined in each section can be constructed by creating and labelling each vertex with a unique label $1$ to $k$ where $k$ is the length of the cycle, then joining the vertices in the sets $T_i$, $X_i$, and $Y_i$ with the relevant labelled vertices of the cycle. For example, we would join the vertices with label (x1, old) to the vertices with labels $1$, $2$, and $3$.

Since we can always add the induced cycle $H$ in this way after the sets $T_i$, $X_i$, and $Y_i$ are labelled, we will ignore the vertices of the induced cycle $H$ when providing the labelling algorithms.

When we refer to a vertex in a particular set $T_i$, $X_i$, or $Y_i$ as having neighbours `outside its own set', we are referring to a vertex in a given set, say $X_1$, having neighbours in any of the other sets $T_i$, $X_i$, or $Y_i$ and not its neighbours in the cycle $H$, since we have established that we will ignore $H$ when labelling.

\begin{Lemma}[\cite{Dai2022}]\label{lem:pair-labelling}
	Given two cliques $S$ and $A$, if a vertex in $S$ is adjacent to at most one vertex in $A$ and vice versa, then the graph $G[S \cup A]$ has bounded clique-width. We will refer to this as `labelling via pairs'.
\end{Lemma}	
\noindent {\it Proof}. Label a vertex in $S$ with label (s, new). If it has a neighbour in $A$, label that neighbour (a, new). Make the join between vertices with label (s, new) and (a, new). Relabel these vertices (s, old) and (a, old), respectively.

Label the next vertex in $S$ with label (s, new). If it has a neighbour in $A$, label that neighbour (a, new). Make the join between vertices with label (s, new) and (a, new). This time, also make the join between vertices with label (s, new) and (s, old) and vertices with (a, new) and (a, old) to construct the cliques $S$ and $A$. Repeat this step for all vertices with neighbours outside their own clique.

Now the remaining unlabelled vertices are those with no neighbours outside their own clique. Give all remaining vertices in $S$ the label (s, new), make the join between vertices with label (s, new) and (s, old), and relabel the vertices (s, new) with label (s, old). Do the same for $A$ with the labels (a, new) and (a, old).

Then the graph has clique-width at most $4$, hence its clique-width is bounded.\qed

\begin{Lemma}[\cite{Dai2022}]\label{lem:nonpair-labelling}
	Given two cliques $S$ and $A$, if a vertex in $S$ is non-adjacent to at most one vertex in $A$ and vice versa, then the graph $G[S \cup A]$ has bounded clique-width. We will refer to this as `labelling via pairs'.
\end{Lemma}	
\noindent {\it Proof}. We can label the graph $G[S \cup A]$ via pairs similarly to what is done in Lemma \ref{lem:pair-labelling}. However, when labelling a non-adjacent pair of vertices (s, new) and (a, new), do not make the join between vertices with label (s, new) and vertices with label (a, new). Instead, make the join between vertices with labels (s, new) and (a, old) and between vertices with labels (a, new) and (s, old).

Then the graph has clique-width at most $4$, hence its clique-width is bounded.\qed

\noindent{\it Note:} Lemmas \ref{lem:pair-labelling} and \ref{lem:nonpair-labelling} were proved in \cite{Dai2022}, but we provide them again here for completeness and clarity of terms used to refer to them.

\begin{Lemma}\label{lem:row-labelling}
	Let $G$ be a graph such that its vertices can be partitioned into $t$ disjoint cliques $S_1, \ldots, S_t$ ($t \geq 3$) such that for each $i$, (a) a vertex in $S_i$ is adjacent to at most one vertex in $S_{i+1}$ and vice versa, where $i$ is taken modulo $t$, or a vertex in $S_i$ is non-adjacent to at most one vertex in $S_{i+1}$ and vice versa and (b) $S_i$ and $S_k$ form a join or co-join for $k\not= i-1, k\not= k+1$. Then $G$ has clique-width at most $2t+1$. We will refer to this as `labelling via rows'.
\end{Lemma}	
\noindent {\it Proof}. Let $a$ be a vertex in $S_i$. If every vertex in $S_i$ is adjacent (respectively, non-adjacent) to at most one vertex in $S_{i+1}$, then let $a'$ be the vertex in $S_{i+1}$ that is adjacent (respectively, non-adjacent) to $a$. We say that $a$ and $a'$ are {\it partners}. A vertex may have two partners, each of them in a different set $S_j$.

We now explain the labelling algorithm for $G$. Construct a graph $G'$ from $G$ on the same vertex set. Two vertices of $G'$ are adjacent if and only if they are partners in $G$. Then $G'$ is a union of paths and cycles. Initially, the vertices of $G$ are not labelled. Let $P$ be a path or cycle of $G'$. Without loss of generality, let the vertices of $P$ be $v_1, v_2, \ldots, v_r$ where $v_j v_{j+1}$ is an edge for $j = 1, \ldots, r$, and $v_r v_1$ is an edge if $P$ is a cycle. Let $i_j$ be the index of the set $S_{i_j}$ that contains $v_i$. First, for $k=1,2,3$, we label $v_k$ with label ($i_k$, new). If $v_1v_2$ is an edge of $G$, then make the join between vertices of label ($i_1$, new) and  ($i_2$, new). Repeat the same process with the edge (or non-edge) $v_2 v_3$. Now relabel $v_2$ with label (2, old). Now, suppose that we have labelled vertices $v_1, \ldots, v_k$ of $P$ such that $v_1, v_k$ have new labels and $v_2, \ldots, v_{k-1 } $ have old labels. We will extend the labelling of this partial part. 

If $v_1$ and $v_k$ are partners in $G$, then we make a join or co-join between the labels ($i_1$, new) and ($i_k$, new), depending on whether $v_1 v_k$ is an edge or not. Then we relabel ($i_1$, new) to ($i_1$, old) and ($i_k$, new) to ($i_k$, old). We have labelled a cycle of $G'$. Now we label the next unlabelled path or cycle of $G'$ (if it exists).

If $v_k$ has no partner in $S_{i_{k +1}}$, then we have reached the end of the path $P$. We relabel ($i_1$, new) to ($i_1$, old) and relabel ($i_k$, new) to ($i_k$, old). We have labelled a path of $G'$. Now we label the next unlabelled path or cycle of $G'$ (if it exists).

Now, we consider the case $v_k$ has an unlabelled partner. Let $v_{k+1}$ be the unlabelled partner of $v_k$. If $v_{k+1} \not \in S_{i_1}$, then we continue our labelling process, that is, we label  $v_{k+1}$ with label ($i_{k+1}$, new) and relabel $v_k$ with label ($i_k$, old). Now, suppose that $v_{k+1} \in S_{i_1}$. At this point, we know that $S_{i_{1}}$ already has a new label with the vertex $v_1$. We need a second `new' label for $v_{k+1}$. We label $v_{k+1}$ with label ($i_1$, new2). We make a join or co-join between the labels ($i_k$, new) and ($i_{k+1}$, new2), depending on whether $v_k v_{k+1}$ is an edge or not. Relabel $v_k$ with label ($i_k$, old). We continue the process with the vertex $v_{k+1}$. Note that $S_{i_1}$ has at most two `new' labels for the entirely of the labelling process. 

After we label all paths and cycles of $G'$, we may assume that the vertices of each $S_k$ are labelled with label  ($k$, old). Now, for each pair of $S_i$ and $S_k$ such that $S_i \circled{1} S_k$, we make a join between labels ($i$, old) and ($k$, old).

This proves that $G$ has clique-width at most $2t+1$.\qed

\subsection{The methodology}
We will prove Theorem~\ref{thm:main} by using clique-width. In particular, our main result will follow from the following theorems that will be proved in later sections.

\noindent 
{\bf Theorem \ref{thm-c7}}
{\it  Let $G$ be a ($claw$,$4K_1$,bridge)-free graph. If $G$ contains a $C_7$, then $G$ has bounded clique-width.}

\noindent
{\bf Theorem \ref{thm-c6}}
{\it 	Let $G$ be a ($claw$,$4K_1$,bridge, $C_4$-twin)-free graph. If $G$ contains a $C_6$, then $G$ has bounded clique-width.}

\noindent
{\bf Theorem \ref{thm-c5}}
{\it 	Let $G$ be a ($claw$,$4K_1$,bridge, $C_4$-twin)-free graph. If $G$ contains a $C_5$, then $G$ has bounded clique-width.}



An induced cycle $C_t$ where $t \ge 8$ contains a $4K_1$. 
The complement of an induced cycle $\overline{C_t}$ where $t \ge 7$ contains a $C_4$-twin. 
So a ($claw$,$4K_1$,bridge, $C_4$-twin)-free graph is either perfect or contains an induced $C_5$ or $C_7$. 
Thus, Theorems~\ref{thm-c7}, \ref{thm-c6}, \ref{thm-c5} imply Theorem~\ref{thm:main}.

\section{Overview of the proofs}\label{sec:proofs}

We investigate two subclasses of $(claw, 4K_1, bridge, C_4$-$twin)$-free graphs: those which contain a $C_5$ and those which contain a $C_7$. In the case of $C_7$, we do not need to exclude the $C_4$-twin. We show that each of these subclasses has bounded clique-width.
When a $(claw, 4K_1, bridge, C_4$-$twin)$-free graph contain a $C_5$, it will sometimes contain a $C_6$. Thus, before we investigate the case of the $C_5$, we will first discuss the case where the graph contains a $C_6$. Doing so allows us to forbid the $C_6$ while investigating the $C_5$. This makes the case of the $C_5$ more manageable.

\subsection{When the graph contains a $C_7$}

Let $G$ be a $(claw, 4K_1, bridge)$-free graph. Assume $G$ contains an induced $C_7$, call it $H$, and let its vertices be $1,2, \dots, 7$, in that order. Throughout this section, indices are taken mod 7. 

We define a fixed number of sets (here 21) which partition the vertices of $G \setminus H$. We prove that each set is a clique, and we describe the edges between the sets. This allows us to then prove that $G$ has bounded clique-width.

For $i \in \{1,2,3,4,5,6,7\}$:

Let $X_i$ be the set of vertices of $G$ that are adjacent to $i, i+1, i+2$. Let $X$ = $X_1 \cup \dots \cup X_7$.

Let $Y_i$ be the set of vertices of $G$ that are adjacent to $i, i+1, i+2, i+3.$ Let $Y$ = $Y_1 \cup \dots \cup Y_7$.

Let $Z_i$ be the set of vertices of $G$ that are adjacent to $i, i+1, i+3, i+4.$ Let $Z$ = $Z_1 \cup \dots \cup Z_7$.\\


When we say a set is {\it big}, we mean it contains at least five vertices. 

\begin{Lemma}\label{lem:big-sets}
    All sets are big.
\end{Lemma}

\noindent {\it Proof}. Let $G$ be a graph. Let $k$ be a constant. It is known that the graph obtained from $G$ by removing any $k$ vertices has bounded clique-width if and only if G has bounded clique-width. Then if any set has fewer than $k$ vertices, we can remove those $k$ vertices and assume the set is empty.\qed

The value of $k$ can be arbitrary. In this paper, we assign it the value 5. Thus, we assume all sets contain at least five vertices,  otherwise they are considered to be empty. 
We will deal with only a constant number of sets.

We first prove that any vertex in $G$ must belong to one of $H$, $X_i$, $Y_i$ or $Z_i$. Under the assumption that the sets are big, the $Z_i$ are all empty, so any vertex in $G$ must actually belong to one of $H$, $X_i$ or $Y_i$ (14 sets not counting $H$). 

\begin{figure}[ht]
    \centering
    \includegraphics{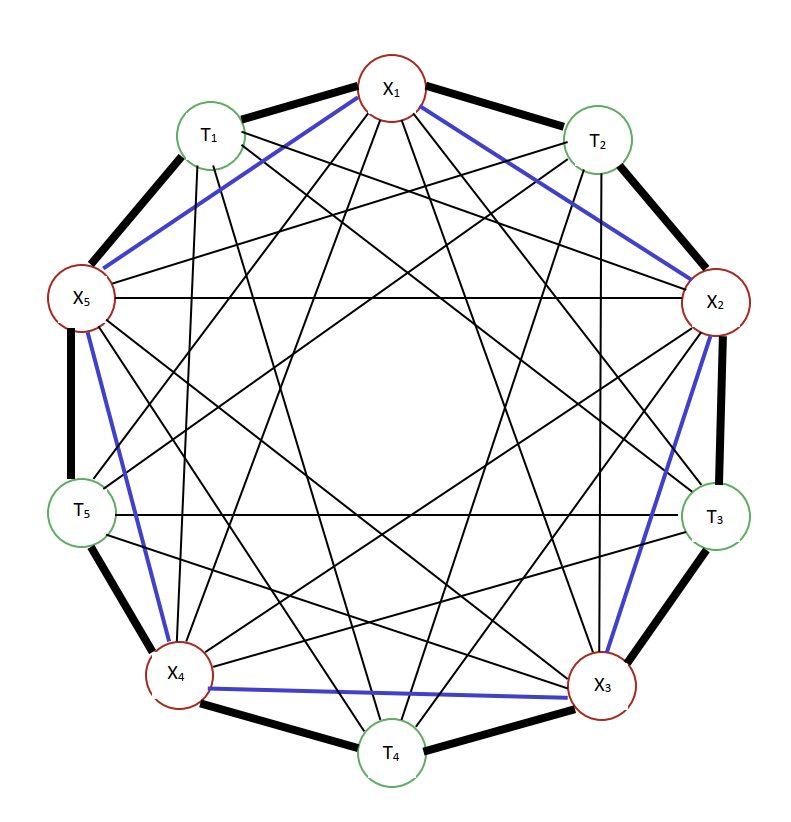}
  \caption{This figure illustrates the the structure of $G \setminus H$ when $G$ contains a $C_7$. Thick lines indicate a join, thin lines indicate that each vertex in one set has at most one neigbhour in the other, blue lines indicate that each vertex in one set has at most one non-neighbour in the other, no line indicates a cojoin, and dashed lines indicate that adjacency is unknown, however, by Properties P5, P7 and P10, not all these sets are non-empty.}
 \label{fig:c7-general}
\end{figure}

\begin{figure}[ht]
    \centering
    \includegraphics{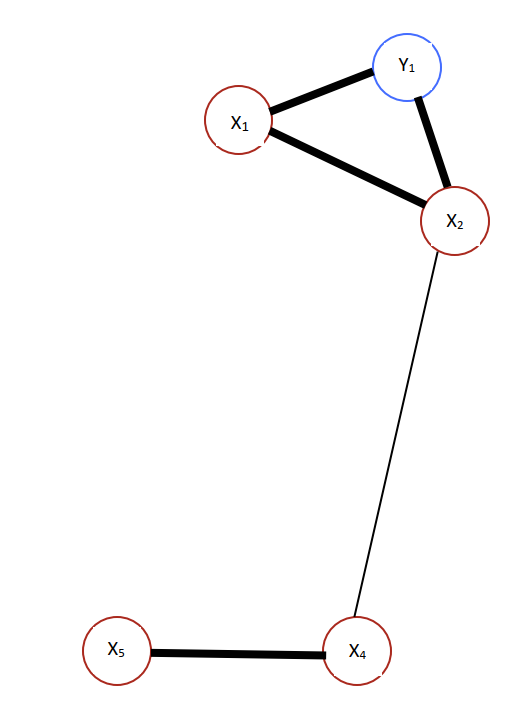}
  \caption{This figure illustrates simplifications to Figure \ref{fig:c7-general} when each of  $X_1, X_2, X_4, X_5,Y_1 \neq \emptyset$. All other sets turn out to be empty.}
 \label{fig:c7-y1nonempty}
\end{figure}

We then study the edges between the sets. See Table \ref{tab:C7} for a summary. Usually, pairs of the sets form a join or co-join. Otherwise, each vertex in one set has at most one neighbour in another particular set or each vertex in one set has at most one non-neighbour in the other  set. Often, one set being non-empty implies that others are empty. From this structure, we prove that $G[X]$ and $G[Y]$ have bounded clique-width. Any set $Y_i$ forms a join with $X_i$ and with $X_{i+1}$. For all other $X_j$, either $Y_i$ form a cojoin with $X_j$ or $X_j$ must be empty if $Y_i$ is non-empty. See Figure \ref{fig:c7-general} for the edges between sets and Figure \ref{fig:c7-y1nonempty} for an example of  simplications when several sets are non-empty. In general, to construct $G[X\cup Y]$, first take the disjoint union of $G[X]$ and $G[Y]$, and then make the required joins. Since the 7-hole $H$ has a fixed number of vertices, bounded clique-width of $G[X\cup Y]$ implies bounded clique-width of $G[X\cup Y \cup H]$. 

\subsection{When the graph contains a $C_6$}

The cases when $G$ contains a $C_6$ or a  $C_5$ are more complex than when it contains a $C_7$, so we found it necessary to forbid another induced subgraph. We chose the $C_4$-$twin$ because it is also one of the nine forbidden induced subgraphs for line-graphs and because, as mentioned before, forbidding it excludes all complements of induced cycles of size at least 7, so in particular, complements of odd induced cycles.

Let $G$ be a $(claw, 4K_1, bridge, C_4$-$twin)$-free graph. Assume $G$ contains an induced $C_6$, call it $H$, and let its vertices be $1,2, \dots, 6$, in that order. Throughout this section, indices are taken mod 6. By Theorem \ref{thm-c7}, we can assume that $G$ is also $C_7$-free. 

For $i=1,2,3,4,5,6,$ We consider the sets $X_i, Y_i,$ and $Z_i$ defined above, and let $X,Y,$ and $Z$ be the union of all $X_i, Y_i,$ and $Z_i$, respectively. 
We also need to define six more sets: 
Let $T_i$ be the set of vertices of $G$ that are adjacent to $i, i+1$. 
Let $T$ = $T_1 \cup \dots \cup T_6$.

We will show that the sets $X_i, Y_i, Z_i$ and $T_i$ partition the vertices of $G \setminus H$. As in the case of the $C_7$, under the assumption that each set is big, we show that all the sets $Z_i$ are empty, so 18 sets partition $G \setminus H$.  We prove that each of the sets $X_i, Y_i,$ and $T_i$ is a clique and we describe the edges between the sets. See Table 3 
for a summary. Again, the edges between sets are usually joins or cojoins, or each vertex in one set has at most one neighbour in another particular set or each vertex in one set has at most one non-neighbour in the other set. Often, one set being non-empty implies that others are empty or that two other sets are related by a join, cojoin, each vertex of one having at most one neighbour or one non-neighbour in another set. Property P36 is different: It states that each vertex in $Y_i$ has at most 2 non-neighbours in $T_{i+1}$.

The structure of $G[T]$ turns out to be fairly simple - at most three sets $T_i$ can be non-empty, they must be $T_i, T_{i+2}$ and $T_{i+4}$, and a vertex of $T_i$ has at most one neighbour in each of $T_{i+2}$ and $T_{i+4}$, and if these neighbours exist, they are adjacent.

We next examine the structure of $G[X]$, $G[Y]$, $G[X \cup T]$ and $G[Y \cup T]$, and $G[X \cup Y]$, showing each has bounded clique-width. When studying $G[X \cup Y \cup T]$, we show that at most two $Y_i$s can be non-empty, and if two are non-empty they must be $Y_i, Y_{i+1}$ or $Y_i, Y_{i+3}$. 

\begin{figure}[ht]
    \centering
    \includegraphics{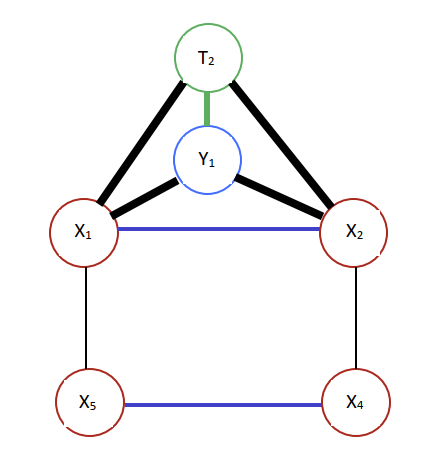}
  \caption{This figure illustrates the the structure of $G \setminus H$ when $G$ contains a $C_6$ and $Y_1$ is the only non-empty $Y_i$ and $T_2$ is the only non-empty $T_i$. Thick lines indicate a join, thin lines indicate that each vertex in one set has at most one neigbhour in the other, blue lines indicate that each vertex in one set has at most one non-neighbour in the other, the green line indicates that each vertex in one set has at most two non-neighbours in the other, and no line indicates a cojoin.}
 \label{fig:c6-y1 only t2}
\end{figure}

\begin{figure}[ht]
    \centering
    \includegraphics{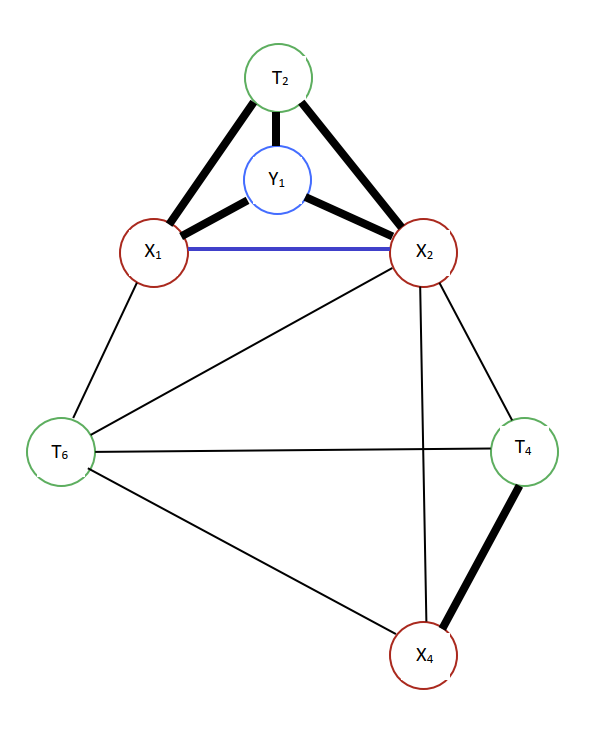}
  \caption {This figure illustrates the the structure of $G \setminus H$ when $G$ contains a $C_6$ and $Y_1$ is the only non-empty $Y_i$ and $T_2, T_4$ and $T_6$ are non-empty. Thick lines indicate a join, thin lines indicate that each vertex in one set has at most one neigbhour in the other, blue lines indicate that each vertex in one set has at most one non-neighbour in the other, and no line indicates a cojoin.}
 \label{fig:c6-y1 3t}
 \end{figure}

This proof is the longest, so we do not try to give more of an overview. However, Figures \ref{fig:c6-y1 only t2} and \ref{fig:c6-y1 3t} give two examples of how the structure simplifies.

\subsection{When the graph contains a $C_5$}

Let $G$ be a $(claw, 4K_1, bridge, C_4$-$twin)$-free graph. Assume $G$ contains an induced $C_5$, call it $H$, and let its vertices be $1,2, \dots, 5$, in that order. Throughout this section, indices are taken mod 5. 

By Theorems \ref{thm-c7} and \ref{thm-c6}, we can assume that $G$ is also $(C_7, C_6)$-free. 

For $i=1,2,3,4,5,$ We consider the sets $X_i, Y_i,$ $Z_i$ and $T_i$ defined above. 
We also need to consider the set $R$ of vertices which have no neighbours on $H$.

We will show that the 21 sets $X_i, Y_i, Z_i$, $T_i$ and $R$ partition the vertices of $G \setminus H$. Similar to the previous cases, under the assumption that each set is big, we show that all the sets $Y_i$ are empty. We prove that each of the sets $X_i, Z_i$, $T_i$ and $R$ is a clique. We show we may assume that $R$ is empty.

\begin{figure}[ht]
    \centering
    \includegraphics{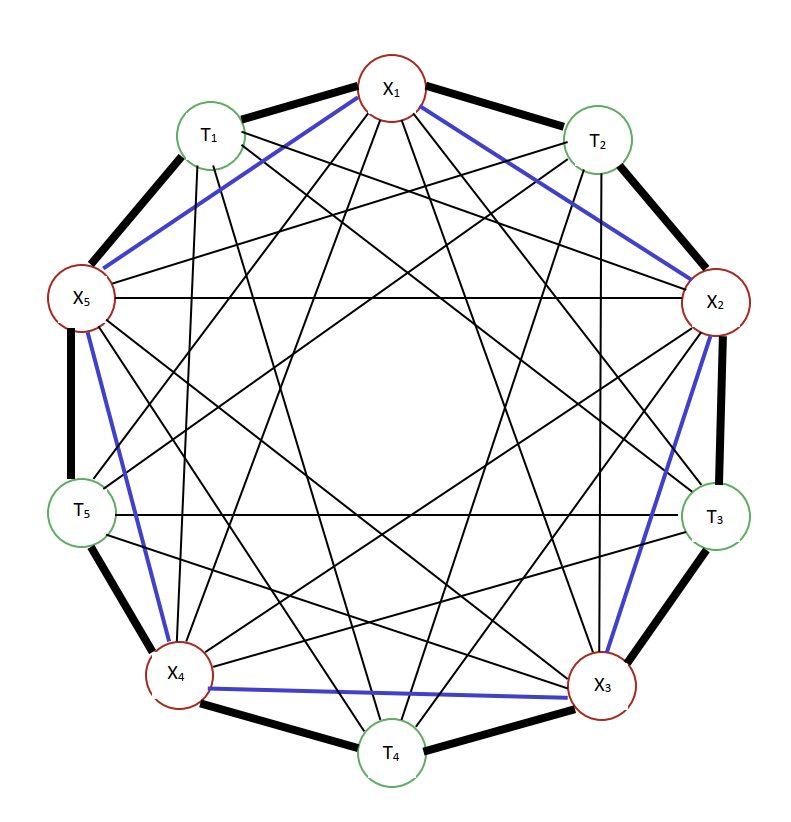}
  \caption{This figure illustrates the the structure of $G \setminus H$ when $G$ contains a $C_5$. Thick lines indicate a join, thin lines indicate that each vertex in one set has at most one neigbhour in the other, blue lines indicate that each vertex in one set has at most one non-neighbour in the other, and no line indicates a cojoin.}
 \label{fig:c5-general}
\end{figure}

We describe the edges between the sets. See Table \ref{tab:C5} for a summary. Again, the edges between sets are usually joins or cojoins, or each vertex in one set has at most one neighbour in another particular set or each vertex in one set has at most one non-neighbour in the other set. Often, one set being non-empty implies that others are empty. See Figure \ref{fig:c5-general} for edges between sets.

We examine the structure of $G$ when $Z \neq \emptyset$,  then $G[T]$, $G[X]$, and $G[X \cup T]$, showing each has bounded clique-width. 

\section{Overview of Sections \ref{sec:c7}, \ref{sec:c6}, and \ref{sec:c5}}

For each case of foundational cycle, we provide a table summarizing the properties that describe the relationships between sets. This serves as a high-level overview of the structure of the graph. The proof of each property follows the same format: We assume for a contradiction that the property does not hold and then obtain one of our forbidden induced subgraphs.

We include only a select handful of proofs of these properties to give the reader an idea of how they proceed. The remaining proofs are omitted in the body but can be found in the Appendix.

Using these properties, we then provide clique-width labelling algorithms using a bounded number of labels for each case.

\begin{table}[h!]
\begin{center}

\begin{tabular}{|c|c|} 
 \hline
 \textbf{Symbol} & \textbf{Definition}\\
 \hline
 $A \;\circled{1}\; B$ & $A$ and $B$ form a join\\
 
 $A \;\circled{$0$}\; B$ & $A$ and $B$ form a co-join\\
 \hline 
 $A \;\circled{$\leq k$}\; B$ & A vertex in $A$ is adjacent to at most $k$ vertices in $B$\\
 
 $A \;\circled{$\leq \bar{k}$}\; B$ & A vertex in $A$ is non-adjacent to at most $k$ vertices in $B$\\
 
 $vv.$ & The property is vice versa.\\
 \hline
\end{tabular}\caption{Summary of Notation Used}
\end{center}
\end{table}

\section{The case of the $C_7$}\label{sec:c7}

Recall the following definitions. 

Let $X_i$ be the set of vertices of $G$ that are adjacent to $i, i+1, i+2$. Let $X$ = $X_1 \cup \dots \cup X_7$.

Let $Y_i$ be the set of vertices of $G$ that are adjacent to $i, i+1, i+2, i+3.$ Let $Y$ = $Y_1 \cup \dots \cup Y_7$.

Let $Z_i$ be the set of vertices of $G$ that are adjacent to $i, i+1, i+3, i+4.$ Let $Z$ = $Z_1 \cup \dots \cup Z_7$.

\noindent Table \ref{tab:C7} summarizes the adjacency between the sets. 

\begin{Observation}\label{obs:c7-sets}
    Any vertex in $G$ must belong to one of $H$, $X_i$, $Y_i$ or $Z_i$.
\end{Observation}	
\noindent {\it Proof}. Consider a vertex $a$ that does not belong to the induced cycle $H$.

If $a$ is not adjacent to any vertex of $H$ or is adjacent to only one vertex of $H$, say $1$, then the set $\{a,2,4,6\}$ induces a $4K_1$.

If $a$ is adjacent to two consecutive vertices of $H$, say $1$ and $2$, then the set $\{a,3,5,7\}$ induces a $4K_1$. Similarly, if $a$ is adjacent to two non-consecutive vertices of $H$, WLOG say $1$ and $3$, then the set $\{a,2,4,6\}$ induces a $4K_1$.

If $a$ is adjacent to four vertices of $H$ such that three of those vertices form a path and one is not adjacent to any of the other three, say $1$, $2$, $3$, $5$, then the set $\{a,1,3,5\}$ induces a $claw$. If $a$ is adjacent to two consecutive vertices and two non-consecutive vertices, say $1,3,5,7$, then the set $\{a,1,3,5\}$ also induces a $claw$.

If $a$ is adjacent to five or more vertices of $H$, then there is a $claw$. The $claw$ exists regardless of the configuration of the five or more vertices to which $a$ is adjacent.

Therefore, if $a$ is not in $H$, it must belong to one of $X_i$, $Y_i$ or $Z_i$.\qed

\begin{Observation}\label{obs:c7-clique}
    The sets $X_i$, $Y_i$, and $Z_i$ are cliques for all $i$.
\end{Observation}	
\noindent {\it Proof}. If $X_i$ contains non-adjacent vertices $a,b$ then $\{i+2, a, b, i+3\}$ is a $claw$. Similarly, if $Y_i$ contains non-adjacent vertices $a,b$ then $\{i+3, i+4, a, b\}$ is a $claw$ and if $Z_i$ contains non-adjacent vertices $a,b$, then $\{i+2, i+5, a,b\}$ is a $4K_1$.\qed

\begin{Observation}\label{obs:c7-z-empty}
    $Z_i$ are all empty.
\end{Observation}
\noindent {\it Proof}. Suppose some $Z_i$ contains two vertices $a,b$. Then the set $\{i,i+1,a,b,i+3,i+4\}$ induces a $bridge$. Then we have $|Z_i | \leq 1$ for all $i$. This contradicts Lemma \ref{lem:big-sets}. Thus, all $Z_i$ must be empty. \qed

\subsection{Properties of the sets $X$ and $Y$}

Properties of the sets $X_i$ and $Y_i$ are given in Table \ref{tab:C7} below.

\begin{table}[ht!]
\begin{center}
\begin{tabular}{ |p{1cm} p{6.5cm} p{6.5cm}| }
\hline
 \textbf{Name} & \textbf{Property} & \textbf{Symmetry} \\
 \hline

\begin{Property} \label{obs:c7-Xi-Xi+1} \vspace*{-1.5em} \end{Property}     & $X_i \;\circled{1}\; X_{i+1}$ &  $X_i \;\circled{1}\; X_{i+6}$ 
 \\

 \begin{Property} \label{obs:c7-Xi-Xi+2} \vspace*{-1.5em} \end{Property}     & $X_i \;\circled{$\leq 1$}\; X_{i+2}$ $ vv.$ & $X_i \;\circled{$\leq 1$}\; X_{i+5}$ $ vv.$ \\
 
 \begin{Property} \label{obs:c7-Xi-Xi+3} \vspace*{-1.5em} \end{Property}     & $X_i \;\circled{0}\; X_{i+3}$ & $X_i \;\circled{0}\; X_{i+4}$ \\

 \begin{Property} \label{obs:c7-xi-xi+2-xi+5} \vspace*{-1.5em} \end{Property}     & $X_i \;\circled{$\leq 1$}\; X_{i+2} \cup X_{i+5}$ &  \\
 
 \begin{Property} \label{obs:c7-no3cons-x} \vspace*{-1.5em} \end{Property}     & There are no three consecutive non-empty sets $X_i$ & \\
 
 \hline
 
 \begin{Property} \label{obs:c7-yi-yi+1} \vspace*{-1.5em} \end{Property}      & $Y_i \;\circled{$\leq \bar{1}$}\; Y_{i+1}$ $vv.$ & $Y_i \;\circled{$\leq \bar{1}$}\; Y_{i+6}$ $vv.$ \\

 \begin{Property} \label{obs:c7-yi-yi+2} \vspace*{-1.5em} \end{Property}       & If $Y_i \not= \emptyset$, then $Y_{i+2} = \emptyset$ & If $Y_i \not= \emptyset$, then $Y_{i+5} = \emptyset$ \\

 \begin{Property} \label{obs:c7-yi-yi+3} \vspace*{-1.5em} \end{Property}       & $Y_i \;\circled{0}\; Y_{i+3}$ & $Y_i \;\circled{0}\; Y_{i+4}$ \\
 
 \hline
 
 \begin{Property} \label{obs:c7-yi-xi} \vspace*{-1.5em} \end{Property}       & $Y_i \;\circled{1}\; X_{i}$ &  $Y_i \;\circled{1}\; X_{i+1}$ \\

 \begin{Property} \label{obs:c7-yi-xi+2} \vspace*{-1.5em} \end{Property}       & If $Y_i \not= \emptyset$, then $X_{i+2} = \emptyset$ & If $Y_i \not= \emptyset$, then $X_{i+6} = \emptyset$ \\

 \begin{Property} \label{obs:c7-yi-xi+3} \vspace*{-1.5em} \end{Property} & $Y_i \;\circled{0}\; X_{i+3}$ & $Y_i \;\circled{0}\; X_{i+5}$ \\

 \begin{Property} \label{obs:c7-yi-xi+4} \vspace*{-1.5em} \end{Property} & $Y_i \;\circled{0}\; X_{i+4}$ & \\
 \hline
\end{tabular}\caption{Properties of $(claw, 4K_1, bridge)$-free graph which contain a $C_7$}
\label{tab:C7}
\end{center}
\end{table}

\subsection{Proofs of some properties}
In the section, we give proofs for some properties in Table \ref{tab:C7}. The other properties have similar proofs and those will be given in the Appendix.

\noindent \textbf{P\ref{obs:c7-Xi-Xi+1}} \textit{$X_i \;\circled{1}\; X_{i+1}$  for all $i$. By symmetry, $X_i \;\circled{1}\; X_{i+6}$.}

\noindent {\it Proof}. If there is a vertex $a \in X_1$ that is non-adjacent to a vertex $b \in X_2$, then the set $\{a,b,5,7\}$ induces a $4K_1$.\qed

\noindent \textbf{P\ref{obs:c7-Xi-Xi+2}} \textit{A vertex in  $X_i$ is adjacent to at most one vertex in $X_{i+2}$ and vice versa. By symmetry, a vertex in  $X_i$ is adjacent to at most one vertex in $X_{i+5}$ and vice versa.}

\noindent {\it Proof}. If there is a vertex $a \in X_1$ that is adjacent to two vertices $b,c \in X_3$, then we claim that every other vertex $d \in X_1$ must be adjacent to exactly one of $\{b,c\}$. Otherwise if $d$ is non-adjacent to both then the set $\{2,d,a,3,b,c\}$ induces a $bridge$ or if $d$ is adjacent to both then the set $\{1,2,a,d,b,c\}$ induces a $bridge$. Without loss of generality, assume $d$ is adjacent to $b$. Then by the assumption that $X_1$ is big, there is another vertex $e \in X_1$ that is also adjacent to $b$. Then, by symmetry, every other vertex in $X_3$ must be adjacent to exactly one of $\{d,e\}$. But $c$ is adjacent to neither. Contradiction.\qed

\noindent \textbf{P\ref{obs:c7-Xi-Xi+3}} \textit{ $X_i \;\circled{0}\; X_{i+3}$  for all $i$. By symmetry, $X_i \;\circled{0}\; X_{i+4}$.}

\noindent {\it Proof}. If there is a vertex $a \in X_1$ that is adjacent to a vertex $b \in X_4$, then the set $\{a,b,1,3\}$ induces a $claw$.\qed

\noindent \textbf{P\ref{obs:c7-no3cons-x}} \textit{There are no three consecutive non-empty sets $X_i$.}

\noindent {\it Proof}. By P\ref{obs:c7-Xi-Xi+1}, any vertex $a \in X_1$ must be adjacent to any vertex $b \in X_2$. Then by P\ref{obs:c7-Xi-Xi+2} and the assumption that $X_7$ is big, we may choose a vertex $c \in X_7$ that is non-adjacent to $b$. By P\ref{obs:c7-Xi-Xi+1}, $a$ is adjacent to $c$. Then the set $\{c,1,a,2,b,3\}$ induces a $bridge$.\qed

\noindent \textbf{P\ref{obs:c7-yi-yi+1}} \textit{A vertex in  $Y_i$ is non-adjacent to at most one vertex in $Y_{i+2}$ and vice versa. By symmetry, a vertex in  $Y_i$ is non-adjacent to at most one vertex in $Y_{i+5}$ and vice versa.}

\noindent {\it Proof}. If there is a vertex $a \in Y_1$ that is non-adjacent to two vertices $b,c \in Y_3$, then we claim that every other vertex $d \in Y_1$ must be adjacent to exactly one of $\{b,c\}$. Otherwise if $d$ is non-adjacent to both then the set $\{d,a,3,4,c,b\}$ induces a $bridge$ or if $d$ is adjacent to both then the set $\{a,2,d,3,c,b\}$ induces a $bridge$. Without loss of generality, assume $d$ is non-adjacent to $b$ and adjacent to $c$. Then by the assumption that $Y_1$ is big, there is another vertex $e \in Y_1$ that is also non-adjacent to $b$ and adjacent to $c$. Then, since $b$ is non-adjacent to both $d$ and $e$, by symmetry, every other vertex in $Y_3$ must be adjacent to exactly one of $\{d,e\}$. But $c$ is adjacent to both. Contradiction.\qed

\noindent \textbf{P\ref{obs:c7-yi-yi+2}} \textit{If $Y_i \not= \emptyset$, then $Y_{i+2} = \emptyset$. By symmetry, if $Y_i \not= \emptyset$, then $Y_{i+5} = \emptyset$.}

\noindent {\it Proof}. Suppose $Y_1 \not= \emptyset$ and $Y_{i+2} \not= \emptyset$. By P\ref{obs:c7-yi-yi+2}a and the assumption that $Y_1$ is big, there is a vertex $a \in Y_1$ that is adjacent to two vertices $c,d \in Y_3$. By P\ref{obs:c7-yi-yi+2}a and the assumption that $Y_3$ is big, there is a vertex $b$ in $Y_1$, different from $a$, that is adjacent to both $c$ and $d$. But now the set $\{1,2,a,b,c,d\}$ induces a $bridge$.\qed

\subsection{Clique-width Labelling}

\subsubsection{The Graph $G[X]$}
\begin{Lemma}\label{lem:c7-xi-xi+1-bcw}
    If $X_i$ and $X_{i+1}$ are non-empty, the graph $G[X]$ has bounded clique-width.
\end{Lemma}
\noindent {\it Proof}. Assume $X_1$ and $X_2$ are non-empty. Then by P\ref{obs:c7-no3cons-x}, $X_3$ and $X_7$ are empty. By P\ref{obs:c7-Xi-Xi+1}, $X_1$ and $X_2$ are a clique.

Suppose only $X_1$ and $X_2$ are non-empty. Then the graph consists of the $C_7$ and the clique $X_1 \cup X_2$, so labelling is trivial.

Now consider which other sets may be non-empty. We may have $X_4$, $X_5$, and $X_6$ non-empty.

Suppose only $X_4$ is non-empty. By P\ref{obs:c7-Xi-Xi+3}, $X_4 \;\circled{0}\; X_1$. By P\ref{obs:c7-Xi-Xi+2}, a vertex in $X_2$ is adjacent to at most one vertex in $X_4$ and vice versa. So we can label $X_2 \cup X_4$ via pairs, label $X_1$, then take the disjoint union of $X_2 \cup X_4$ and $X_1$ and make the join between vertices with label (x1, old) and (x2, old). Thus the graph has bounded clique-width.

Suppose only $X_5$ is non-empty. By P\ref{obs:c7-Xi-Xi+3}, $X_5 \;\circled{0}\; X_1 \cup X_2$. So the graph has bounded clique-width.

The case when only $X_6$ is non-empty is symmetric to the case when only $X_4$ is non-empty.

By P\ref{obs:c7-no3cons-x}, only two of $X_4$, $X_5$, and $X_6$ may be non-empty at a time.

Suppose $X_4$ and $X_5$ are non-empty. By P\ref{obs:c7-Xi-Xi+1}, $X_4 \;\circled{1}\; X_5$. Take the disjoint union of labelled graphs $X_1 \cup X_2 \cup X_4$ and $X_5$, then make the join between vertices with labels (x4, old) and (x5, old). Thus the graph has bounded clique-width.

The case when $X_5$ and $X_6$ are non-empty is symmetric to the case when $X_4$ and $X_5$ are non-empty.

Suppose $X_4$ and $X_6$ are non-empty. By P\ref{obs:c7-Xi-Xi+2}, a vertex in $X_4$ is adjacent to at most one vertex in $X_6$ and vice versa. Since a vertex in $X_1$ is adjacent to at most one in $X_6$ and vice versa and a vertex in $X_2$ is adjacent to at most one in $X_4$ and vice versa, we can label via rows. Thus this graph also has bounded clique-width.

Therefore, if $X_1$ and $X_2$ are non-empty, the graph $G[X]$ has bounded clique-width.\qed

\begin{Lemma}\label{lem:c7-xi-xi+2-bcw}
    If $X_i$ and $X_{i+2}$ are non-empty, the graph $G[X]$ has bounded clique-width.
\end{Lemma}
\noindent {\it Proof}. Assume $X_1$ and $X_3$ are non-empty. By Lemma \ref{lem:c7-xi-xi+1-bcw}, we can assume no two consecutive sets $X_i$ and $X_{i+1}$ are both non-empty. Thus $X_2$, $X_4$, and $X_7$ are all empty. 

Suppose only $X_1$ and $X_3$ are non-empty. By P\ref{obs:c7-Xi-Xi+2}, a vertex in $X_1$ is adjacent to at most one vertex in $X_3$ and vice versa. Then we can label via pairs.

Now consider which other sets may be non-empty. We may have $X_5$ or $X_6$ non-empty, but not both since otherwise there would be two consecutive non-empty sets. $X_5$ and $X_6$ are symmetric, so we will only consider $X_5$.

Suppose $X_5$ is non-empty. By P\ref{obs:c7-Xi-Xi+2}, a vertex in $X_3$ is adjacent to at most one vertex in $X_5$. By P\ref{obs:c7-Xi-Xi+3}, $X_1 \;\circled{0}\; X_5$. Since a vertex in $X_1$ is adjacent to at most one vertex in $X_3$ and vice versa, we can label via rows.

Therefore, if $X_1$ and $X_3$ are non-empty, the graph has bounded clique-width.\qed

\begin{Lemma}\label{lem:c7-xi-xi+3-bcw}
    If $X_i$ and $X_{i+3}$ are non-empty, the graph $G[X]$ has bounded clique-width.
\end{Lemma}
\noindent {\it Proof}. Assume $X_1$ and $X_4$ are non-empty. By Lemma \ref{lem:c7-xi-xi+1-bcw}, we can assume no two consecutive sets $X_i$ and $X_{i+1}$ are both non-empty. Thus all other sets $X_i$ are empty. By P\ref{obs:c7-Xi-Xi+3}, $X_1 \;\circled{0}\; X_4$. Then $G[X]$ consists of two disjoint cliques. Therefore it has bounded clique-width.\qed

\begin{Lemma}\label{lem:c7-GX}
    The graph $G[X]$ has bounded clique-width.
\end{Lemma}
\noindent {\it Proof}. Since we can take any non-empty set $X_i$ to be $X_1$, without loss of generality assume $X_1$ is non-empty. Then if $X_2$ or $X_7$ is non-empty, refer to Lemma \ref{lem:c7-xi-xi+1-bcw} and the graph has bounded clique-width. If $X_3$ or $X_6$ is non-empty, refer to Lemma \ref{lem:c7-xi-xi+2-bcw} and the graph has bounded clique-width. If $X_4$ or $X_5$ is non-empty, refer to Lemma \ref{lem:c7-xi-xi+3-bcw} and the graph has bounded clique-width. Therefore, $G[X]$ has bounded clique-width.\qed

\subsubsection{The Graph $G[Y]$}
\begin{Lemma}\label{lem:c7-GY}
The graph $G[Y]$ has bounded clique-width.
\end{Lemma}
\noindent {\it Proof}. Since we can take any non-empty set $Y_i$ to be $Y_1$, without loss of generality assume $Y_1$ is non-empty. Then by P\ref{obs:c7-yi-yi+2}, $Y_3$ and $Y_6$ are empty.

Suppose $Y_2$ is non-empty. Then by P\ref{obs:c7-yi-yi+2}, $Y_4$ and $Y_7$ are empty. So only $Y_5$ may also be non-empty. By P\ref{obs:c7-yi-yi+1}, a vertex in $Y_1$ is non-adjacent to at most one vertex in $Y_2$ and vice versa. So we can label the graph $Y_1 \cup Y_2$ via pairs.

Suppose $Y_5$ is also non-empty. By P\ref{obs:c7-yi-yi+3}, $Y_1 \;\circled{0}\; Y_5$ and $Y_2 \;\circled{0}\; Y_5$. So we can take the disjoint union of the labelled graphs $Y_5$ and $Y_1 \cup Y_2$. This holds whether any of the three sets $Y_1, Y_2, Y_5$ are empty.

All other combinations of sets $Y_i$ are symmetric to when $Y_1$ and $Y_2$ are non-empty. That is, we can only have two consecutive non-empty sets $Y_i$ and $Y_{i+1}$ and one non-consecutive set $Y_{i+4}$.

Therefore the graph $G[Y]$ has bounded clique-width.\qed

\subsubsection{The Graph $G[X\cup Y]$}
\begin{Lemma}\label{lem:c7-GXY}
The graph $G[X\cup Y]$ has bounded clique-width.
\end{Lemma}

\noindent {\it Proof}. For any set $Y_i$, we have $Y_i \;\circled{1}\; X_i$ and $Y_i \;\circled{1}\; X_i+1$. For all other $X_j$, either $Y_i \;\circled{0}\; X_j$ or $X_j$ must be empty if $Y_i$ is non-empty. Then to construct $G[X\cup Y]$, first take the disjoint union of $G[X]$ and $G[Y]$. For all $Y_i$, make the join between labels (yi, old) and (xi, old) and between labels (yi, old) and (xi+1, old).

Therefore, the graph $G[X \cup Y]$ has bounded clique-width.\qed

\subsection{Graphs containing a $C_7$}
\begin{Theorem}\label{thm-c7}
    Let $G$ be a $(claw, 4K_1, bridge)$-free graph. If $G$ contains a $C_7$, then $G$ has bounded clique-width.
\end{Theorem}

\noindent {\it Proof}. In Lemma \ref{lem:c7-GX}, we proved that if $G$ consists of only the $C_7$ and the sets $X_i$, then $G$ has bounded clique-width. In Lemma \ref{lem:c7-GY}, we proved that if $G$ consists of only the $C_7$ and the sets $Y_i$, then $G$ has bounded clique-width. In Lemma \ref{lem:c7-GXY}, we proved that if $G$ consists of the $C_7$ and the sets $X_i$ and $Y_i$, then $G$ has bounded clique-width. From Observations \ref{obs:c7-sets} and \ref{obs:c7-z-empty}, we know no other configuration of vertices is possible. Therefore, if $G$ contains a $C_7$, $G$ has bounded clique-width. \qed

\section{The case of the $C_6$}\label{sec:c6}
Recall the definitions of the sets of neighbours of the $C_6$.
For $i=1,2,3,4,5,6,$, we consider the sets $X_i, Y_i,$ and $Z_i$ defined as in the case of the $C_7$, 
and let $X,Y,$ and $Z$ be the union of all $X_i, Y_i,$ and $Z_i$, respectively. 
Let $T_i$ be the set of vertices of $G$ that are adjacent to $i, i+1$. 
Let $T$ = $T_1 \cup \dots \cup T_6$.

\begin{Observation}\label{obs:c6-sets}
    Any vertex in $G$ must belong to one of $H$, $T_i$, $X_i$, $Y_i$ or $Z_i$.
\end{Observation}	
\noindent {\it Proof}. Consider a vertex $a$ that does not belong to the induced cycle $H$.

If $a$ is not adjacent to any vertex of $H$ or is adjacent to only one vertex of $H$, say $1$, then the set $\{a,2,4,6\}$ induces a $4K_1$.

If $a$ is adjacent to two non-consecutive vertices of $H$, say $1$ and $3$, then the set $\{a,2,4,6\}$ induces a $4K_1$.

If $a$ is adjacent to three non-consecutive vertices of $H$, say $1$, $3$ and $5$, then the set $\{a,2,4,6\}$ induces a $4K_1$. If $a$ is adjacent to two consecutive vertices and one non-consecutive vertex of $H$, say $1$, $2$, $4$, then the set $\{4,a,3,5\}$ induces a $claw$.

If $a$ is adjacent to four vertices of $H$ such that three of those vertices form a path and one is not adjacent to any of the other three, say $1$, $2$, $3$, $5$, then the set $\{a,1,3,5\}$ induces a $claw$.

If $a$ is adjacent to five or more vertices of $H$, say $a$ is adjacent to all vertices of $H$ except possibly $6$, then the set $\{a,1,3,5\}$ induces a $claw$.

Therefore, if $a$ is not in $H$, it must belong to one of $T_i$, $X_i$, $Y_i$ or $Z_i$.\qed

\begin{Observation}\label{obs:c6-clique}
    The sets $T_i$, $X_i$, $Y_i$, and $Z_i$ are cliques for all $i$.
\end{Observation}	
\noindent {\it Proof}. If $T_i$ contains non-adjacent vertices $a,b$ then $\{i+1, a, b, i+2\}$ is a $claw$. If $X_i$ contains non-adjacent vertices $a,b$ then $\{i+2, a, b, i+3\}$ is a $claw$. Similarly, if $Y_i$ contains non-adjacent vertices $a,b$ then $\{i+3, i+4, a, b\}$ is a $claw$ and if $Z_i$ contains non-adjacent vertices $a,b$, then $\{i+2, i+5, a,b\}$ is a $claw$.\qed

\begin{Observation}\label{obs:c6-z-empty}
    $Z_i$ are all empty.
\end{Observation}
\noindent {\it Proof}. Suppose some $Z_i$ contains two vertices $a,b$. Then the set $\{i,i+1,a,b,i+3,i+4\}$ induces a $bridge$. Then we have $|Z_i | \leq 1$ for all $i$. This contradicts Lemma \ref{lem:big-sets}. Thus, all $Z_i$ must be empty.\qed

\subsection{Properties of the sets $T$, $X$, and $Y$}
We summarize the properties of the sets in Table \ref{tab:C6} below. Since it is a routine but tedious matter to verify these properties, we will give the proofs in the Appendix.
\begin{center}
\begin{longtable}[h!]{ |p{1cm} p{6.5cm} p{6.5cm}| }
\hline
 \textbf{Name} & \textbf{Property} & \textbf{Symmetry} \\
 \hline
 
 \begin{Property} \vspace*{1em} \label{obs:c6-ti-ti+1} \vspace*{-1.5em} \end{Property} & If $T_i \not= \emptyset$, then $T_{i+1} = \emptyset$ &  If $T_i \not= \emptyset$, then $T_{i+5} = \emptyset$ \\
 
 \begin{Property}  \label{obs:c6-ti-ti+2} \vspace*{-1.5em} \end{Property} & $T_i \;\circled{$\leq 1$}\; T_{i+2}$ $vv.$ & $T_i \;\circled{$\leq 1$}\; T_{i+4}$ $vv.$ \\
 
 \begin{Property} \label{obs:c6-ti-ti+3} \vspace*{-1.5em} \end{Property} & If $T_i \not= \emptyset$, then $T_{i+3} = \emptyset$ & \\
 
 \begin{Property} \label{obs:c6-ti-ti+2-ti+4} \vspace*{-1.5em} \end{Property} & If a vertex in $T_i$ is adjacent to a vertex in $T_{i+2}$ and a vertex in $T_{i+4}$, then those vertices must also be adjacent & $vv.$ holds for all orderings of $T_i$, $T_{i+2}$, and $T_{i+4}$ in the statement\\
 
 \hline
 
 \begin{Property} \label{obs:c6-Xi-Xi+1} \vspace*{-1.5em} \end{Property} & $X_i \;\circled{$\leq \bar{1}$}\; X_{i+1}$ $vv.$ & $X_i \;\circled{$\leq \bar{1}$}\; X_{i+5}$ $vv.$ \\

 \begin{Property} \label{obs:c6-Xi-Xi+2} \vspace*{-1.5em} \end{Property} & $X_i \;\circled{$\leq 1$}\; X_{i+2}$ $vv.$ & $X_i \;\circled{$\leq 1$}\; X_{i+4}$ $vv.$\\

 \begin{Property} \label{obs:c6-Xi-Xi+3} \vspace*{-1.5em} \end{Property} & $X_i \;\circled{0}\; X_{i+3}$ & \\

 \begin{Property} \label{obs:c6-no3cons-x} \vspace*{-1.5em} \end{Property} & There are no three consecutive non-empty sets $X_i$ &\\
 
 \hline
 
 \begin{Property} \label{obs:c6-yi-yi+1} \vspace*{-1.5em} \end{Property} & $Y_i \;\circled{$\leq \bar{1}$}\; Y_{i+1}$ $vv.$& $Y_i \;\circled{$\leq \bar{1}$}\; Y_{i+5}$ $vv.$\\

 \begin{Property} \label{obs:c6-yi-yi+2} \vspace*{-1.5em} \end{Property} & If $Y_i \not= \emptyset$, then $Y_{i+2} = \emptyset$ & If $Y_i \not= \emptyset$, then $Y_{i+4} = \emptyset$ \\

 \begin{Property} \label{obs:c6-yi-yi+3} \vspace*{-1.5em} \end{Property} & $Y_i \;\circled{$\leq 1$}\; Y_{i+3}$ $vv.$& \\
 
 \hline
 
 \begin{Property} \label{obs:c6-xi-ti} \vspace*{-1.5em} \end{Property} & $X_i \;\circled{1}\; T_i$ & $X_i \;\circled{1}\; T_{i+1}$\\

 \begin{Property} \label{obs:c6-xi-ti+2} \vspace*{-1.5em} \end{Property} & $X_i \;\circled{$\leq 1$}\; T_{i+2}$ $vv.$& $X_i \;\circled{$\leq 1$}\; T_{i+5}$ $vv.$\\

 \begin{Property} \label{obs:c6-xi-ti+3} \vspace*{-1.5em} \end{Property} & $X_i \;\circled{0}\; T_{i+3}$ & $X_i \;\circled{0}\; T_{i+4}$\\

 \begin{Property} \label{obs:c6-xi-ti-ti+2} \vspace*{-1.5em} \end{Property} & If $X_i \not= \emptyset$, then $T_i \;\circled{0}\; T_{i+2}$ & If $X_i \not= \emptyset$, then $T_{i+1} \;\circled{0}\; T_{i+5}$\\

 \begin{Property} \label{obs:c6-xi-xi+1-ti} \vspace*{-1.5em} \end{Property} & If $X_i \not= \emptyset$ and $X_{i+1} \not= \emptyset$, then $T_i = \emptyset$ & If $X_i \not= \emptyset$ and $X_{i+1} \not= \emptyset$, then $T_{i+2} = \emptyset$\\

 \begin{Property} \label{obs:c6-ti-xi-ti+2} \vspace*{-1.5em} \end{Property} & If $T_i \not= \emptyset$, then $X_i \;\circled{0}\; T_{i+2}$ & \\

 \begin{Property} \label{obs:c6-ti-xi+1-xi+2} \vspace*{-1.5em} \end{Property} & If $T_i \not= \emptyset$, then $X_{i+1} \;\circled{1}\; X_{i+2}$ & \\

 \begin{Property} \label{obs:c6-xi-xi+2-ti+2} \vspace*{-1.5em} \end{Property} & A vertex in $X_i$ cannot be adjacent to both a vertex in $T_{i+2}$ and a vertex in $X_{i+2}$ & A vertex in $X_{i+2}$ cannot be adjacent to both a vertex in $T_{i+1}$ and a vertex in $X_i$ \\
 
 \begin{Property} \label{obs:c6-xi-xi+2-ti+4} \vspace*{-1.5em} \end{Property} & A vertex in $X_i$ cannot be adjacent to both a vertex in $X_{i+2}$ and a vertex in $T_{i+4}$ &\\

  \begin{Property} \label{obs:c6-xi+2-xi-ti+4} \vspace*{-1.5em} \end{Property} & A vertex in $X_{i+2}$ cannot be adjacent to both a vertex in $T_{i+4}$ and a vertex in $X_i$ &\\

  \begin{Property} \label{obs:c6-xi-ti+1-ti+5-a} \vspace*{-1.5em} \end{Property} & If $X_i$, $T_{i+1}$, and $T_{i+5}$ are all non-empty, then $T_{i+1} \;\circled{0}\; T_{i+5}$ &\\

 \hline

 \begin{Property} \label{obs:c6-yi-ti} \vspace*{-1.5em} \end{Property} & $Y_i \;\circled{1}\; T_i$ & $Y_i \;\circled{1}\; T_{i+2}$ \\

 \begin{Property} \label{obs:c6-yi-ti+1} \vspace*{-1.5em} \end{Property} & $Y_i \;\circled{$\leq \bar{2}$}\; T_{i+1}$ & \\

 \begin{Property} \label{obs:c6-yi-ti+1-b} \vspace*{-1.5em} \end{Property} & If a vertex in $Y_i$ is non-adjacent to two vertices in $T_{i+1}$, one of those vertices is adjacent to all of $Y_i$ while the other is non-adjacent to all of $Y_i$ & \\

 \begin{Property} \label{obs:c6-yi-ti+3} \vspace*{-1.5em} \end{Property} & $Y_i \;\circled{0}\; T_{i+3}$ & $Y_i \;\circled{0}\; T_{i+5}$ \\

 \begin{Property} \label{obs:c6-yi-ti+4} \vspace*{-1.5em} \end{Property} & $Y_i \;\circled{0}\; T_{i+4}$ & \\

 \begin{Property} \label{obs:c6-yi-ti-ti+2} \vspace*{-1.5em} \end{Property} & If $Y_i \not= \emptyset$, then one of $T_i$ or $T_{i+2}$ is empty &  \\

 \begin{Property} \label{obs:c6-yi-yi+1-ti} \vspace*{-1.5em} \end{Property} & If $Y_i \not= \emptyset$ and $Y_{i+1} \not= \emptyset$, then $T_i = \emptyset$ & If $Y_i \not= \emptyset$ and $Y_{i+1} \not= \emptyset$, then $T_{i+3} = \emptyset$ \\

 \begin{Property} \label{obs:c6-yi-ti+1-ti+3} \vspace*{-1.5em} \end{Property} & If $Y_i \not= \emptyset$ and $T_{i+3} \not= \emptyset$, then $Y_i \;\circled{1}\; T_{i+1}$ & If $Y_i \not= \emptyset$ and $T_{i+5} \not= \emptyset$, then $Y_i \;\circled{1}\; T_{i+1}$ \\

 \begin{Property} \label{obs:c6-yi-yi+1-ti+1} \vspace*{-1.5em} \end{Property} & If $Y_{i+1} \not= \emptyset$, then $Y_i \;\circled{$ \leq 1 $}\; T_{i+1}$ & \\

 \begin{Property} \label{obs:c6-yi-ti+1-ti+5} \vspace*{-1.5em} \end{Property} & If $Y_i \not= \emptyset$, then $T_{i+1} \;\circled{0}\; T_{i+5}$ & If $Y_i \not= \emptyset$, then $T_{i+1} \;\circled{0}\; T_{i+3}$ \\
 
 \hline
 
 \begin{Property} \label{obs:c6-yi-xi} \vspace*{-1.5em} \end{Property} & $Y_i \;\circled{1}\; X_i$ & $Y_i \;\circled{1}\; X_{i+1}$ \\

 \begin{Property} \label{obs:c6-yi-xi+2} \vspace*{-1.5em} \end{Property} & If $Y_i \not= \emptyset$, then $X_{i+2} = \emptyset$ & If $Y_i \not= \emptyset$, then $X_{i+5} = \emptyset$\\
 
 \begin{Property} \label{obs:c6-yi-xi+3} \vspace*{-1.5em} \end{Property} & $Y_i \;\circled{0}\; X_{i+3}$ & $Y_i \;\circled{0}\; X_{i+4}$\\
 
 \hline
 
 \begin{Property} \label{obs:c6-yi-xi+1-ti} \vspace*{-1.5em} \end{Property} & If $Y_i \not= \emptyset$, then one of $X_{i+2}$ or $T_i$ is empty & If $Y_i \not= \emptyset$, then one of $X_i$ or $T_{i+2}$ is empty\\
 
 \hline
 \caption{Properties of $(claw, 4K_1, bridge, C_4$-$twin, C_7)$-free graphs which contain a $C_6$ }\label{tab:C6}\\
\end{longtable}

\end{center}

\subsection{Clique-width Labelling}

\subsubsection{The Graph $G[T]$}
\begin{Lemma}\label{lem:c6-GT}
    The graph $G[T]$ has bounded clique width.
\end{Lemma}
\noindent {\it Proof}. Since we can take any non-empty set $T_i$ to be $T_1$, without loss of generality assume $T_1$ is non-empty. Then by Observations \ref{obs:c6-ti-ti+1} and \ref{obs:c6-ti-ti+3}, $T_2$, $T_4$, and $T_6$ are empty. Thus if $T_1$ is non-empty, only $T_3$ and $T_5$ may be non-empty.

By Observations \ref{obs:c6-ti-ti+2} and \ref{obs:c6-ti-ti+2-ti+4}, the graph $G[T_1 \cup T_3 \cup T_5]$ satisfies the conditions of \ref{lem:cliquewidth-2k}. This holds whether either of $T_3$ and $T_5$ are empty.

Therefore, the graph $G[T]$ has bounded clique-width.\qed

\subsubsection{The Graph $G[X]$}

\begin{Lemma}\label{lem:c6-xi-xi+1-bcw}
    If $X_i$ and $X_{i+1}$ are non-empty, $G[X]$ has bounded clique-width.
\end{Lemma}
\noindent {\it Proof}. Assume $X_1$ and $X_2$ are non-empty. By P\ref{obs:c6-Xi-Xi+1}, a vertex in $X_1$ is non-adjacent to at most one vertex in $X_2$ and vice versa. Then we can label the graph $X_1 \cup X_2$ via pairs, thus it has bounded clique-width.

Now consider which other sets may be non-empty. By P\ref{obs:c6-no3cons-x}, $X_3$ and $X_6$ must be empty. We may have $X_4$ and $X_5$ non-empty.

Suppose $X_4$ is non-empty. By P\ref{obs:c6-Xi-Xi+3}, $X_1 \;\circled{0}\; X_4$. By P\ref{obs:c6-Xi-Xi+1}, a vertex in $X_2$ is adjacent to at most one vertex in $X_4$ and vice versa. Then we can label the graph $X_1 \cup X_2 \cup X_4$ via rows.

The case when $X_5$ is non-empty is symmetric to the case when $X_4$ is non-empty.

Suppose both $X_4$ and $X_5$ are non-empty. Then we can label via rows, thus the graph has bounded clique-width.

Therefore, if $X_1$ and $X_2$ are non-empty, the graph $G[X]$ has bounded clique-width.\qed

\begin{Lemma}\label{lem:c6-xi-xi+2-bcw}
    If $X_i$ and $X_{i+2}$ are not empty, $G[X]$ has bounded clique-width.
\end{Lemma}
\noindent {\it Proof}. Assume $X_1$ and $X_3$ are non-empty. By Lemma \ref{lem:c6-xi-xi+1-bcw}, we can assume no two consecutive sets $X_i$ and $X_{i+1}$ are both non-empty.  Thus $X_2$, $X_4$ and $X_6$ are all empty. Then we may have only $X_5$ non-empty.

Suppose only $X_1$ and $X_3$ are non-empty. By P\ref{obs:c6-Xi-Xi+2}, a vertex in $X_1$ is adjacent to at most one vertex in $X_3$ and vice versa. So we can label $X_1 \cup X_3$ via pairs.

Now assume $X_5$ is also non-empty. By P\ref{obs:c6-Xi-Xi+2}, a vertex in $X_1$ is adjacent to at most one vertex in $X_5$ and vice versa and a vertex in $X_3$ is adjacent to at most one vertex in $X_5$ and vice versa. Suppose there is a vertex $a \in X_1$ that is adjacent to vertices $b \in X_3$ and $c \in X_5$, then $bc$ is an edge. If $bc$ is not an edge, then $\{a, b, c, 2\}$ induces a $claw$. Thus $bc$ is an edge. Now the graph satisfies the hypothesis of Lemma~\ref{lem:cliquewidth-2k}, so it has bounded clique-width.

Therefore, if $X_1$ and $X_3$ are non-empty, the graph $G[X]$ has bounded clique-width.\qed

\begin{Lemma}\label{lem:c6-GX}
    The graph $G[X]$ has bounded clique-width.
\end{Lemma}
\noindent {\it Proof}. Since we can take any non-empty set $X_i$ to be $X_1$, without loss of generality assume $X_1$ is non-empty. Then if $X_2$ or $X_6$ is non-empty, refer to Lemma \ref{lem:c6-xi-xi+1-bcw} and the graph has bounded clique-width. If $X_3$ or $X_5$ is non-empty, refer to Lemma \ref{lem:c6-xi-xi+2-bcw} and the graph has bounded clique-width. So only $X_4$ may be non-empty. By P\ref{obs:c6-Xi-Xi+3}, $X_1 \;\circled{0}\; X_4$, so the graph consists of two independent cliques $X_1$ and $X_4$, for which labelling is trivial.

Therefore $G[X]$ has bounded clique-width.\qed

\subsubsection{The Graph $G[Y]$}

\begin{Lemma}\label{lem:c6-GY}
    The graph G[Y] has bounded clique width.
\end{Lemma}
\noindent {\it Proof}. Since we can take any non-empty set $Y_i$ to be $Y_1$, without loss of generality assume $Y_1$ is non-empty. Then by P\ref{obs:c6-yi-yi+2}, $Y_3$ and $Y_5$ are empty. So we may have only $Y_2$, $Y_4$, and $Y_6$ non-empty.

Suppose $Y_2$ is non-empty. Then by P\ref{obs:c6-yi-yi+2}, $Y_4$ and $Y_6$ are empty. By P\ref{obs:c6-yi-yi+1}, a vertex in $Y_1$ is non-adjacent to at most one vertex in $Y_2$ and vice versa. So we can label the graph $Y_1 \cup Y_2$ via pairs.

The case when $Y_6$ is non-empty is symmetric to the case when $Y_2$ is non-empty.

Suppose $Y_4$ is non-empty. By P\ref{obs:c6-yi-yi+3}, a vertex in $Y_1$ is adjacent to at most one vertex in $Y_4$ and vice versa. Then the graph satisfies the conditions of Lemma \ref{lem:cliquewidth-2k} and thus has bounded clique-width.

Therefore, the graph $G[Y]$ has bounded clique-width.\qed

\subsubsection{The Graph $G[X \cup T]$}

\begin{Lemma}\label{lem:c6-xi-xi+1-t-bcw}
    If $X_i$ and $X_{i+1}$ are non-empty, the graph G$[X \cup T]$ has bounded clique-width.
\end{Lemma}
\noindent {\it Proof}. Assume $X_1$ and $X_2$ are non-empty. Then by P\ref{obs:c6-xi-xi+1-ti}, $T_1$ and $T_3$ must be empty. By P\ref{obs:c6-no3cons-x}, $X_3$ and $X_6$ must also be empty.

Suppose $X_4$ and $X_5$ are non-empty. Then by P\ref{obs:c6-xi-xi+1-ti}, $T_4$ and $T_6$ are empty. By P\ref{obs:c6-ti-ti+3}, only one of $T_2$ or $T_5$ may be non-empty. The cases are symmetric so without loss of generality assume $T_2$ is non-empty. By Observations \ref{obs:c6-xi-ti} and \ref{obs:c6-xi-ti+3}, $T_2$ forms either a join or a co-join with each of $X_1, X_2, X_4, X_5$. So we can label the graph $G[X]$, label $T_2$, and then make the join between vertices with labels (t2, old) and (x1, old) and between vertices with label (t2, old) and (x2, old). Thus the graph has bounded clique-width.

Now suppose one of $X_4$ and $X_5$ is empty. The cases are symmetric, so without loss of generality assume $X_4$ is non-empty. We know $T_1$ and $T_3$ are empty. If $T_2$ is non-empty, then by P\ref{obs:c6-ti-ti+3}, $T_5$ is empty. If $T_4$ is non-empty, then by P\ref{obs:c6-ti-ti+1}, $T_5$ is empty. Thus if either $T_2$ or $T_4$ is non-empty, $T_5$ is empty.

Suppose $T_2$ is non-empty, so we have the sets $X_1, X_2, X_4, T_2$ non-empty. If $T_2$ is the only set $T_i$ which is non-empty, it is a join or co-join with each of $X_1, X_2$ and $X_4$. So suppose $T_4$ is also non-empty. By Observations \ref{obs:c6-xi-ti}, \ref{obs:c6-xi-ti+3}, and \ref{obs:c6-ti-xi-ti+2}, $T_2$ and $T_4$ are either a join or a co-join with each of $X_1, X_2, X_4$. So we can take the disjoint union of graphs $G[X]$ and $G[T]$ and then make the joins and co-joins between the labelled sets $X_i$ and $T_i$ as needed.

Now suppose $T_6$ is also non-empty such that we have $X_1, X_2, X_4, T_2, T_4, T_6$ non-empty, then by P\ref{obs:c6-ti-xi-ti+2},  $X_4 \;\circled{0}\; T_6$. By P\ref{obs:c6-xi-ti-ti+2}, $T_2 \;\circled{0}\; T_4$ and $T_4 \;\circled{0}\; T_6$. By P\ref{obs:c6-xi-ti+1-ti+5-a}, $T_2 \;\circled{0}\; T_6$. Then consider a vertex $a \in X_1$, $b \in X_2$, and $c \in T_6$. By Observations \ref{obs:c6-Xi-Xi+2} and \ref{obs:c6-xi-ti+2} and the assumption that sets are big, a choice of $d \in X_4$ exists such that $d$ is non-adjacent to both $a$ and $b$. By P\ref{obs:c6-Xi-Xi+3}, $a$ and $d$ are also non-adjacent. Then if $a$ and $b$ are not adjacent, the set $\{a,b,c,d\}$ induces a $4K_1$. So $X_1 \;\circled{1}\; X_2$. Then the only relations that are not a join or co-join are between $X_1$ and $T_6$ and between $X_2$ and $X_4$. Both of these relations are that a vertex in one set is adjacent to at most one vertex in the other set and vice versa. Then we can label the graphs $X_1 \cup T_6$ and $X_2 \cup X_4$ via pairs first, then take the necessary joins with the other sets $T_i$ and $X_i$.

Now suppose $T_2$ is empty, but $X_4$ is still non-empty. So we have the sets $X_1, X_2, X_4$ non-empty and only $T_4$, $T_5$ or $T_6$ may also be non-empty. By P\ref{obs:c6-ti-ti+1}, only one of $T_4$ and $T_5$ may be non-empty at a time. The cases are symmetric, so without loss of generality assume $T_4$ is non-empty. So we have the sets $X_1, X_2, X_4, T_4$ non-empty. Then the only relations that are not a join or a co-join are $X_1 \cup X_2$, $X_2 \cup X_4$, and $X_2 \cup T_4$. By P\ref{obs:c6-xi-xi+2-ti+2}, a vertex in $X_2$ cannot be adjacent to a vertex in $T_4$ and a vertex in $X_4$. So a vertex in $X_2$ may be non-adjacent to at most one in $X_1$ and adjacent to at most one in $X_4 \cup T_4$. We will refer to $X_4 \cup T_4$ as one clique called $X_4^*$. Then we can label the graph $X_1 \cup X_2 \cup X_4^*$, via rows.

Assume now that $T_4$ is empty, but $T_6$ is non-empty. Then we have sets $X_1, X_2, X_4, T_6$ non-empty. By P\ref{obs:c6-xi+2-xi-ti+4}, a vertex in $X_4$ is adjacent to at most one vertex $X_2 \cup T_6$. By P\ref{obs:c6-ti-xi+1-xi+2}, $X_1 \;\circled{1}\; X_2$. Then we have that a vertex in $X_1$ is adjacent to at most one vertex in $T_6$ and vice versa (by Observations \ref{obs:c6-xi-ti+2} and \ref{obs:c6-xi-ti+2}), a vertex in $T_6$ is adjacent to at most one vertex in $X_4$ and vice versa (by Observations \ref{obs:c6-xi-ti+2} and \ref{obs:c6-xi-ti+2}), and a vertex in $X_4$ is adjacent to at most one vertex in $X_2$ and vice versa (by P\ref{obs:c6-Xi-Xi+2}). Any other relations are a join or a co-join. Thus we can label via rows and the graph has bounded clique-width.

Suppose now that both $T_4$ and $T_6$ are non-empty. Then we have sets $X_1, X_2, X_4, T_4, T_6$ non-empty. By P\ref{obs:c6-xi-ti-ti+2}, $T_4 \;\circled{0}\; T_6$. We will refer to the clique $X_4 \cup T_4$ as $X_4^*$. Then we have that a vertex in $X_1$ is adjacent to at most one vertex in $T_6$ and vice versa (by Observations \ref{obs:c6-xi-ti+2} and \ref{obs:c6-xi-ti+2}), a vertex in $T_6$ is adjacent to at most one vertex in $X_4$ and vice versa (by Observations \ref{obs:c6-xi-ti+2} and \ref{obs:c6-xi-ti+2}), a vertex in $X_4$ is adjacent to at most one vertex in $X_2$ and vice versa (by P\ref{obs:c6-Xi-Xi+2}), and a vertex in $X_2$ is adjacent to at most one vertex in $T_4$ and vice versa (by Observations \ref{obs:c6-xi-ti+2} and \ref{obs:c6-xi-ti+2}). Any other relations are a join or a co-join. Thus we can label via rows and the graph has bounded clique-width.

Now suppose $X_4$ is empty, so we have only $X_1$ and $X_2$ non-empty and all other sets $X_i$ are empty. We know $T_1$ and $T_3$ remain empty.

If $T_2$ is non-empty, then by P\ref{obs:c6-ti-ti+3}, $T_5$ is empty and, by P\ref{obs:c6-xi-ti}, $T_2$ is a join with both $X_1$ and $X_2$, so the graph has bounded clique-width.

If one of $T_4$ or $T_6$ is non-empty, by Observations \ref{obs:c6-xi-ti+2}, \ref{obs:c6-xi-ti+2} and \ref{obs:c6-xi-ti+3}, a vertex in $T_4$ or $T_6$ is adjacent to at most one vertex in $X_1 \cup X_2$, so we can label via rows. The cases of either $T_4$ or $T_6$ being non-empty are symmetric.

If $T_2$ and $T_4$ are both non-empty, by Observations \ref{obs:c6-xi-ti}, \ref{obs:c6-xi-ti+3}, and \ref{obs:c6-ti-xi-ti+2}, they are each either a join or a co-join with $X_1$ and $X_2$, so the graph has bounded clique-width. This is symmetric to the case when $T_2$ and $T_6$ are non-empty.

If all three sets $T_2, T_4$ and $T_6$ are non-empty, by Observations \ref{obs:c6-xi-ti}, \ref{obs:c6-xi-ti+3}, and \ref{obs:c6-ti-xi-ti+2}, they are each either a join or a co-join with $X_1$ and $X_2$. Then we can take the disjoint union of labelled graphs $G[X]$ and $G[T]$, and then make the join between vertices with labels (xi, old) and (ti, old) as needed. Thus if $X_4$ is empty, the graph has bounded clique-width.

Therefore, if $X_1$ and $X_2$ are non-empty, the graph $G[X \cup T]$ has bounded clique-width.\qed

\begin{Lemma}\label{lem:c6-xi-xi+2-t-bcw}
    If $X_i$ and $X_{i+2}$ are non-empty, the graph G$[X \cup T]$ has bounded clique-width.
\end{Lemma}
\noindent {\it Proof}. Assume $X_1$ and $X_3$ are non-empty. By P\ref{obs:c6-Xi-Xi+2}, a vertex in $X_1$ is adjacent to at most one vertex in $X_3$ and vice versa.

By Lemma \ref{lem:c6-xi-xi+1-t-bcw}, we can assume no two consecutive sets $X_i$ and $X_{i+1}$ are both non-empty. Thus $X_2$, $X_4$, and $X_6$ are all empty. So of the sets $X_i$, if $X_1$ and $X_3$ are non-empty, only $X_5$ may also be non-empty.

First assume $X_5$ is empty so we have only $X_1$ and $X_3$ non-empty.

Suppose $T_1$ is non-empty. Then by Observations \ref{obs:c6-ti-ti+1} and \ref{obs:c6-ti-ti+3}, $T_2$, $T_4$, and $T_6$ are all empty. By P\ref{obs:c6-xi-ti}, $T_1 \;\circled{1}\; X_1$. By P\ref{obs:c6-xi-ti+3}, $T_1 \;\circled{0}\; X_3$. If only $T_1$ is non-empty, take the disjoint union of labelled graphs $G[X]$ and $T_1$, then make the join between vertices with labels (x1, old) and (t1, old). Thus $X_1 \cup X_3 \cup T_1$ has bounded clique-width.

Suppose both $T_1$ and $T_3$ are non-empty. By P\ref{obs:c6-ti-xi-ti+2}, since $T_1$ is not empty, $X_1 \;\circled{0}\; T_3$. By P\ref{obs:c6-xi-ti}, $X_3 \;\circled{1}\; T_3$. So take the disjoint union of labelled graphs $G[X]$ and $G[T]$, then make the join between vertices with labels (x1, old) and (t1, old) and vertices with labels (x3, old) and (t3, old). Thus $X_1 \cup X_3 \cup T_1 \cup T_3$ has bounded clique-width.

Suppose both $T_1$ and $T_5$ are non-empty. Then by P\ref{obs:c6-xi+2-xi-ti+4}, a vertex in $X_3$ is adjacent to at most one vertex in $X_1 \cup T_5$. Label via pairs first, then the rest of the relations between the sets are either a join or a co-join. Thus $X_1 \cup X_3 \cup T_1 \cup T_5$ has bounded clique-width.

Suppose $T_1$, $T_3$, $T_5$ are all non-empty. We can label via pairs for $X_3$ and $X_1 \cup T_5$ first. Then the remaining graph is either a join or a co-join between the sets. Thus $X_1 \cup X_3 \cup T_1 \cup T_3 \cup T_5$ has bounded clique-width.

Now suppose $T_1$ is empty. Note that $T_1$ being empty and $T_4$ being non-empty is symmetric to when $T_1$ is non-empty and $T_4$ is empty. So we may assume $T_4$ is also empty. 

Suppose $T_2$ is non-empty. By P\ref{obs:c6-xi-ti}, $X_1 \;\circled{1}\; T_1$. By P\ref{obs:c6-xi-xi+2-ti+2}, a vertex in $X_3$ is adjacent to at most one vertex in $X_1 \cup T_2$. So label via pairs for $X_3$ and $X_1 \cup T_2$ first. Then the rest of the graph consists of only joins or co-joins between the sets. Thus $X_1 \cup X_3 \cup T_2$ has bounded clique-width.

Suppose $T_2$ is non-empty. By Observations \ref{obs:c6-ti-ti+1} and \ref{obs:c6-ti-ti+3}, $T_3$ and $T_5$ are empty. So only $T_6$ may be non-empty. Suppose $T_2$ and $T_6$ are both non-empty. Since $X_1$ is non-empty, by P\ref{obs:c6-xi-ti-ti+2}, $T_2 \;\circled{0}\; T_6$. By P\ref{obs:c6-xi-xi+2-ti+2}, a vertex in $X_3$ is adjacent to at most one vertex in $X_1 \cup T_2$. Similarly, a vertex in $X_1$ is adjacent to at most one vertex in $X_3 \cup T_6$. So label the adjacent pairs in $X_1$ and $X_3$, $X_1$ and $T_6$, and $X_3$ and $T_2$. Then the rest of the graph is either a join or a co-join between the sets. Thus $X_1 \cup X_3 \cup T_2 \cup T_6$ has bounded clique-width.

Now assume $X_5$, $T_1$, $T_2$ are all empty. The case when $T_3$ is non-empty is symmetric to the case when $T_2$ is non-empty, so we may assume $T_3$ is also empty. The case when $T_4$ is non-empty is symmetric to the case when $T_1$ is non-empty, so we may assume $T_4$ is also empty. Then of the sets $T_i$, only $T_5$ and $T_6$ may be non-empty.

By P\ref{obs:c6-ti-ti+1}, only one of $T_5$ and $T_6$ may be non-empty at a time. The cases are symmetric, so without loss of generality assume $T_5$ is non-empty. We know a vertex in $X_1$ is adjacent to at most one vertex in $X_3 \cup T_6$. So label the adjacent pairs, then make the join or co-join between the sets for the rest of the graph.

Therefore if $X_1$ and $X_3$ are non-empty but $X_5$ is empty, the graph $G[X \cup T]$ has bounded clique-width.

Now assume $X_5$ is non-empty.

Suppose $T_1$ is non-empty. Then by Observations \ref{obs:c6-ti-ti+1} and \ref{obs:c6-ti-ti+3}, $T_2$, $T_4$, and $T_6$ are all empty. So we may have only $T_1, T_3, T_5$ non-empty. If we have only $T_1$ non-empty, label $T_1 \cup X_5$ via pairs. Then the remaining sets are either a join or a co-join with each other and with $T_1 \cup X_5$. Thus $X_1 \cup X_3 \cup X_5 \cup T_1$ has bounded clique-width.

Suppose both $T_1$ and $T_3$ are non-empty. By P\ref{obs:c6-ti-xi-ti+2}, $X_1 \;\circled{0}\; T_3$. So label the graph $G[X \cup T_1]$ as before, take the disjoint union with $T_3$, then make the join between vertices with labels (x3, old) and (t3, old).

The case when both $T_1$ and $T_5$ are non-empty is symmetric to the case when $T_1$ and $T_3$ are non-empty.

Suppose $T_1$, $T_3$, $T_5$ are all non-empty. By P\ref{obs:c6-ti-xi-ti+2}, all relationships between the sets $T_i$ and $X_i$ are either a join or a co-join. Take the disjoint union of $G[X]$ and $G[T]$, then make the join between the labels (ti, old) and (xi, old) for all $i$. Thus $X_1 \cup X_3 \cup X_5 \cup T_1 \cup T_3 \cup T_5$ has bounded clique-width.

The case when any other set $T_i$ is non-empty is symmetric to the case when $T_1$ is non-empty.

Therefore, if $X_1$, $X_3$ and $X_5$ are non-empty, the graph $G[X \cup T]$ has bounded clique-width.

Therefore if $X_1$ and $X_3$ are non-empty, the graph $G[X \cup T]$ has bounded clique-width.\qed

\begin{Lemma}\label{lem:c6-xi-xi+3-t-bcw}
    If $X_i$ and $X_{i+3}$ are non-empty, the graph G$[X \cup T]$ has bounded clique-width.
\end{Lemma}
\noindent {\it Proof}. Assume $X_1$ and $X_4$ are non-empty. By P\ref{obs:c6-Xi-Xi+3}, $X_1 \;\circled{0}\; X_4$.

By Lemma \ref{lem:c6-xi-xi+1-t-bcw}, we can assume no two consecutive sets $X_i$ and $X_{i+1}$ are both non-empty. Thus $X_2$, $X_3$, $X_5$ and $X_6$ are all empty. So of the sets $X_i$, we have only $X_1$ and $X_4$ non-empty.

Suppose $T_1$ is non-empty. Then by Observations \ref{obs:c6-ti-ti+1} and \ref{obs:c6-ti-ti+3}, $T_2$, $T_4$, and $T_6$ are all empty. So only $T_3$ and $T_5$ may be non-empty. If only $T_1$ is non-empty, then by Observations \ref{obs:c6-Xi-Xi+3}, \ref{obs:c6-xi-ti}, and \ref{obs:c6-xi-ti+3}, between any pairs of the three sets $X_1, X_4, T_1$ is a join or a co-join, so the graph $G[X \cup T_1]$ has bounded clique-width.

Suppose only $T_3$ is non-empty. By P\ref{obs:c6-xi-ti+2}, a vertex in $T_3$ is adjacent to at most one vertex in $X_1$ and, by symmetry, at most one in $X_4$. Then the graph $G[X \cup T_3]$ can be labelled via rows and thus has bounded clique-width.

The case when only $T_5$ is non-empty is symmetric to the case when only $T_1$ is non-empty.

Suppose both $T_1$ and $T_3$ are non-empty. By P\ref{obs:c6-ti-xi-ti+2}, $X_1 \;\circled{0}\; T_3$ and, by P\ref{obs:c6-xi-ti-ti+2}, $T_1 \;\circled{0}\; T_3$. Take the disjoint union of labelled graphs $T_1$ and $G[X \cup T_3]$, then make the join between vertices with labels (t1, old) and (x1, old). Thus $G[X \cup T_1 \cup T_3]$ has bounded clique-width.

Suppose both $T_1$ and $T_5$ are non-empty. Take the disjoint union of labelled graphs $G[X]$ and $G[T_1 \cup T_5]$. Then the sets $X_1$ and $X_4$ and the sets $T_1$ and $T_5$ all have either a join or a co-join between them, so labelling is trivial.

Now suppose $T_1$, $T_3$, and $T_5$ are all non-empty. Since $X_1$ and $X_4$ are non-empty, by P\ref{obs:c6-xi-ti-ti+2}, $T_1 \;\circled{0}\; T_3$ and $T_3 \;\circled{0}\; T_5$. To label this graph, take the disjoint union of labelled graphs $G[X \cup T_3]$ and $G[T_1 \cup T_5]$ then make the join between vertices with labels (x1, old) and (t1, old) and vertices with labels (x4, old) and (t5, old).

Now suppose that $T_1$ is empty. Then all cases are symmetric to when $T_1$ is non-empty, with the exception of $T_3$. From the case when only $T_3$ is non-empty, we know how to label $G[X \cup T_3]$.

Therefore, if $X_1$ and $X_4$ are non-empty, $G[X \cup T]$ has bounded clique-width.\qed

\begin{Lemma}\label{lem:c6-GXT}
    The graph $G[X \cup T]$ has bounded clique-width.
\end{Lemma}
\noindent {\it Proof}. Since we can take any non-empty set $X_i$ to be $X_1$, without loss of generality assume $X_1$ is non-empty. Then by Lemmas \ref{lem:c6-xi-xi+1-t-bcw}, \ref{lem:c6-xi-xi+2-t-bcw} and \ref{lem:c6-xi-xi+3-t-bcw}, all other sets $X_i$ must be empty.

Suppose $T_3$ and $T_6$ are both empty. By Observations \ref{obs:c6-xi-ti} and \ref{obs:c6-xi-ti+3}, $X_1$ forms either a join or a co-join with all other sets $T_i$. So take the disjoint union of labelled graphs $G[X]$ and $G[T]$, then make the join between vertices with labels (x1, old) and (t1, old) and vertices with labels (x1, old) and (t2, old), if either $T_1$ or $T_2$ is non-empty.

Now consider $T_3$ and $T_6$. By P\ref{obs:c6-ti-ti+3}, only one of $T_3$ and $T_6$ may be non-empty at a time. The cases are symmetric, so without loss of generality suppose $T_3$ is non-empty. Then by P\ref{obs:c6-ti-ti+1}, $T_2$ and $T_4$ are empty. So only $T_1$ and $T_5$ may be non-empty.

Suppose $T_1$ is non-empty. So we have the sets $X_1, T_1, T_3$ non-empty. By P\ref{obs:c6-xi-ti-ti+2}, since $X_1$ is non-empty, $T_1 \;\circled{0}\; T_3$. By P\ref{obs:c6-ti-xi-ti+2}, since $T_1$ is non-empty, $X_1 \;\circled{0}\; T_3$. To label the graph, take the disjoint union of labelled graphs $G[X]$ and $G[T]$, then make the join between vertices with labels (x1, old) and (t1, old). Thus $G[X_1 \cup T_1 \cup T_3]$ has bounded clique-width.

Now suppose $T_1$ is empty but $T_5$ is non-empty. So we have the sets $X_1, T_3, T_5$ non-empty. By Observations \ref{obs:c6-ti-ti+2} and \ref{obs:c6-xi-ti+2}, a vertex in $T_3$ is adjacent to at most one vertex in $T_5$ and one vertex in $X_1$. Then we can label via rows, thus the graph has bounded clique-width.

Now suppose $T_1$ and $T_5$ are both non-empty. So we have the sets $X_1, T_1, T_3, T_5$ non-empty. Take the disjoint union of labelled graphs $G[X]$ and $G[T]$. Then by Observations \ref{obs:c6-xi-ti}, \ref{obs:c6-xi-ti+3}, and \ref{obs:c6-ti-xi-ti+2}, $X_1$ is either a join or a co-join with the three sets $T_1, T_3, T_5$, so labelling is trivial.

Therefore, the graph $G[X \cup T]$ has bounded clique-width.\qed

\subsubsection{The Graph $G[T \cup Y]$}

\begin{Observation}\label{obs:c6-yi-ti+1-bcw}
    The graph $G[Y_i \cup T_{i+1}]$ has bounded clique-width.
\end{Observation}
\noindent {\it Proof}. Assume $Y_1$ and $T_2$ are non-empty. By P\ref{obs:c6-yi-ti+1}, a vertex in $Y_1$ is non-adjacent to at most two vertices in $T_2$, but not necessarily vice versa. By P\ref{obs:c6-yi-ti+1-b}, if a vertex $a \in Y_1$ is non-adjacent to two vertices $b,c \in T_2$, one of $\{b,c\}$ must be adjacent to all of $Y_1 - a$ while the other must be non-adjacent to all of $Y_1 - a$. Without loss of generality, assume in this case $b$ is adjacent to $Y_1 - a$ while $c$ is non-adjacent to $Y_1 - a$.

Label all vertices in $T_2$ that are adjacent to none in $Y_1$, similar to and including $c$, with the label (t2, old) first. Now every remaining vertex in $Y_1$ is non-adjacent to at most one vertex in $T_2$. So we can label via pairs.\qed

\begin{Observation}\label{obs:c6-yi-yi+1-ti+1-bcw}
    The graph $G[Y_i \cup Y_{i+1} \cup T_{i+1}]$ has bounded clique-width.
\end{Observation}
\noindent {\it Proof}. Assume $Y_1$, $Y_2$, and $T_2$ are non-empty. By P\ref{obs:c6-yi-yi+1-ti+1}, a vertex in $Y_1$ is non-adjacent to at most one vertex in $T_2$. By P\ref{obs:c6-yi-yi+1}, a vertex in $Y_1$ is also non-adjacent to at most one vertex in $Y_2$. However, a vertex in $T_2$ may be non-adjacent to many vertices in $Y_1$. We will label these vertices in $T_2$ first.

Label a vertex in $T_2$ that is non-adjacent to some vertices in $Y_1$ with the label (t2, new). For each of the vertices it is non-adjacent to in $Y_1$: label them (y1, new), label the vertex in $Y_2$ that it is non-adjacent to (y2, new) if it exists, make the join between (t2, new) and (y2, new), between (y1, old) and (y1, new), and between (y2, old) and (y2, new). Then relabel (y2, new) to (y2, old) and relabel (y1, new) to (y1, old). Repeat for each vertex (t2, new) is non-adjacent to. Now make the join between vertices with labels (t2, old) and (t2, new) and relabel (t2, new) to (t2, old). Repeat with the next $T_2$ vertex that is non-adjacent to some vertices in $Y_1$.

Label any remaining $Y_1$ and $Y_2$ non-adjacencies via pairs. Now the rest of the relations between the cliques are either a join or a co-join, for which labelling is trivial.\qed

\begin{Lemma}\label{lem:GY1Y2T}
    If $Y_i$ and $Y_{i+1}$ are non-empty, $G[Y \cup T]$ has bounded clique-width.
\end{Lemma}
\noindent {\it Proof}. Assume $Y_1$ and $Y_2$ are non-empty. Then by P\ref{obs:c6-yi-yi+2}, since $Y_1$ is non-empty, $Y_3$ and $Y_5$ must be empty. Similarly, since $Y_2$ is non-empty, $Y_4$ and $Y_6$ must also be empty. So of the sets $Y_i$, we have only $Y_1$ and $Y_2$ non-empty.

Suppose $T_1$ is non-empty. Since $Y_1$ and $T_1$ are non-empty, by P\ref{obs:c6-yi-ti-ti+2}, $T_3$ is empty. By P\ref{obs:c6-ti-ti+1}, $T_2$ and $T_6$ are empty. By P\ref{obs:c6-ti-ti+3}, $T_4$ is also empty. So only $T_5$ may be non-empty. By P\ref{obs:c6-yi-ti}, $Y_1 \;\circled{1}\; T_1$. By Observations \ref{obs:c6-yi-ti+3} and \ref{obs:c6-yi-ti+4}, all other relations between $Y_1, Y_2, T_1, T_5$ are a co-join, so labelling is trivial. Thus $Y_1 \cup Y_2 \cup T_1 \cup T_5$ has bounded clique-width. This holds as long as $T_1$ is non-empty, regardless of whether $T_5$ is empty, since we have not used that $T_5$ is non-empty in our argument. 

Now suppose $T_1$ is empty but $T_2$ is non-empty. Since $Y_2$ and $T_2$ are non-empty, by P\ref{obs:c6-yi-ti-ti+2}, $T_4$ is empty. By P\ref{obs:c6-ti-ti+1}, $T_3$ is empty. By P\ref{obs:c6-ti-ti+3}, $T_5$ is also empty. So only $T_6$ may be non-empty. Assume it is non-empty since otherwise, by P\ref{obs:c6-yi-yi+1-ti+1-bcw}, the graph has bounded clique-width. By Observations \ref{obs:c6-yi-ti+3} and \ref{obs:c6-yi-ti+4}, $T_6 \;\circled{0}\; Y_1 \cup Y_2$. By P\ref{obs:c6-yi-ti+1-ti+3}, since $T_6$ is non-empty, $Y_1 \;\circled{1}\; T_2$. By P\ref{obs:c6-yi-ti}, $Y_2 \;\circled{1}\; T_2$. So we can take the disjoint union of graphs $G[Y]$ and $G[X]$ and then make the join between vertices with labels (y1, old) and (t2, old) and vertices with labels (y2, old) and (t2, old). Thus $Y_1 \cup Y_2 \cup T_2 \cup T_6$ has bounded clique-width.

Suppose both $T_1$ and $T_2$ are empty. Then the case when $T_3$ is non-empty is symmetric to the case when $T_2$ is non-empty and the case when $T_4$ is non-empty is symmetric to the case when $T_1$ is non-empty.

By P\ref{obs:c6-ti-ti+1}, $T_5$ and $T_6$ may not both be non-empty at a time. The cases are symmetric, so without loss of generality suppose $T_5$ is non-empty. By P\ref{obs:c6-yi-ti+4} and \ref{obs:c6-yi-ti+3}, $T_5 \;\circled{0}\; Y_1 \cup Y_2$, so we need only take the disjoint union of labelled graphs $T_5$ and $G[Y]$ and we are done. Thus $Y_1 \cup Y_2 \cup T_5$ has bounded clique-width.

Thus if $Y_2$ is non-empty, $G[Y \cup T]$ has bounded clique-width.\qed

\begin{Observation}\label{obs:c6-yi-yi+3-ti+1-bcw}
    The graph $G[Y_i \cup Y_{i+3} \cup T_{i+1}]$ has bounded clique-width.
\end{Observation}
\noindent {\it Proof}. Assume $Y_1$, $Y_4$, and $T_2$ are non-empty. By P\ref{obs:c6-yi-ti+4}, $Y_4 \;\circled{0}\; T_2$. By P\ref{obs:c6-yi-ti+1}, a vertex in $Y_1$ is non-adjacent to at most two vertices in $T_2$. By P\ref{obs:c6-yi-yi+3}, a vertex in $Y_1$ is adjacent to at most one vertex in $Y_4$ and vice versa.

Similar to what is done in P\ref{obs:c6-yi-ti+1-bcw}, label the vertices in $T_2$ that are adjacent to none in $Y_1$ first. Now every remaining vertex in $Y_1$ is non-adjacent to at most one vertex in $T_2$ and vice versa. We can label the non-adjacent pairs of $Y_1$ and $Y_2$ and the adjacent pairs of $Y_1$ and $Y_4$ via rows. Then the remaining vertices of the sets are all either a join or a co-join with each other, for which labelling is trivial.

Thus $G[Y_i \cup Y_{i+3} \cup T_{i+1}]$ has bounded clique-width.\qed

\begin{Lemma}\label{lem:GY1Y4T}
    If $Y_i$ and $Y_{i+3}$ are non-empty, $G[Y \cup T]$ has bounded clique-width.
\end{Lemma}
\noindent {\it Proof}. Assume $Y_1$ and $Y_4$ are non-empty. Then by P\ref{obs:c6-yi-yi+2}, $Y_2$, $Y_3$, $Y_5$, and $Y_6$ must all be empty. So of the sets $Y_i$, we have only $Y_1$ and $Y_4$ non-empty. By P\ref{obs:c6-yi-yi+3}, a vertex in $Y_1$ is adjacent to at most one vertex in $Y_4$ and vice versa.

Suppose $T_1$ is non-empty. Since $Y_1$ and $T_1$ are non-empty, by P\ref{obs:c6-yi-ti-ti+2}, $T_3$ is empty. By P\ref{obs:c6-ti-ti+1}, $T_2$ and $T_6$ are empty. So only $T_5$ may also be non-empty. By P\ref{obs:c6-ti-ti+3}, $T_4$ is also empty. If only $T_1$ is non-empty, it is a join with $Y_1$ and a co-join with $Y_4$, for which labelling is trivial.

The case when only $T_5$ is non-empty is symmetric to P\ref{obs:c6-yi-yi+3-ti+1-bcw} and thus the graph has bounded clique-width in this case.

Suppose both $T_1$ and $T_5$ are non-empty. By P\ref{obs:c6-ti-ti+2}, a vertex in $T_1$ is adjacent to at most one vertex in $T_5$ and vice versa. Then we can label via rows, thus $Y_1 \cup Y_4 \cup T_1 \cup T_5$ has bounded clique-width.

Now suppose $T_1$ is empty but $T_2$ is non-empty. By P\ref{obs:c6-yi-yi+3-ti+1-bcw}, if only $T_2$ is non-empty, the graph has bounded clique-width. By P\ref{obs:c6-ti-ti+1}, $T_3$ is empty. By P\ref{obs:c6-ti-ti+3}, $T_5$ is also empty. So only $T_4$ and $T_6$ may also be non-empty.

Suppose $T_2$ and $T_4$ are both non-empty. By P\ref{obs:c6-yi-ti+1-ti+3}, $Y_1 \;\circled{1}\; T_2$. By Observations \ref{obs:c6-yi-ti+3} and \ref{obs:c6-yi-ti}, all other relations between the sets $Y$ and $T$ are either a join or a co-join. So we can take the disjoint union of the graphs $G[T]$ and $G[Y]$, then make the join between vertices with labels (y4, old) and (t4, old) and vertices with labels (y1, old) and (t2, old). Thus $Y_1 \cup Y_4 \cup T_2 \cup T_4$ has bounded clique-width.

Suppose $T_2$ and $T_6$ are both non-empty. Then by P\ref{obs:c6-yi-ti-ti+2}, $T_4$ is empty and the graph is symmetric to $G[Y_1 \cup Y_4 \cup T_1 \cup T_5]$. Thus the graph has bounded clique-width as explained above.

Now suppose both $T_1$ and $T_2$ are empty. Then all other sets $T_i$ are symmetric to the case of either $T_1$ or $T_2$ being non-empty.

Therefore, if $Y_1$ and $Y_4$ are non-empty, the graph $G[Y \cup T]$ has bounded clique-width.\qed

\begin{Lemma}\label{lem:c6-GTY}
    The graph $G[Y \cup T]$ has bounded clique-width.
\end{Lemma}
\noindent {\it Proof}. Since we can take any non-empty set $Y_i$ to be $Y_1$, without loss of generality assume $Y_1$ is non-empty. By Lemmas \ref{lem:GY1Y2T} and \ref{lem:GY1Y4T} and by P\ref{obs:c6-yi-yi+2}, we may assume all other sets $Y_i$ are empty.

Suppose $T_2$ is empty. By Observations \ref{obs:c6-yi-ti}, \ref{obs:c6-yi-ti+3}, and \ref{obs:c6-yi-ti+4}, $Y_1$ is either a join or a co-join with all other sets $T_i$. Thus regardless of which sets $T_i$ are non-empty, we can take the disjoint union of labelled graphs $Y_1$ and $G[T]$, then make join between the sets $Y_1$ and $T_i$ as needed. Thus, if $T_2$ is empty, $G[Y \cup T]$ has bounded clique-width.

Suppose $T_2$ is non-empty. By P\ref{obs:c6-ti-ti+1}, $T_1$ and $T_3$ are empty. By P\ref{obs:c6-ti-ti+3}, $T_5$ is also empty. So only $T_4$ and $T_6$ may be non-empty. By P\ref{obs:c6-yi-ti+1-bcw}, if only $T_2$ is non-empty, the graph has bounded clique-width.

Suppose $T_2$ and $T_4$ are both non-empty. By P\ref{obs:c6-yi-ti+1-ti+5}, $T_2 \;\circled{0}\; T_4$. So we can take the disjoint union of labelled graphs $Y_1 \cup T_2$ and $T_4$, thus the graph has bounded clique-width.

The case when $T_2$ and $T_6$ are both non-empty is symmetric to the case when $T_2$ and $T_4$ are both non-empty.

Suppose $T_2$, $T_4$, and $T_6$ are all non-empty. We know $T_2 \;\circled{0}\; T_4 \cup T_6$ and $Y_1 \;\circled{0}\; T_4 \cup T_6$. Then to label the graph, we need only take the disjoint union of labelled graphs $Y_1 \cup T_2$ and $T_4 \cup T_6$. Thus the graph has bounded clique-width.

Therefore, the graph $G[Y \cup T]$ has bounded clique-width.\qed

\subsubsection{The Graph $G[X \cup Y]$}
\begin{Lemma}\label{lem:c6-GXY}
    The graph $G[X \cup Y]$ has bounded clique-width.
\end{Lemma}
\noindent {\it Proof}. By Observations \ref{obs:c6-yi-xi}, \ref{obs:c6-yi-xi+2}, and \ref{obs:c6-yi-xi+3}, all sets $Y_i$ and $X_i$ are either related by a join, a co-join, or if one is non-empty then the other is empty. Then to label $G[X \cup Y]$, first take the disjoint union of labelled graphs $G[X]$ and $G[Y]$. Then for all $Y_i$ and $X_i$, make the join between labels (yi, old) and (xi, old) as necessary to satisfy P\ref{obs:c6-yi-xi}. Thus the graph $G[X \cup Y]$ has bounded clique-width.\qed

\subsubsection{The Graph $G[T \cup X \cup Y]$}
For simplicity, we will examine the graph in terms of the sets $Y_i$. Since we can take any non-empty set $Y_i$ to be $Y_1$, without loss of generality assume $Y_1$ is non-empty. Then by P\ref{obs:c6-yi-yi+2}, we have three possible combinations of sets $Y_i$:

Only one set $Y_i$ is non-empty. Assume $Y_1$ is non-empty.

Two consecutive sets $Y_i$ and $Y_{i+1}$ are non-empty. Assume $Y_1$ and $Y_2$ are non-empty.

Two non-consecutive sets $Y_i$ and $Y_{i+3}$ are non-empty. Assume $Y_1$ and $Y_4$ are non-empty.

\begin{Lemma}\label{lem:c6-GTXY1}
    If $Y_i$ is non-empty, the graph $G[T \cup X \cup Y]$ has bounded clique-width.
\end{Lemma}
\noindent {\it Proof}. Assume $Y_1$ is non-empty. Then by P\ref{obs:c6-yi-xi+2}, $X_3$ and $X_6$ are empty.

Then the only relationship between $Y_1$ and any of the sets $X_i$ or $T_i$ that is not a join or a co-join is between $Y_1$ and $T_2$. So if $T_2$ is empty, take the disjoint union of labelled graphs $G[Y]$ and $G[T \cup X]$, then make the join between the labelled vertices of the sets as necessary. Thus, if $T_2$ is empty, the graph $G[Y_1 \cup X \cup T]$ has bounded clique-width.

Now suppose $T_2$ is non-empty. Since we know $X_3$ and $X_6$ are empty, by Observations \ref{obs:c6-xi-ti} and \ref{obs:c6-xi-ti+3}, $T_2$ forms either a join or a co-join with all non-empty sets $X_i$. Then when of the sets $T_i$ only $T_2$ is non-empty, we can take the disjoint union of labelled graphs $G[Y_1 \cup T_2]$ and $G[X]$, then make the join between labelled vertices of the sets as necessary. Thus, $G[Y_1 \cup T_2 \cup X]$ has bounded clique-width.

Suppose $T_2$ is not the only non-empty set $T_i$. By P\ref{obs:c6-ti-ti+1}, $T_1$ and $T_3$ are empty. By P\ref{obs:c6-ti-ti+3}, $T_5$ is also empty. So only $T_4$ and $T_6$ may be non-empty.

Suppose $T_2$ and $T_4$ are both non-empty. Then by P\ref{obs:c6-yi-ti+1-ti+3}, $Y_1 \;\circled{1}\; T_2$. Since $Y_1$ is non-empty and by symmetry, by P\ref{obs:c6-yi-ti+1-ti+5}, $T_2 \;\circled{0}\; T_4$. Then since $Y_1$ and $T_2$ form either a join or a co-join with all other sets, we can take the disjoint union of labelled graphs $G[T_4 \cup X]$, $Y_1$, and $T_2$, then make the joins between the labelled vertices of the sets as necessary. Thus $G[Y_1 \cup T_2 \cup T_4 \cup X]$ has bounded clique-width.

The case when $T_2$ and $T_6$ are both non-empty is symmetric to the case when $T_2$ and $T_4$ are non-empty.

Now suppose $T_2$, $T_4$, and $T_6$ are all non-empty. By P\ref{obs:c6-yi-ti+1-ti+5}, we know $T_2 \;\circled{0}\; T_4 \cup T_6$ and, by P\ref{obs:c6-yi-ti+3}, $Y_1 \;\circled{0}\;T_4 \cup T_6$. So we can take the disjoint union of labelled graphs $G[T_4 \cup T_6 \cup X]$, $T_2$ and $Y_1$, then make the joins between the labelled vertices of the sets as necessary. Thus $G[Y_1 \cup T_2 \cup T_4 \cup T_6 \cup X]$ has bounded clique-width.

Therefore, if $Y_1$ is non-empty, the graph $G[T \cup X \cup Y]$ has bounded clique-width.\qed

\begin{Lemma}\label{lem:c6-GTXY1Y2}
    If $Y_i$ and $Y_{i+1}$ are non-empty, the graph $G[T \cup X \cup Y]$ has bounded clique-width.
\end{Lemma}
\noindent {\it Proof}. Assume $Y_1$ and $Y_2$ are non-empty. Then by P\ref{obs:c6-yi-xi+2}, $X_1$, $X_3$, $X_4$ and $X_6$ are all empty. So of the sets $X_i$, we need only consider when $X_2$ and $X_5$ are non-empty. 

Suppose $X_2$ is non-empty. By P\ref{obs:c6-yi-xi}, $Y_1 \;\circled{1}\; X_2$ and $Y_2 \;\circled{1}\; X_2$. Since $Y_1$ and $Y_2$ are non-empty, by P\ref{obs:c6-yi-xi+1-ti}, $T_1$ and $T_4$ are empty. Then by Observations \ref{obs:c6-xi-ti} and \ref{obs:c6-xi-ti+3}, $X_2$ forms either a join or a co-join with all non-empty sets $T_i$. So we can take the disjoint union of labelled graphs $X_2$ and $G[T \cup Y]$, then make the join between $X_2$ and the labelled vertices of the other sets as necessary. Thus, $G[X_2 \cup T \cup Y]$ has bounded clique-width.

Suppose $X_2$ is empty but $X_5$ is non-empty. By P\ref{obs:c6-yi-xi+3}, $Y_1 \;\circled{0}\; X_5$ and $Y_2 \;\circled{0}\; X_5$. By P\ref{obs:c6-yi-yi+1-ti}, $T_1$ and $T_4$ are empty. Then by Observations \ref{obs:c6-xi-ti} and \ref{obs:c6-xi-ti+3}, $X_5$ forms either a join or a co-join with all non-empty sets $T_i$. So we can take the disjoint union of labelled graphs $X_5$ and $G[T \cup Y]$, then make the join between $X_5$ and the labelled vertices of the other sets as necessary. Thus, $G[X_5 \cup T \cup Y]$ has bounded clique-width.

Suppose $X_2$ and $X_5$ are both non-empty. Then by Observations \ref{obs:c6-yi-xi}, \ref{obs:c6-yi-xi+3}, \ref{obs:c6-xi-ti}, and \ref{obs:c6-xi-ti+3}, $X_2$ and $X_5$ are either a join or a co-join with all other non-empty sets. So we can take the disjoint union of labelled graphs $G[X_2 \cup X_5]$ and $G[T \cup Y]$, then make the join between the labelled vertices of the sets as necessary.

Thus, when $Y_1$ and $Y_2$ are non-empty, the graph $G[T \cup X \cup Y]$ has bounded clique-width.\qed

\begin{Lemma}\label{lem:c6-GTXY1Y4}
    If $Y_i$ and $Y_{i+3}$ are non-empty, the graph $G[T \cup X \cup Y]$ has bounded clique-width.
\end{Lemma}
\noindent {\it Proof}. Assume $Y_1$ and $Y_4$ are non-empty. Then by P\ref{obs:c6-yi-xi+2}, $X_3$ and $X_6$ are empty.

Then the only relationships between $Y_1$ and $Y_4$ and the other potentially non-empty sets that are not a join or a co-join are between $Y_1 \cup T_2$ and $Y_4 \cup T_5$. Then if $T_2$ and $T_5$ are both empty, we can take the disjoint union of labelled graphs $Y_1 \cup Y_4$ and $G[T \cup X]$, then make the join between the labelled vertices of the sets as necessary. Thus if $T_2$ and $T_5$ are both empty, the graph has bounded clique-width.

Suppose $T_2$ and $T_5$ are not both empty. By P\ref{obs:c6-ti-ti+3}, only one of $T_2$ and $T_5$ may be non-empty at a time. The cases are symmetric, so without loss of generality suppose $T_2$ is non-empty.

Suppose $T_2$ is the only non-empty set $T_i$. By Observations \ref{obs:c6-xi-ti} and \ref{obs:c6-xi-ti+3}, since $X_3$ and $X_6$ are both empty, it forms a join or a co-join with all non-empty sets $X_i$. By Observations \ref{obs:c6-yi-xi} and \ref{obs:c6-yi-xi+3}, $Y_1$ and $Y_4$ also both form a join or a co-join with all non-empty sets $X_i$. Then we can take the disjoint union of labelled graphs $T_2 \cup Y_1 \cup Y_4$ and $G[X]$, then make the join between the labelled vertices of the sets as necessary. Thus $G[T_2 \cup Y \cup X]$ has bounded clique-width.

Now consider when $T_2$ is not the only non-empty set $T_i$. By Observations \ref{obs:c6-ti-ti+1}, $T_1$ and $T_3$ are empty. So of the sets $T_i$, we may have only $T_4$ and $T_6$ non-empty.

Suppose $T_2$ and $T_4$ are both non-empty. By P\ref{obs:c6-yi-ti+1-ti+5}, $T_4 \;\circled{0}\; T_2$. By P\ref{obs:c6-yi-ti}, $T_4 \;\circled{1}\; Y_4$. By P\ref{obs:c6-yi-ti+3}, $T_4 \;\circled{0}\; Y_1$. We know $T_2$, $Y_1$, and $Y_4$ all form either a join or a co-join with any non-empty sets $X_i$. So we can take the disjoint union of labelled graphs $G[T_2 \cup Y_1 \cup Y_4]$ and $G[T_4 \cup X]$, then make the join between labelled vertices of the sets as necessary. Thus, $G[T_2 \cup T_4 \cup Y \cup X]$ has bounded clique-width.

The case when $T_2$ and $T_6$ are both non-empty is symmetric to the case when $T_2$ and $T_4$ are non-empty.

Suppose $T_2$, $T_4$, and $T_6$ are all non-empty. We know $T_2$, $Y_1$, and $Y_4$ all form either a join or a co-join with any non-empty sets $X_i$. We also know $T_4$ and $T_6$ form either a join or a co-join with all three sets $T_2$, $Y_1$, and $Y_4$. So we can take the disjoint union of labelled graphs $G[T_2 \cup Y_1 \cup Y_4]$ and $G[T_4 \cup T_6 \cup X]$, then make the join between labelled vertices of the sets as necessary. Thus, $G[T_2 \cup T_4 \cup T_6 \cup Y \cup X]$ has bounded clique-width.

Therefore, when $Y_1$ and $Y_4$ are non-empty, the graph $G[T \cup X \cup Y]$ has bounded clique-width.\qed

\subsection{Graphs containing a $C_6$}

The purpose of the previous results was to prove the following theorem:

\begin{Theorem}\label{thm-c6} Let $G$ be a $(claw, 4K_1, bridge, C_4$-$twin$)-free graph. If $G$ contains a $C_6$, then $G$ has bounded clique-width.
\end{Theorem}
\noindent {\it Proof}. 
By Theorem \ref{thm-c7}, we may assume that $G$ contains no $C_7$. In Lemma \ref{lem:c6-GT}, we proved that if $G$ consists of only the $C_6$ and the sets $T_i$, then $G$ has bounded clique-width. In Lemma \ref{lem:c6-GX}, we proved that if $G$ consists of only the $C_6$ and the sets $X_i$, then $G$ has bounded clique-width. In Lemma \ref{lem:c6-GY}, we proved that if $G$ consists of only the $C_6$ and the sets $Y_i$, then $G$ has bounded clique-width.

In Lemma \ref{lem:c6-GXT}, we proved that if $G$ consists of only the $C_6$ and the sets $X_i$ and $T_i$, then $G$ has bounded clique-width. In Lemma \ref{lem:c6-GTY}, we proved that if $G$ consists of only the $C_6$ and the sets $T_i$ and $Y_i$, then $G$ has bounded clique-width. In Lemma \ref{lem:c6-GXY}, we proved that if $G$ consists of only the $C_6$ and the sets $X_i$ and $Y_i$, then $G$ has bounded clique-width. In Lemmas \ref{lem:c6-GTXY1Y2} and \ref{lem:c6-GTXY1Y4}, we proved that if $G$ consists of the $C_6$ and any combination of sets $T_i$, $X_i$, and $Y_i$, then $G$ has bounded clique-width.

From Observations \ref{obs:c6-sets} and \ref{obs:c6-z-empty}, we know no other configuration of vertices is possible. Therefore, if $G$ contains a $C_6$, $G$ has bounded clique-width. \qed

\section{The case of the $C_5$}\label{sec:c5}
By Lemma \ref{lem:big-sets}, we will assume each of the sets $T_i, X_i, Y_i, Z_i, R_i$ is big, that is, they have at least $5$ vertices. WLOG, we assume that if a set is non-empty, it is big.

In addition to the $claw$, $4K_1$, $bridge$, and $C_4$-$twin$, by Theorems \ref{thm-c7} and \ref{thm-c6}, we may also forbid the $C_7$ and $C_6$ as induced subgraphs.

\begin{Observation}\label{obs:c5-sets}
    Any vertex in $G$ must belong to one of $H$, $T_i$, $X_i$, $Y_i$, $Z$, or $R$.
\end{Observation}	
\noindent {\it Proof}. Consider a vertex $a$ that does not belong to the induced cycle $H$.

If $a$ is not adjacent to any vertices of the cycle $H$, it is in $R$.

If $a$ is adjacent to one vertex in the cycle $H$, say vertex $1$, then the set $\{1,5,2,a\}$ induces a $claw$. Similarly if $a$ is adjacent to two non-consecutive vertices in the cycle $H$, say vertices $1$ and $3$, or if $a$ is adjacent to three non-consecutive vertices, say $1$, $3$, and $4$, then the same set $\{1,5,2,a\}$ induces a $claw$. Thus if $a$ is adjacent to two vertices, it must be in $T_i$ and if $a$ is adjacent to three vertices, it must be in $X_i$.

If $a$ is adjacent to four vertices of the cycle $H$, it is in $Y_i$. If $a$ is adjacent to all five vertices of the cycle $H$, it is in $Z$.

Therefore, if $a$ is not in $H$, it must belong to one of $T_i$, $X_i$, $Y_i$, $Z$, or $R$.\qed

\begin{Observation}\label{obs:c5-clique}
    The sets $T_i$, $X_i$, $Y_i$, $Z$, and $R$ are cliques for all $i$.
\end{Observation}	
\noindent {\it Proof}. If $T_i$ contains non-adjacent vertices $a,b$ then $\{i+1, a, b, i+2\}$ is a $claw$. If $X_i$ contains non-adjacent vertices $a,b$ then $\{i+2, a, b, i+3\}$ is a $claw$. Similarly, if $Y_i$ contains non-adjacent vertices $a,b$ then $\{i+3, i+4, a, b\}$ is a $claw$ and if $Z$ contains non-adjacent vertices $a,b$, then $\{i+2, i+5, a,b\}$ is a $claw$. If $R$ contains non-adjacent vertices $a,b$, then $\{i, i+2, a,b\}$ is a $4K_1$.\qed

\begin{Observation}\label{obs:c5-y-empty}
    $Y_i$ are all empty.
\end{Observation}
\noindent {\it Proof}. Suppose some $Y_i$ contains two vertices $a,b$. Then the set $\{i,a,i+3,i+4,b\}$ induces a $C_4$-$twin$. Then we have $|Y_i | \leq 1$ for all $i$. This contradicts Lemma \ref{lem:big-sets}. Thus, all $Y_i$ must be empty.\qed

\subsection{Properties of the set $R$}

\begin{Observation}\label{obs:c5-rt}
    $T_i \;\circled{1}\; R$ for all $i$.
\end{Observation}
\noindent {\it Proof}. If there is a vertex $a \in T_1$ that is non-adjacent to a vertex $b \in R$, then the set $\{a,b,3,5\}$ induces a $4K_1$.\qed

\begin{Observation}\label{obs:c5-r-t}
    If $R$ is non-empty, $T_i$ is empty for all $i$.
\end{Observation}
\noindent {\it Proof}. By the assumption that sets are big, both $T_1$ and $R$ contain at least two vertices $a,b \in T_1$ and $c,d \in R$. By Observation \ref{obs:c5-r-t}, $a,b,c,d$ are all adjacent. Then the set $\{1,2,a,b,c,d\}$ induces a $bridge$.\qed

\begin{Observation}\label{obs:c5-r-x}
    $X_i \;\circled{0}\; R$ for all $i$.
\end{Observation}
\noindent {\it Proof}. If there is a vertex $a \in X_1$ that is adjacent to a vertex $b \in R$, then the set $\{a,b,1,3\}$ induces a $claw$.\qed

\begin{Observation}\label{obs:c5-r-z}
    $Z \;\circled{0}\; R$
\end{Observation}
\noindent {\it Proof}. If there is a vertex $a \in Z$ that is adjacent to a vertex $b \in R$, then the set $\{a,b,1,3\}$ induces a $claw$.\qed

\begin{Observation}\label{obs:c5-r}
    The set $R$ is empty.
\end{Observation}
\noindent {\it Proof}. By Observations \ref{obs:c5-y-empty}, \ref{obs:c5-r-t}, \ref{obs:c5-r-x} and \ref{obs:c5-r-z}, if $R$ is non-empty then the graph $G$ is disconnected. Since $R$ is a clique, it has bounded clique-width. Then when $R$ is removed, if $G$ has bounded clique-width, so does $G \cup R$. Therefore, we need only consider the graph $G$ when $R$ is empty.\qed

\subsection{Properties of the set $Z$}

\begin{Observation}\label{obs:c5-z-t}
    If $Z$ is non-empty, then $T_i$ is empty for all $i$.
\end{Observation}
\noindent {\it Proof}. If there is a vertex $a \in Z$ that is adjacent to a vertex $b \in T_i$, then the set $\{a,b,3,5\}$ induces a claw. So $Z \;\circled{0}\; T_i$. By the assumption that sets are big, there is another vertex $c \in T_i$ and $d \in Z$, then the set $\{a,b,c,d,1,2\}$ induces a $bridge$.\qed

\begin{Observation}\label{obs:c5-z-x}
    If $Z$ is non-empty, then $X_i$ is empty for all $i$.
\end{Observation}
\noindent {\it Proof}. If there is a vertex $a \in X_i$ that is adjacent to two vertices $b, c \in Z$, then the set $\{a,2,4,5,b,c\}$ induces a $bridge$. Then any vertex in $X_i$ is adjacent to at most one vertex in $Z$. Since by Lemma \ref{lem:big-sets} $Z$ is big, $a$ is non-adjacent to at least four vertices in $Z$. Then there will be two vertices $a,b \in X_i$ which are non-adjacent to vertices $c,d \in Z$. Then the set $\{a,b,c,d,1,2\}$ induces a $bridge$.\qed

\begin{Observation}\label{obs:c5-z}
    If $Z$ is non-empty, then $G$ has bounded clique-width.
\end{Observation}
\noindent {\it Proof}. If $Z$ is non-empty, then the graph $G$ consists of the $C_5$ and the set $Z$. Since $Z$ is a clique, labelling is trivial.\qed

\subsection{Properties of the sets $T$ and $X$}
Properties of the sets $T_i$ and $X_i$ are given in Table \ref{tab:C5} below. Since it is a routine but tedious matter to verify these properties, we will give the proofs in the Appendix.
\begin{table}[h!]
\begin{center}
\begin{tabular}{ |p{1cm} p{6.5cm} p{6.5cm}| }
\hline
 \textbf{Name} & \textbf{Property} & \textbf{Symmetry}\\
 \hline
 \begin{Property} \label{obs:c5-ti-ti+1} \vspace*{-1.5em} \end{Property} & $T_i \;\circled{0}\; T_{i+1}$ &  $T_i \;\circled{0}\; T_{i+4}$\\

 \begin{Property} \label{obs:c5-ti-ti+2} \vspace*{-1.5em} \end{Property} & $T_i \;\circled{$\leq 1$}\; T_{i+2}$ $vv.$& $T_i \;\circled{$\leq 1$}\; T_{i+3}$ $vv.$\\
 
 \begin{Property} \label{obs:c5-ti-ti+2-ti+3} \vspace*{-1.5em} \end{Property} & $T_i \;\circled{$\leq 1$}\; T_{i+2}\bigcup T_{i+3}$ & \\

 \begin{Property} \label{obs:c5-no3cons-t} \vspace*{-1.5em} \end{Property} & There are no 3 consecutive non-empty sets $T_i$ & \\
 
 \hline
 \begin{Property} \label{obs:c5-xi-xi+1} \vspace*{-1.5em} \end{Property} & $X_i \;\circled{$\leq \bar{1}$}\; X_{i+1}$ $vv.$& $X_i \;\circled{$\leq \bar{1}$}\; X_{i+4}$ $vv.$\\

 \begin{Property} \label{obs:c5-xi-xi+2} \vspace*{-1.5em} \end{Property} & $X_i \;\circled{$\leq 1$}\; X_{i+2}$ $vv.$& $X_i \;\circled{$\leq 1$}\; X_{i+3}$ $vv.$\\

 \begin{Property} \label{obs:c5-no3cons-x} \vspace*{-1.5em} \end{Property}  & There are no 3 consecutive non-empty sets $X_i$ & \\

 \begin{Property} \label{obs:c5-xi-xi+1-xi+3} \vspace*{-1.5em} \end{Property}  & If $X_i \not= \emptyset$ and $X_{i+1} \not= \emptyset$, then $X_{i+3} \;\circled{0}\; X_i\bigcup X_{i+1}$ &\\

 \hline
 \begin{Property} \label{obs:c5-ti-xi} \vspace*{-1.5em} \end{Property}  & $T_i \;\circled{1}\; X_i$ &  $T_i \;\circled{1}\; X_{i+4}$\\

 \begin{Property} \label{obs:c5-ti-xi+1} \vspace*{-1.5em} \end{Property} & $T_i \;\circled{$\leq 1$}\; X_{i+1}$ $vv.$ & $T_i \;\circled{$\leq 1$}\; X_{i+3}$ $vv.$\\

 \begin{Property} \label{obs:c5-ti-xi+2} \vspace*{-1.5em} \end{Property} & $T_i \;\circled{0}\; X_{i+2}$ & \\

 \begin{Property} \label{obs:c5-ti-xi-xi+1} \vspace*{-1.5em} \end{Property} & One of the sets $\{T_i, X_i, X_{i+1}\}$ is empty & One of the sets $\{T_{i+2}, X_i, X_{i+1}\}$ is empty \\

 \begin{Property} \label{obs:c5-ti-ti+1-xi} \vspace*{-1.5em} \end{Property} & One of the sets $\{T_i, X_i, T_{i+1}\}$ is empty &\\

 \begin{Property} \label{obs:c5-ti+1-xi-xi+1} \vspace*{-1.5em} \end{Property} & If $T_{i+1} \not= \emptyset$, then $X_i \;\circled{1}\; X_{i+1}$ & \\

 \begin{Property} \label{obs:c5-ti+1-xi-xi+2} \vspace*{-1.5em} \end{Property}& If $T_{i+1} \not= \emptyset$, then $X_i \;\circled{0}\; X_{i+2}$ & If $T_{i+2} \not= \emptyset$, then $X_i \;\circled{0}\; X_{i+2}$ \\

 \begin{Property} \label{obs:c5-xi-ti-ti+2} \vspace*{-1.5em} \end{Property} & If $X_i \neq \emptyset$, then $T_i \;\circled{0}\; T_{i+2}$& \\
 
 \begin{Property} \label{obs:c5-xi-xi+2-ti+4} \vspace*{-1.5em} \end{Property} & If a vertex in $X_i$ is adjacent to vertices in $X_{i+2}$ and in $T_{i+4}$, then its neighbours in $X_{i+2} \cup T_{i+4}$ are a clique & \\

 \begin{Property} \label{obs:c5-xi-xi+1-ti+1-clique} \vspace*{-1.5em} \end{Property} & If $X_i$, $X_{i+1}$, $T_{i+1}$ are all non-empty, then $X_i \cup X_{i+1} \cup T_{i+1}$ is a clique & \\
 \hline
\end{tabular}
\end{center}
\caption{Properties of $(claw, 4K_1, bridge, C_4$-$twin, C_7, C_6)$-free graphs which contain a $C_5$}
\label{tab:C5}
\end{table}

\subsection{Clique-width Labelling}

\subsubsection{The Graph $G[T]$}

\begin{Lemma}\label{lem:c5-GT}
    The graph $G[T]$ has bounded clique-width.
\end{Lemma}
\noindent {\it Proof}. Since we can take any non-empty set $T_i$ to be $T_1$, without loss of generality assume $T_1$ is non-empty. If only two sets $T_i$ are non-empty, by Observations \ref{obs:c5-ti-ti+1} and \ref{obs:c5-ti-ti+2}, they form either a co-join or a vertex in one set is adjacent to at most one vertex in the other set and can therefore be labelled via pairs. Thus, if only two sets $T_i$ are non-empty, the graph $G[T]$ has bounded clique-width.

Now consider when more than two sets $T_i$ are non-empty. Suppose $T_2$ is non-empty. Then by P\ref{obs:c5-no3cons-t}, $T_3$ and $T_5$ are empty. So only $T_4$ may be non-empty. 

By Property P\ref{obs:c5-ti-ti+1}, $T_1 \;\circled{0}\; T_2$. By Property P\ref{obs:c5-ti-ti+2-ti+3}, a vertex in $T_4$ is adjacent to at most one vertex in $T_1 \cup T_2$ and vice versa. So we can label via pairs.

Any other combination of three sets $T_i$ is symmetric to $T_1 \cup T_2 \cup T_4$. Therefore, the graph $G[T]$ has bounded clique-width.\qed

\subsubsection{The Graph $G[X]$}

\begin{Lemma}\label{lem:c5-GX}
    The graph $G[X]$ has bounded clique-width.
\end{Lemma}
\noindent {\it Proof}. Since we can take any non-empty set $X_i$ to be $X_1$, without loss of generality assume $X_1$ is non-empty. If any two sets $X_i$ are non-empty, by Observations \ref{obs:c5-xi-xi+1} and \ref{obs:c5-xi-xi+2}, a vertex in the first set is either adjacent or non-adjacent to at most one vertex in the second set and thus the union of the sets can be labelled via pairs.

Now consider more than two non-empty sets $X_i$. Suppose $X_2$ is non-empty. Then by Property P\ref{obs:c5-no3cons-x}, $X_3$ and $X_5$ are empty. So only $X_4$ may be non-empty.  

By Property P\ref{obs:c5-xi-xi+1}, a vertex in $X_1$ is non-adjacent to at most one vertex in $X_2$. By Property P\ref{obs:c5-xi-xi+1-xi+3}, $X_1 \;\circled{0}\; X_4$ and $X_2 \;\circled{0}\; X_4$. So we can label $X_1 \cup X_2$ via pairs then take the disjoint union with $X_4$.

Any other combination of three sets $X_i$ is symmetric to $X_1 \cup X_2 \cup X_4$. Therefore, $G[X]$ has bounded clique-width.\qed

\subsubsection{The Graph $G[T \cup X]$}

\begin{Lemma}\label{lem:c5-GX2T}
    If $X_i$ and $X_{i+1}$ are non-empty, the graph $G[T \cup X]$ has bounded clique-width.
\end{Lemma}
\noindent {\it Proof}. Assume $X_1$ and $X_2$ are non-empty. By Property P\ref{obs:c5-xi-xi+1}, a vertex in $X_1$ is non-adjacent to at most one vertex in $X_2$ and vice versa. By Property P\ref{obs:c5-no3cons-x}, $X_3$ and $X_5$ are empty. So of the sets $X_i$, only $X_4$ may also be non-empty.

Consider the sets $T_i$. By Property P\ref{obs:c5-ti-xi-xi+1}, $T_1$ and $T_3$ are empty.

Suppose $T_2$ is non-empty. By Property P\ref{obs:c5-xi-xi+1-ti+1-clique}, $T_1 \cup X_1 \cup X_2$ is a clique. Thus if only $X_1, X_2$ and $T_2$ are non-empty, the graph has bounded clique-width. We will refer to the clique formed by $X_1 \cup X_2 \cup T_2$ as $X_1^*$.

Suppose $X_4$ is also non-empty. By Observations \ref{obs:c5-ti-xi+2} and \ref{obs:c5-xi-xi+1-xi+3}, $X_4 \;\circled{0}\; X_1^*$. Then if both $T_4$ and $T_5$ are empty, the graph consists of two disjoint cliques $X_4$ and $X_1^*$, and labelling is trivial. Thus, we may assume $T_4$ or $T_5$ is non-empty.

Suppose $T_4$ is also non-empty, so we have sets $X_1, X_2, X_4, T_2, T_4$ non-empty. By Property P\ref{obs:c5-ti-xi}, $T_4 \;\circled{1}\; X_4$. By Property P\ref{obs:c5-ti-xi+2}, $T_4 \;\circled{0}\; X_1$. Since $X_2$ is non-empty, by Property P\ref{obs:c5-xi-ti-ti+2} and symmetry, $T_4 \;\circled{0}\; T_2$. By Property P\ref{obs:c5-ti-xi+1}, a vertex in $T_4$ is adjacent to at most one vertex in $X_2$ and vice versa. Then the graph $G[T_4 \cup X_1^*]$ satisfies Lemma \ref{lem:cliquewidth-2k} and thus has bounded clique-width. Then when $X_4$ is non-empty, take the disjoint union of labelled graphs $G[T_4 \cup X_1^*]$ and $X_4$ and make the join between vertices with labels (x4, old) and (t4, old). So this graph also has bounded clique-width.

The case when $T_5$ is non-empty is symmetric to the case when $T_4$ is non-empty. By Property P\ref{obs:c5-ti-ti+1-xi} and symmetry, both cannot be non-empty at the same time.

Now suppose that $T_2$ is empty and only $T_4$ of the sets $T_i$ is non-empty, so we have sets $X_1, X_2, T_4$ non-empty. Then $X_1 \cup X_2$ is no longer a clique. By Property P\ref{obs:c5-xi-xi+1}, a vertex in $X_1$ is non-adjacent to at most one vertex in $X_2$, and vice versa. By Property P\ref{obs:c5-ti-xi+2}, there are no edges between $T_4$ and $X_1$. By Property P\ref{obs:c5-ti-xi+1}, a vertex in $T_4$ is adjacent to at most one vertex in $X_2$. Thus we can label the vertices of $X_1 \cup X_2 \cup T_4$ via rows.

Since if $X_4$ is also non-empty it forms a join with $T_4$ and a co-join with $X_1$ and $X_2$, we can take the disjoint union of labelled graphs $X_1 \cup X_2 \cup T_4$ and $X_4$ then make the join between vertices with labels (x4, old) and (t4, old). Thus this graph has bounded clique-width.

Therefore if $X_1$ and $X_2$ are non-empty, the graph $G[T \cup X]$ has bounded clique-width.\qed

\begin{Lemma}\label{lem:c5-GX3T}
    If $X_i$ and $X_{i+2}$ are non-empty, the graph $G[T \cup X]$ has bounded clique-width.
\end{Lemma}
\noindent {\it Proof}. Assume $X_1$ and $X_3$ are non-empty. By Lemma \ref{lem:c5-GX2T}, we can assume no two consecutive sets $X_i$ and $X_{i+1}$ are non-empty. Thus $X_2$, $X_4$, and $X_5$ are all empty. So we have only $X_1$ and $X_3$ non-empty. By Property P\ref{obs:c5-xi-xi+2}, a vertex in $X_1$ is adjacent to at most one vertex in $X_3$ and vice versa.

Suppose $T_1$ is non-empty. Then by Property P\ref{obs:c5-ti-ti+1-xi}, $T_2$ is empty. By Property P\ref{obs:c5-ti-xi}, $T_1 \;\circled{1}\; X_1$. By Property P\ref{obs:c5-ti-xi+2}, $T_1 \;\circled{0}\; X_3$. So if only $T_1$ is non-empty, take the disjoint union of labelled graphs $G[X]$ and $T_1$, then make the join between vertices with labels (t1, old) and (x1, old). Thus $X_1 \cup X_3 \cup T_1$ has bounded clique-width.

Suppose $T_1$ is empty but $T_3$ is non-empty. Then by Property P\ref{obs:c5-ti-ti+1-xi}, $T_4$ is empty. By Property P\ref{obs:c5-ti-xi}, $T_3 \;\circled{1}\; X_3$. By Property P\ref{obs:c5-ti+1-xi-xi+2}, $X_1 \;\circled{0}\; X_3$. By Observations \ref{obs:c5-ti-xi+1} and \ref{obs:c5-ti-xi+1}, a vertex in $T_3$ is adjacent to at most one in $X_1$ and vice versa. Then the graph $G[T_3 \cup X_1]$ satisfies Lemma \ref{lem:cliquewidth-2k} and thus has bounded clique-width. So take the disjoint union of labelled graphs $G[T_3 \cup X_1]$ and $X_3$, then make the join between vertices with labels (t3, old) and (x3, old). Thus $X_1 \cup X_3 \cup T_3$ has bounded clique-width.

Now suppose both $T_1$ and $T_3$ are non-empty. By Property P\ref{obs:c5-xi-ti-ti+2}, $T_1 \;\circled{0}\; T_3$. By Property P\ref{obs:c5-ti-xi+2}, we have $T_1 \circled{0} X_3$. Recall that $T_1 \circled{1} X_1$ by Property P\ref{obs:c5-ti-xi}. Take the disjoint union of labelled graphs $X_1 \cup X_3 \cup T_3$ and $T_1$ then make the join between vertices with labels (t1, old) and (x1, old). Thus, $X_1 \cup X_3 \cup T_1 \cup  T_3$ has bounded clique-width.

Suppose $T_1$ and $T_3$ are empty, but $T_5$ is non-empty. When we have only the sets $X_1$, $X_3$, $T_5$ non-empty, all three sets are related such that a vertex in one set is adjacent to at most one vertex in another set. By Property P\ref{obs:c5-xi-xi+2-ti+4}, if a vertex in either $X_1$ or $X_3$ has neighbours in both the other set $X_i$ and in $T_5$, its neighbourhood is a clique. Then we can label via rows. Thus $X_1 \cup X_3 \cup T_5$ has bounded clique-width.

Suppose both $T_1$ and $T_5$ are non-empty. By Property P\ref{obs:c5-ti-ti+1}, $T_1 \;\circled{0}\; T_5$. By Property P\ref{obs:c5-ti-xi+2}, $T_1 \;\circled{0}\; X_3$. Then we can take the disjoint union of labelled graphs $X_1 \cup X_3 \cup T_5$ and $T_1$, then make the join between vertices with labels (t1, old) and (x1, old). Thus $X_1 \cup X_3 \cup T_1 \cup T_5$ has bounded clique-width.

Suppose both $T_3$ and $T_5$ are non-empty. By Property P\ref{obs:c5-xi-ti-ti+2}, $T_3 \;\circled{0}\; T_5$. By Property P\ref{obs:c5-ti+1-xi-xi+2}, $X_1 \;\circled{0}\; X_3$.  Then we can label via rows. Thus $X_1 \cup X_3 \cup T_3 \cup T_5$ has bounded clique-width.

Suppose all three sets $T_1, T_3, T_5$ are non-empty. Take the disjoint union of labelled graphs $X_1 \cup X_3 \cup T_3 \cup T_5$, then make the join between vertices with labels (t1, old) and (x1, old). Thus $X_1 \cup X_3 \cup T_1 \cup T_3 \cup T_5$ has bounded clique-width.

The case when $T_2$ is non-empty is symmetric to the case when $T_3$ is non-empty. The case when $T_4$ is non-empty is symmetric to the case when $T_1$ is non-empty.

Therefore if $X_1$ and $X_3$ are non-empty, the graph $G[T \cup X]$  has bounded clique-width.\qed

\begin{Lemma}\label{lem:c5-GTX}
    The graph $G[T \cup X]$ has bounded clique-width.
\end{Lemma}
\noindent {\it Proof}. Since we can take any non-empty set $X_i$ to be $X_1$, without loss of generality assume $X_1$ is non-empty. If $X_2$ or $X_5$ are non-empty, refer to Lemma \ref{lem:c5-GX2T}. If $X_3$ or $X_4$ are non-empty, refer to Lemma \ref{lem:c5-GX3T}. So of the sets $X_i$, only $X_1$ is non-empty.

If there is only one set $T_i$ that is non-empty, by Observations \ref{obs:c5-ti-xi}, \ref{obs:c5-ti-xi+1}, \ref{obs:c5-ti-xi+1}, and \ref{obs:c5-ti-xi+2}, it forms either a join or co-join with $X_1$, or a vertex in $T_i$ is adjacent to at most one vertex in $X_1$ and vice versa, in which case we can label via rows. Thus $G[T \cup X]$ has bounded clique-width.

Consider when $T_1$ is non-empty. Then by Property P\ref{obs:c5-ti-ti+1-xi}, $T_2$ is empty. By Property P\ref{obs:c5-ti-xi}, $T_1 \;\circled{1}\; X_1$.

Suppose $T_1$ and $T_3$ are non-empty. By Property P\ref{obs:c5-xi-ti-ti+2}, $T_1 \;\circled{0}\; T_3$. By Observations \ref{obs:c5-ti-xi+1} and \ref{obs:c5-ti-xi+1}, a vertex in $T_3$ is adjacent to at most one vertex in $X_1$ and vice versa. So label $T_3 \cup X_1$ via rows, take the disjoint union with labelled $T_1$, then make the join between vertices with labels (t1, old) and (x1, old). Thus $X_1 \cup T_1 \cup T_3$ has bounded clique-width.

Suppose $T_1$ and $T_5$ are non-empty. By Property P\ref{obs:c5-ti-ti+1}, $T_1 \;\circled{0}\; T_5$. By Observations \ref{obs:c5-ti-xi+1} and \ref{obs:c5-ti-xi+1}, a vertex in $T_5$ is adjacent to at most one vertex in $X_1$ and vice versa. So label $T_5 \cup X_1$ via rows, take the disjoint union with labelled $T_1$, then make the join between vertices with labels (t1, old) and (x1, old). Thus $X_1 \cup T_1 \cup T_5$ has bounded clique-width.

Suppose $T_1$, $T_3$, and $T_5$ are all non-empty. Then we can label $X_1 \cup T_3 \cup T_5$ via rows, label $T_1$, then make the join between vertices with labels (t1, old) and (x1, old). Thus $X_1 \cup T_1 \cup T_3 \cup T_5$ has bounded clique-width.

Now consider when $T_4$ is non-empty. By Property P\ref{obs:c5-ti-xi+2}, $T_4 \;\circled{0}\; X_1$.

Suppose $T_1$ is also non-empty. By Property P\ref{obs:c5-ti-ti+2}, a vertex in $T_1$ is adjacent to at most one vertex in $T_4$ and vice versa. Then we can label $T_1 \cup T_4$ via pairs, take the disjoint union with labelled $X_1$, then make the join between vertices with labels (t1, old) and (x1, old).

By Property P\ref{obs:c5-no3cons-t}, one of $T_3$ or $T_5$ must be empty. The cases are symmetric, so without loss of generality assume $T_3$ is non-empty. By Property P\ref{obs:c5-ti-ti+1}, $T_4 \;\circled{0}\; T_3$.  If $T_1$ is empty, then we can take the disjoint union of labelled graphs $G[T_3 \cup X_1]$ and $T_4$. Thus $X_1 \cup T_3 \cup T_4$ has bounded clique-width. If $T_1$ is non-empty, then we can take the disjoint union of labelled graphs $G[T_3 \cup X_1]$ and $G[T_1 \cup T_4]$ labelled via pairs, then make the join between labels (t1, old) and (x1, old). Thus $X_1 \cup T_1 \cup T_3 \cup T_4$ has bounded clique-width

The case when $T_2$ is non-empty is symmetric to the case when $T_1$ is non-empty.

Therefore the graph $G[T \cup X]$ has bounded clique-width.\qed

\subsection{Graphs containing a $C_5$}

The purpose of the previous results was to prove the following theorem:
\begin{Theorem}\label{thm-c5}
 Let $G$ be a $(claw, 4K_1, bridge, C_4$-$twin)$-free graph. If $G$ contains a $C_5$, then $G$ has bounded clique-width.
\end{Theorem}
\noindent {\it Proof.} By Theorems \ref{thm-c7} and \ref{thm-c6}, 
we may assume that $G$ contains no $C_7$ and no $C_6$. In Lemma \ref{lem:c5-GT}, we proved that if $G$ consists of only the $C_5$ and the sets $T_i$, then $G$ has bounded clique-width. In Lemma \ref{lem:c5-GX}, we proved that if $G$ consists of only the $C_5$ and the sets $X_i$, then $G$ has bounded clique-width. In Lemma \ref{lem:c5-GTX}, we proved that if $G$ consists of the $C_5$ and the sets $X_i$ and $T_i$, then $G$ has bounded clique-width. In Observation \ref{obs:c5-z}, we showed that if the set $Z$ is not empty, then the graph has bounded clique-width.

From Observations \ref{obs:c5-sets} and \ref{obs:c5-y-empty}, we know no other configuration of vertices is possible. Therefore, if $G$ contains a $C_5$, $G$ has bounded clique-width. \qed

\section{Proof of the main result}\label{sec:proof-of-main-results}

The purpose of this section is to synthesize the previous results to show that $(claw, 4K_1, bridge, C_4$-$twin)$-free graphs can be coloured in polynomial time using the structure results from Sections \ref{sec:c7}, \ref{sec:c6}, and \ref{sec:c7}. We will prove the two theorems mentioned in Section \ref{sec:intro}.

\noindent \textbf{Theorem \ref{thm:structure} {\it 
	Let $G$ be a  \CCD-free graph. Then one of the following holds.
	\begin{description}
		\item[(i)] $G$ is perfect. 
		\item[(ii)] $G$ has bounded clique-width.
	\end{description} }
}

\noindent {\it Proof of Theorem \ref{thm:structure}.} Let $G$ be a \CCD-free graph. We may assume that $G$ is not perfect, for otherwise we 
are done. Suppose $G$ contains an odd hole $H$. Since $G$ is $4 K_1$-free, $H$ is a $C_5$ or $C_7$. From Theorems \ref{thm-c5} and \ref{thm-c7}, we know that $G$ has bounded clique-width.  Now, we may assume that $G$ contains an odd antihole $F$. Since the $C_5$ and its complement are isomorphic, if $F$ has five vertices, by Theorem \ref{thm-c5}, we are done. So we can assume $F$ is of size at least $7$. Let $\{1, 2, 3, 4, 5, 6, 7, \dots, k\}$ be the vertices of $F$ where the non-edges are $(i, i+1)$ with the subscript taken modulo $k$. Then the set $\{1,6,2,5,3\}$ induces a $C_4$-$twin$, a contradiction. \qed

\noindent \textbf{Theorem \ref{thm:main}} {\it There is a polynomial-time algorithm to colour a $(claw, 4K_1, bridge, C_4$-$twin)$-free graph.}

\noindent {\it Proof}. Let $G$ be a $(claw, 4K_1, bridge, C_4$-$twin)$-free graph. If $G$ is perfect, then by the result of Hsu from \cite{Hsu1981}, $G$ can be coloured in polynomial time, we are done. So by Theorem \ref{thm:structure}, $G$ has bounded clique-width. Then we can colour $G$ in polynomial time using the algorithm from \cite{Rao}.  \qed

\section{Conclusion and open problems}\label{sec:conclusions}

The complexity of vertex colouring for $(claw, 4K_1)$-free graphs is one of only three open problems for four-vertex forbidden induced graphs. The problem has been solved for some subclasses with additional forbidden subgraphs. In this paper, we proved colouring is polynomial-time solvable for the subclass of $(claw, 4K_1, bridge, C_4$-$twin)$-free graphs by proving graphs in this class are either perfect or have bounded clique-width and can thus be coloured using either Hsu's algorithm from \cite{Hsu1981} or Rao's algorithm from \cite{Rao}.

The $C_4$-$twin$ was used in six observations in the case of the $C_6$ and six observations in the case of the $C_5$, for a total of twelve observations. We did not use the $C_4$-$twin$ in the case of the $C_7$. Thus, we believe that omitting the $C_4$-$twin$ entirely and showing that $(claw, 4K_1, bridge)$-free graphs have bounded clique-width is possible. Thus we conjecture the following:

\begin{Conjecture}
    Let G be a $(claw, 4K_1, bridge)$-free graph. Then one of the following holds.
    \begin{itemize}
        \item G is perfect.
        \item G has bounded clique-width.
    \end{itemize}
\end{Conjecture}

Since we have shown this is true when the graph contains a $C_7$, it would only be a matter of proving it for when the graph contains a $C_5$ or a $C_6$. Since the $C_4$-$twin$ is not used in every observation regarding the properties of the graphs when they contain a $C_5$ or a $C_6$, many of the observations in this paper hold regardless of whether the $C_4$-$twin$ is no longer included as a forbidden induced subgraph. New observations would need to be made only in place of the twelve occurrences of the $C_4$-$twin$. It may be possible to prove it while omitting an investigation of the case of the $C_6$ as well, since the $C_6$ was only added to make the case of the $C_5$ more manageable. It may be possible to investigate the $C_5$ without forbidding the $C_6$.

To conclude, we note that the complexity of colouring $(claw, 4K_1)$-free graphs in general has remained an open problem, along with the two other open problems for four-vertex forbidden induced subgraphs outlined in \cite{Lozin} and \cite{DHH}, namely $(claw, 4K_1, co$-$diamond)$-free and $(4K_1, C_4)$-free graphs.


\begin{appendices}
In this appendix, we give proofs for properties (in the three tables) that were not given in the main body of the paper. 

\section{Properties of sets in the case of the $C_7$}



    

\noindent \textbf{P\ref{obs:c7-yi-yi+3}}\textit{
    $Y_i \;\circled{0}\; Y_{i+3}$  for all $i$. By symmetry, $Y_i \;\circled{0}\; Y_{i+4}$.
}

\noindent {\it Proof}. If there is a vertex $a \in Y_1$ that is adjacent to a vertex $b \in Y_4$, then the set $\{a,b,1,3\}$ induces a $claw$.\qed

\noindent \textbf{P\ref{obs:c7-yi-xi}}\textit{
	$Y_i \;\circled{1}\; X_i$  for all $i$. By symmetry, $Y_i \;\circled{1}\; X_{i+1}$.
}

\noindent {\it Proof}. If there is a vertex $a \in Y_1$ that is non-adjacent to a vertex $b \in X_1$, then the set $\{a,b,5,7\}$ induces a $4K_1$.\qed

\noindent \textbf{P\ref{obs:c7-yi-xi+2}}\textit{
	If $Y_i \not= \emptyset$, then $X_{i+2} = \emptyset$. By symmetry, if $Y_i \not= \emptyset$, then $X_{i+6} = \emptyset$.
}

\noindent {\it Proof}. Suppose there are two vertices $a,b \in X_3$ and two vertices $c,d \in Y_1$. If ${a,b} \;\circled{0}\; {c,d}$, then the set $\{a,b,3,4,c,d\}$ induces a $bridge$. If ${a,b} \;\circled{1}\; {c,d}$, then the set $\{a,b,c,d,1,2\}$ induces a $bridge$.

Suppose a vertex $a \in X_3$ is adjacent to two vertices $b,c \in Y_1$. Then the set $\{1,2,b,c,a,4\}$ induces a $bridge$. Therefore, a vertex in $X_3$ is adjacent to at most one vertex in $Y_1$.

Suppose a vertex $a \in Y_1$ is adjacent to two vertices $b,c \in X_3$. There exists a vertex $d \in Y_1$ which is not adjacent to either of $b,c$. Then the set $\{b,c,a,3,d,2\}$ induces a $bridge$. So, we know that

\begin{equation}\label{equ:y1-x3-see-one}
\mbox{A vertex in $Y_i$ is adjacent to at most one vertex in $X_{i+2}$ and vice versa.}
\end{equation}

Consider a vertex $a \in Y_1$. Then by (\ref{equ:y1-x3-see-one}), for some vertex $a' \in X_3$, $a$ misses all of $X_3 - a'$. Consider a vertex $b \in Y_1$. By (\ref{equ:y1-x3-see-one}), $b$ is non-adjacent to at least two vertices $c,d \in X_3$. Then ${a,b} \;\circled{0}\; {c,d}$, so the set $\{a,b,3,4,c,d\}$ induces a $bridge$. Contradiction. Thus the first part of the observation holds.\qed

\noindent \textbf{P\ref{obs:c7-yi-xi+3}}\textit{
	$Y_i \;\circled{0}\; X_{i+3}$  for all $i$. By symmetry, $Y_i \;\circled{0}\; X_{i+5}$.
}

\noindent {\it Proof}. If there is a vertex $a \in Y_1$ that is adjacent to a vertex $b \in X_4$, then the set $\{a,b,1,3\}$ induces a $claw$.\qed

\noindent \textbf{P\ref{obs:c7-yi-xi+4}}\textit{
	$Y_i \;\circled{0}\; X_{i+4}$  for all $i$.
}

\noindent {\it Proof}. If there is a vertex $a \in Y_1$ that is adjacent to a vertex $b \in X_5$, then the set $\{a,b,1,3\}$ induces a $claw$.\qed

\section{Properties of sets in the case of the $C_6$}

\noindent \textbf{P\ref{obs:c6-ti-ti+1}}\textit{
    If $T_i \not= \emptyset$, then $T_{i+1} = \emptyset$ and vice versa. By symmetry, if $T_i \not= \emptyset$, then $T_{i+5} = \emptyset$ and vice versa.
}

\noindent {\it Proof}. If there is a vertex $a \in T_1$ that is adjacent to a vertex $b \in T_2$, then the set $\{1,a,b,3,4,5,6\}$ induces a $C_7$. Conversely, if $a$ and $b$ are not adjacent, then the set $\{a,b,4,6\}$ induces a $4K_1$.\qed

\noindent \textbf{P\ref{obs:c6-ti-ti+2}}\textit{
    A vertex in $T_i$ is adjacent to at most one vertex in $T_{i+2}$ and vice versa. By symmetry, a vertex in $T_i$ is adjacent to at most one vertex in $T_{i+4}$ and vice versa.
}

\noindent {\it Proof}. If there is a vertex $a \in T_1$ that is adjacent to two vertices $b,c \in T_3$, then the set $\{a,2,3,b,c\}$ induces a $C_4$-$twin$.\qed

\noindent \textbf{P\ref{obs:c6-ti-ti+3}}\textit{
    If $T_i \not= \emptyset$, then $T_{i+3} = \emptyset$.
}

\noindent {\it Proof}. If there is a vertex $a \in T_1$ that is not adjacent to a vertex $b \in T_4$, then the set $\{a,3,b,6\}$ induces a $4K_1$. So we must have $T_1 \;\circled{1}\; T_4$. Since $T_1$ and $T_4$ are big, there are two vertices $c \in T_1$ and $d \in T_4$ such that $\{a,b,c,d\}$ is a clique. Then the set $\{1,2,a,c,b,d\}$ induces a $bridge$.\qed

\noindent \textbf{P\ref{obs:c6-ti-ti+2-ti+4}}\textit{
    If a vertex in $T_i$ is adjacent to a vertex in $T_{i+2}$ and a vertex in $T_{i+4}$, those vertices must also be adjacent. By symmetry, the same is true if a vertex in $T_{i+2}$ is adjacent to a vertex in $T_i$ and a vertex in $T_{i+4}$ and if a vertex in $T_{i+4}$ is adjacent to a vertex in $T_i$ and a vertex in $T_{i+2}$.
}

\noindent {\it Proof}. Suppose there is a vertex $a \in T_1$ that is adjacent to two vertices $b \in T_3$ and $c \in T_5$. If $b$ and $c$ are not adjacent, then the set $\{a,2,b,c\}$ induces a $claw$.\qed

\noindent \textbf{P\ref{obs:c6-Xi-Xi+1}}\textit{
    A vertex in $X_i$ is non-adjacent to at most one vertex in $X_{i+1}$ and vice versa. By symmetry, a vertex in $X_i$ is non-adjacent to at most one vertex in $X_{i+5}$ and vice versa.
}

\noindent {\it Proof}. If there is a vertex $a \in X_1$ that is non-adjacent to two vertices $b,c \in X_2$, then we claim that every other vertex $d \in X_1$ must be adjacent to exactly one of $\{b,c\}$. Otherwise if $d$ is non-adjacent to both, the set $\{a,d,2,3,b,c\}$ induces a $bridge$ or, if $d$ is adjacent to both, the set $\{1,a,d,2,b,c\}$ induces a $bridge$. WLOG, assume $d$ is adjacent to $b$. Then by the assumption that $X_1$ is big, there is another vertex $e \in X_1$ that is also adjacent to $b$. Then the set $\{c,4,b,3,e,d\}$ induces a $bridge$.\qed

\noindent \textbf{P\ref{obs:c6-Xi-Xi+2}}\textit{
    A vertex in $X_i$ is adjacent to at most one vertex in $X_{i+2}$ and vice versa. By symmetry, a vertex in $X_i$ is adjacent to at most one vertex in $X_{i+4}$ and vice versa.
}

\noindent {\it Proof}. If there is a vertex $a \in X_1$ that is adjacent to two vertices $b,c \in X_3$, then we claim that every other vertex $d \in X_1$ must be adjacent to exactly one of $\{b,c\}$. Otherwise if $d$ is non-adjacent to both, then the set $\{d,2,a,3,b,c\}$ induces a $bridge$ or if $d$ is adjacent to both then the set $\{1,2,a,d,b,c\}$ induces a $bridge$ . WLOG, assume $d$ is adjacent to $b$. Then by the assumption that $X_1$ is big, there is another vertex $e \in X_1$ that is also adjacent to $b$. Then the set $\{d,e,b,3,c,4\}$ induces a $bridge$.\qed

\noindent \textbf{P\ref{obs:c6-Xi-Xi+3}}\textit{
    $X_i \;\circled{0}\; X_{i+3}$ for all $i$. 
}

\noindent {\it Proof}. If there is a vertex $a \in X_1$ that is adjacent to a vertex $b \in X_4$, then the set $\{a,b,1,3\}$ induces a $claw$.\qed

\noindent \textbf{P\ref{obs:c6-no3cons-x}}\textit{
    There are no 3 consecutive non-empty sets $X_i$.
}

\noindent {\it Proof}. Suppose $X_1, X_2, X_3$ are non-empty. Since the sets $X_i$ are all big, by Property P\ref{obs:c6-Xi-Xi+2}, there are vertices $a \in X_1$ and $c \in X_3$ that are not adjacent. Since the set $X_2$ is big, by Property P\ref{obs:c6-Xi-Xi+1}, there are many vertices $b \in X_2$ that are adjacent to both $a$ and $c$. Then the set $\{a,2,b,3,c,4\}$ induces a $bridge$.\qed

\noindent \textbf{P\ref{obs:c6-yi-yi+1}}\textit{
    A vertex in $Y_i$ is non-adjacent to at most one vertex in $Y_{i+1}$ and vice versa. By symmetry, a vertex in $Y_i$ is non-adjacent to at most one vertex in $Y_{i+5}$ and vice versa.
}

\noindent {\it Proof}. If there is a vertex $a \in Y_1$ that is non-adjacent to two vertices $b,c \in Y_2$, then we claim that every other vertex $d \in Y_1$ must be adjacent to exactly one of $\{b,c\}$. Otherwise, if $d$ is non-adjacent to both, then the set $\{d,a,2,3,c,b\}$ induces a $bridge$ or, if $d$ is adjacent to both, then the set $\{a,1,d,2,c,b\}$ induces a $bridge$. WLOG, assume $d$ is non-adjacent to $b$ and adjacent to $c$. Then by the assumption that $Y_1$ is big, there is another vertex $e \in Y_1$ that is also non-adjacent to $b$ and adjacent to $c$. Then since $b$ is non-adjacent to both $d$ and $e$, by symmetry, every other vertex in $Y_2$ must be adjacent to exactly one of $\{d,e\}$. But $c$ is adjacent to both. Contradiction.\qed

\noindent \textbf{P\ref{obs:c6-yi-yi+2}a}\textit{
    A vertex in $Y_i$ is adjacent to at most one vertex in $Y_{i+2}$ and vice versa. By symmetry, a vertex in $Y_i$ is adjacent to at most one vertex in $Y_{i+4}$ and vice versa.
}

\noindent {\it Proof}. If there is a vertex $a \in Y_1$ that is adjacent to two vertices $b,c \in Y_3$, then the set $\{a,b,c,1,6\}$ induces a $C_4$-$twin$.\qed

\noindent \textbf{P\ref{obs:c6-yi-yi+2}}\textit{
    If $Y_i \not= \emptyset$, then $Y_{i+2} = \emptyset$ and vice versa. By symmetry, if $Y_i \not= \emptyset$, then $Y_{i+4} = \emptyset$ and vice versa.
}

\noindent {\it Proof}. Consider a vertex $a \in Y_1$. By Property P\ref{obs:c6-yi-yi+2}a, $a$ is non-adjacent to some vertices $b,c \in Y_1$. Consider a second vertex $d \in Y_1$. If $d$ is non-adjacent to both $b$ and $c$, the set $\{a,d,3,4,b,c\}$ induces a bridge. This along with Property P\ref{obs:c6-yi-yi+2}a means $d$ must be adjacent to exactly one of $\{b,c\}$. This is true of every vertex in $Y_1 - a$. Then since $Y_1$ is big, at least one of $b,c$ is adjacent to more than one vertex in $Y_1$, contradicting Property P\ref{obs:c6-yi-yi+2}a. Therefore, both sets cannot be non-empty.\qed

\noindent \textbf{P\ref{obs:c6-yi-yi+3}}\textit{
    A vertex in $Y_i$ is adjacent to at most one vertex in $Y_{i+3}$ and vice versa.
}

\noindent {\it Proof}. If there is a vertex $a \in Y_1$ that is adjacent to two vertices $b,c \in Y_4$, then we claim that every other vertex $d \in Y_1$ must be adjacent to exactly one of $\{b,c\}$. Otherwise, if $d$ is non-adjacent to both, then the set $\{d,3,a,4,b,c\}$ induces a $bridge$ or, if $d$ is adjacent to both, then the set $\{2,3,a,d,b,c\}$ induces a $bridge$. WLOG, assume $d$ is non-adjacent to $b$ and adjacent to $c$. Then by the assumption that $Y_1$ is big, there is another vertex $e \in Y_1$ that is also non-adjacent to $b$ and adjacent to $c$. Then since $b$ is non-adjacent to both $d$ and $e$, by symmetry, every other vertex in $Y_4$ must be adjacent to exactly one of $\{d,e\}$. But $c$ is adjacent to both. Contradiction.\qed

\noindent \textbf{P\ref{obs:c6-xi-ti}}\textit{
    $X_i \;\circled{1}\; T_i$ for all $i$. By symmetry, $X_i \;\circled{1}\; T_{i+1}$ for all $i$.
}

\noindent {\it Proof}. If there is a vertex $a \in X_1$ that is not adjacent to a vertex $b \in T_1$, then the set $\{a,b,4,6\}$ induces a $4K_1$.\qed

\noindent \textbf{P\ref{obs:c6-xi-ti+2}}\textit{
    A vertex in $X_i$ is adjacent to at most one vertex in $T_{i+2}$. By symmetry, a vertex in $X_i$ is adjacent to at most one vertex in $T_{i+5}$.
}

\noindent {\it Proof}. If there is a vertex $a \in X_1$ that is adjacent to two vertices $b,c \in T_3$, then we claim that any other vertex $d \in X_1$ must be adjacent to exactly one of $b,c$. Otherwise, if $d$ is non-adjacent to both, then the set $\{d,2,a,3,b,c\}$ induces a $bridge$ or, if $d$ is adjacent to both, then the set $\{1,2,a,d,b,c\}$ induces a $bridge$. WLOG, assume $d$ is non-adjacent to $b$ and adjacent to $c$. Then by the assumption that $X_1$ is big, there is another vertex $e \in X_1$ that is also non-adjacent to $b$ and adjacent to $c$. Then the set $\{c,4,b,3,d,e\}$ induces a $bridge$.\qed

\noindent \textbf{P\ref{obs:c6-xi-ti+2}vv.}\textit{
    A vertex in $T_{i+2}$ is adjacent to at most one vertex in $X_i$. By symmetry, a vertex in $T_{i+5}$ is adjacent to at most one vertex in $X_i$.
}

\noindent {\it Proof}. Suppose there is a vertex $a \in T_3$ that is adjacent to two vertices $b,c \in X_1$. Since $T_3$ is big, there exists two more vertices $d,e \in T_3$. From Property P\ref{obs:c6-xi-ti+2}, $d,e$ are not adjacent to $b$ or $c$. Then the set $\{b,c,3,a,d,e\}$ induces a $bridge$.\qed

\noindent \textbf{P\ref{obs:c6-xi-ti+3}}\textit{
    $X_i \;\circled{0}\; T_{i+3}$ for all $i$. By symmetry, $X_i \;\circled{0}\; T_{i+4}$ for all $i$.
}

\noindent {\it Proof}. If there is a vertex $a \in X_1$ that is adjacent to a vertex $b \in T_4$, then the set $\{a,1,3,b\}$ induces a $claw$.\qed

\noindent \textbf{P\ref{obs:c6-xi-ti-ti+2}}\textit{
    If $X_i$ is non-empty, then $T_i \;\circled{0}\; T_{i+2}$. By symmetry, if $X_i$ is non-empty, then $T_{i+1} \;\circled{0}\; T_{i+5}$.
}

\noindent {\it Proof}. If there are vertices $a \in T_1$, $b \in X_1$, then we may choose $c \in T_3$ to be any of the many vertices that, by Property P\ref{obs:c6-xi-ti+2}, are not adjacent to $b$. If $a$ and $c$ are adjacent, then the set $\{a,b,3,c,2\}$ induces a $C_4$-$twin$.\qed

\noindent \textbf{P\ref{obs:c6-xi-xi+1-ti}}\textit{
    If $X_i$ and $X_{i+1}$ are both non-empty, then $T_i$ is empty. By symmetry, if $X_i$ and $X_{i+1}$ are non-empty, then $T_{i+2}$ is empty.
}

\noindent {\it Proof}. Suppose there are vertices $a \in T_1$ and $b \in X_1$. Choose $c \in X_2$ as one of the many vertices that, by Property P\ref{obs:c6-Xi-Xi+1}, is adjacent to $b$. By Property P\ref{obs:c6-xi-ti}, $a$ and $b$ are adjacent. By Property P\ref{obs:c6-xi-ti+2}, $a$ and $c$ are not adjacent. Then the set $\{a,1,b,2,c,3\}$ induces a $bridge$.\qed

\noindent \textbf{P\ref{obs:c6-ti-xi-ti+2}}\textit{
    If $T_i$ is non-empty, then $X_i \;\circled{0}\; T_{i+2}$.
}

\noindent {\it Proof}. Suppose there are vertices $a \in T_1$, $b \in X_1$, and $c \in T_3$. By Property P\ref{obs:c6-xi-ti-ti+2}, $a$ and $c$ are not adjacent. By Property P\ref{obs:c6-xi-ti}, $a$ and $b$ are adjacent. Then if $b$ and $c$ are adjacent, the set $\{a,1,b,2,c,3\}$ induces a $bridge$.\qed

\noindent \textbf{P\ref{obs:c6-ti-xi+1-xi+2}}\textit{
    If $T_i$ is non-empty, then $X_{i+1} \;\circled{1}\; X_{i+2}$.
}

\noindent {\it Proof}. Consider $a \in X_2$ and $b \in X_3$. By Property P\ref{obs:c6-xi-ti+2}, a choice of $c \in T_1$ exists such that $c$ is non-adjacent to $a$. By Property P\ref{obs:c6-xi-ti+3}, $c$ is also non-adjacent to $b$. Then if $a$ and $b$ are non-adjacent, the set $\{a,b,c,6\}$ induces a $4K_1$.\qed

\noindent \textbf{P\ref{obs:c6-xi-xi+2-ti+2}}\textit{
    A vertex in $X_i$ cannot be adjacent to both a vertex in $T_{i+2}$ and a vertex in $X_{i+2}$. By symmetry, a vertex in $X_{i+2}$ cannot be adjacent to both a vertex in $T_{i+1}$ and a vertex in $X_i$.
}

\noindent {\it Proof}. Consider $a, b \in X_1$, $c \in T_3$ and $d \in X_3$. By Property P\ref{obs:c6-xi-ti}, $c$ and $d$ are adjacent. If $a$ is adjacent to both $c$ and $d$, then by Observations \ref{obs:c6-Xi-Xi+2} and \ref{obs:c6-xi-ti+2}vv., $b$ is adjacent to both $c$ and $d$. Then the set $\{b,2,a,3,c,d\}$ induces a $bridge$.\qed

\noindent \textbf{P\ref{obs:c6-xi-xi+2-ti+4}}\textit{
    A vertex in $X_i$ cannot be adjacent to both a vertex in $X_{i+2}$ and a vertex in $T_{i+4}$.
}

\noindent {\it Proof}. Consider $a \in X_1$, $b \in X_3$ and $c \in T_5$. By Property P\ref{obs:c6-xi-ti+3}, $a$ and $c$ are not adjacent. If $b$ is adjacent to both $a$ and $c$, then the set $\{b,a,4,c\}$ induces a $claw$.\qed

\noindent \textbf{P\ref{obs:c6-xi+2-xi-ti+4}}\textit{
    A vertex in $X_{i+2}$ cannot be adjacent to both a vertex in $T_{i+4}$ and a vertex in $X_i$.
}

\noindent {\it Proof}. Consider $a \in X_3$, $b \in X_1$ and $c \in T_5$. By Property P\ref{obs:c6-xi-ti+3}, $b$ and $c$ are not adjacent. If $a$ is adjacent to both $b$ and $c$, then the set $\{a,b,4,c\}$ induces a $claw$.\qed

\noindent \textbf{P\ref{obs:c6-xi-ti+1-ti+5-a}}\textit{
    If $X_i, T_{i+1}, T_{i+5}$ are all non-empty, $T_{i+1} \;\circled{0}\; T_{i+5}$.
}

\noindent {\it Proof}. Consider $a \in T_2$ and $b \in T_6$. Choose $c \in X_1$ to be one of the many vertices that, by Property P\ref{obs:c6-xi-ti+2}, is not adjacent to $b$. By Property P\ref{obs:c6-xi-ti}, $a$ and $c$ are adjacent. If $a$ and $b$ are adjacent, then the set $\{1,2,a,b,c\}$ induces a $C_4$-$twin$.\qed

\noindent \textbf{P\ref{obs:c6-yi-ti}}\textit{
    $Y_i \;\circled{1}\; T_i$ for all $i$. By symmetry, $Y_i \;\circled{1}\; T_{i+2}$ for all $i$.
}

\noindent {\it Proof}. If there is a vertex $a \in Y_1$ that is not adjacent to a vertex $b \in T_1$, then the set $\{1,a,b,6\}$ induces a $claw$.\qed

\noindent \textbf{P\ref{obs:c6-yi-ti+1}}\textit{
    A vertex in $Y_i$ is non-adjacent to at most two vertices in $T_{i+1}$.
}

\noindent {\it Proof}. If there is a vertex $a \in Y_1$ that is non-adjacent to two vertices $b,c \in T_2$, then we claim that every other vertex $d \in Y_1$ must be adjacent to exactly one of $\{b,c\}$. Otherwise, if $d$ is non-adjacent to both, then the set $\{a,d,2,3,b,c\}$ induces a $bridge$ or, if $d$ is adjacent to both, then the set $\{a,1,2,d,b,c\}$ induces a $bridge$. Then if $a$ is non-adjacent to a third vertex $e \in T_2$, any other vertex $d$ must be adjacent to exactly one of $\{b,c\}$, one of $\{b,e\}$, and one of $\{e,c\}$, which is not possible.\qed

\noindent \textbf{P\ref{obs:c6-yi-ti+1-b}}\textit{
    If a vertex in $Y_i$ is non-adjacent to two vertices in $T_{i+1}$, one of those vertices must be adjacent to all of $Y_i$ while the other must be non-adjacent to all of $Y_i$.
}

\noindent {\it Proof}. If there is a vertex in $a \in Y_1$ that is non-adjacent to two vertices $b,c \in T_2$, then we know from Property P\ref{obs:c6-yi-ti+1} that every other vertex in $Y_1$ is adjacent to exactly one of $\{b,c\}$. We will show that it must be the same one.

Consider $d,e \in Y_1$. Let $d$ be adjacent to $b$ and non-adjacent to $c$ while $e$ is non-adjacent to $b$ and adjacent to $c$. Now consider another vertex $f \in Y_1$. We know $f$ must also be adjacent to exactly one of $\{b,c\}$. WLOG, assume it is adjacent to $b$ and non-adjacent to $c$. Then the set $\{b,c,e,d,f\}$ induces a $C_4$-$twin$.\qed

\noindent \textbf{P\ref{obs:c6-yi-ti+3}}\textit{
    $Y_i \;\circled{0}\; T_{i+3}$ for all $i$. By symmetry,  $Y_i \;\circled{0}\; T_{i+5}$ for all $i$.
}

\noindent {\it Proof}. If there is a vertex $a \in Y_1$ that is adjacent to a vertex $b \in T_4$, then the set $\{a,1,3,b\}$ induces a $claw$.\qed

\noindent \textbf{P\ref{obs:c6-yi-ti+4}}\textit{
    $Y_i \;\circled{0}\; T_{i+4}$ for all $i$.
}

\noindent {\it Proof}. If there is a vertex $a \in Y_1$ that is adjacent to a vertex $b \in T_5$, then the set $\{a,1,3,b\}$ induces a $claw$.\qed

\noindent \textbf{P\ref{obs:c6-yi-ti-ti+2}}\textit{
    If $Y_i \not= \emptyset$, then one of $T_i$ or $T_{i+2} = \emptyset$.
}

\noindent {\it Proof}. Suppose there are vertices $a \in T_1$, $b,c \in Y_1$. By Property P\ref{obs:c6-ti-ti+2}, a choice of $d \in T_3$ exists such that it is non-adjacent to $a$. By Property P\ref{obs:c6-yi-ti}, $b,c$ are adjacent to both $a$ and $d$. Then the set $\{a,1,b,c,d,3\}$ induces a $bridge$.\qed

\noindent \textbf{P\ref{obs:c6-yi-yi+1-ti}}\textit{
    If $Y_i \not= \emptyset$ and $Y_{i+1} \not= \emptyset$ , then $T_i = \emptyset$. By symmetry, if $Y_i \not= \emptyset$ and $Y_{i+1} \not= \emptyset$ , then $T_{i+3} = \emptyset$.
}

\noindent {\it Proof}. Suppose there are vertices $a \in T_1$ and $b \in Y_1$. By Property P\ref{obs:c6-yi-yi+1}, a choice of $c \in Y_2$ exists such that it is adjacent to $b$. By Property P\ref{obs:c6-yi-ti}, $a$ and $b$ are adjacent. By Property P\ref{obs:c6-yi-ti+3}, $a$ and $c$ are not adjacent. Then the set $\{a,1,b,2,c,3\}$ induces a $bridge$.\qed

\noindent \textbf{P\ref{obs:c6-yi-ti+1-ti+3}}\textit{
    If $Y_i$ and $T_{i+3}$ are both non-empty, $Y_i \;\circled{1}\; T_{i+1}$. By symmetry, if $Y_i$ and $T_{i+5}$ are both non-empty, $Y_i \;\circled{1}\; T_{i+1}$.
}

\noindent {\it Proof}. If there are vertices $a \in Y_1$, $b \in T_4$ and $c \in T_2$ such that $a$ and $c$ are not adjacent, then the set $\{a,b,c,6\}$ induces a $4K_1$.\qed

\noindent \textbf{P\ref{obs:c6-yi-yi+1-ti+1}}\textit{
    If $Y_{i+1} \not= \emptyset$, a vertex in $Y_i$ is non-adjacent to at most one vertex in $T_{i+1}$.
}

\noindent {\it Proof}. If there is a vertex $a \in Y_1$ that is non-adjacent to two vertices $b,c \in T_2$, consider a vertex $d \in Y_2$. By Property P\ref{obs:c6-yi-yi+1}, a choice of $d$ exists such that it is adjacent to $a$. By Property P\ref{obs:c6-yi-ti}, $d$ is adjacent to both $b$ and $c$. Then the set $\{a,4,d,3,b,c\}$ induces a $bridge$.\qed

\noindent \textbf{P\ref{obs:c6-yi-ti+1-ti+5}}\textit{
    If $Y_i \not= \emptyset$, $T_{i+1} \;\circled{0}\; T_{i+5}$. By symmetry, if $Y_i \not= \emptyset$, $T_{i+1} \;\circled{0}\; T_{i+3}$.
}

\noindent {\it Proof}. Suppose there is a vertex $a \in T_2$ that is adjacent to a vertex $b \in T_6$. Consider a vertex $c \in Y_1$. By Property P\ref{obs:c6-yi-ti+1-ti+3}, $a$ and $c$ are adjacent. By Property P\ref{obs:c6-yi-ti+3}, $b$ and $c$ are not adjacent. Then the set $\{1,2,a,b,c\}$ induces a $C_4$-$twin$.\qed

\noindent \textbf{P\ref{obs:c6-yi-xi}}\textit{
    $Y_{i} \;\circled{1}\; X_{i}$ for all $i$. By symmetry, $Y_{i} \;\circled{1}\; X_{i+1}$ for all $i$.
}

\noindent {\it Proof}. If there is a vertex $a \in Y_1$ that is not adjacent to a vertex $b \in X_1$, then the set $\{1,a,b,6\}$ induces a $claw$.

\noindent \textbf{P\ref{obs:c6-yi-xi+2}a}\textit{
    A vertex in $Y_i$ is non-adjacent to at most one vertex in $X_{i+2}$. By symmetry, a vertex in $Y_i$ is non-adjacent to at most one vertex in $X_{i+5}$.
}

\noindent {\it Proof}. If there is a vertex $a \in Y_1$ that is non-adjacent to two vertices $b,c \in X_3$, then we claim that every other vertex $d \in Y_1$ must be adjacent to exactly one of $\{b,c\}$. Otherwise, if $d$ is non-adjacent to both, then the set $\{a,d,3,4,b,c\}$ induces a $bridge$ or, if $d$ is adjacent to both, then the set $\{a,2,d,3,b,c\}$ induces a $bridge$. WLOG, assume $d$ is non-adjacent to $b$ and adjacent to $c$. Then by the assumption that $Y_1$ is big, there is another vertex $e \in Y_1$ that is also non-adjacent to $b$ and adjacent to $c$. Then the set $\{d,e,c,4,b,5\}$ induces a $bridge$.\qed

\noindent \textbf{P\ref{obs:c6-yi-xi+2}vv.}\textit{
    A vertex in $X_{i+2}$ is non-adjacent to at most one vertex in $Y_i$. By symmetry, a vertex in $X_{i+5}$ is non-adjacent to at most one vertex in $Y_i$.
}

\noindent {\it Proof}. If there is a vertex $a \in X_3$ that is non-adjacent to two vertices $b,c \in Y_1$, then we claim that every other vertex $d \in X_3$ must be adjacent to exactly one of $\{b,c\}$. Otherwise, if $d$ is non-adjacent to both, then the set $\{a,d,3,4,b,c\}$ induces a $bridge$ or, if $d$ is adjacent to both, then the set $\{b,c,d,4,a,5\}$ induces a $bridge$. WLOG, assume $d$ is non-adjacent to $b$ and adjacent to $c$. Then by the assumption that $X_3$ is big, there is another vertex $e \in X_3$ that is also non-adjacent to $b$ and adjacent to $c$. Then the set $\{b,2,c,3,d,e\}$ induces a $bridge$.\qed

\noindent \textbf{P\ref{obs:c6-yi-xi+2}}\textit{
    If $Y_i \not= \emptyset$, then $X_{i+2} = \emptyset$ and vice versa. By symmetry, if $Y_i \not= \emptyset$, then $X_{i+5} = \emptyset$ and vice versa.
}

\noindent {\it Proof}. From Observations \ref{obs:c6-yi-xi+2}a and \ref{obs:c6-yi-xi+2}vv. and the assumption that $Y_1$ and $X_3$ are big, there exist vertices $a,b \in Y_1$ and $c,d \in X_3$ such that $a,b,c,d$ form a clique. Then the set $\{1,2,a,b,c,d\}$ induces a $bridge$.\qed

\noindent \textbf{P\ref{obs:c6-yi-xi+3}}\textit{
    $Y_{i} \;\circled{0}\; X_{i+3}$ for all $i$. By symmetry, $Y_{i} \;\circled{0}\; X_{i+4}$ for all $i$.
}

\noindent {\it Proof}. If there is a vertex $a \in Y_1$ that is adjacent to a vertex $b \in X_4$, then the set $\{a,1,3,b\}$ induces a $claw$.\qed

\noindent \textbf{P\ref{obs:c6-yi-xi+1-ti}}\textit{
    If $Y_i$ is non-empty, then either $X_{i+1}$ or $T_i$ is empty. By symmetry, if $Y_i$ is non-empty, then either $X_i$ or $T_{i+2}$ is empty.
}

\noindent {\it Proof}. Consider $a \in Y_1$, $b \in X_2$, and $c \in T_1$. By Property P\ref{obs:c6-xi-ti+2}, a choice of $c$ exists such that $b$ and $c$ are not adjacent. By Property P\ref{obs:c6-yi-xi}, $a$ and $b$ are adjacent. By Property P\ref{obs:c6-yi-ti}, $a$ and $c$ are also adjacent. Then the set $\{c,1,a,2,b,3\}$ induces a $bridge$.\qed

\section{Properties of sets in the case of the $C_5$}

\noindent \textbf{P\ref{obs:c5-ti-ti+1}}\textit{
    $T_i \;\circled{0}\; T_{i+1}$ for all $i$. By symmetry, $T_i \;\circled{0}\; T_{i+4}$ for all $i$.
}

\noindent {\it Proof}. If there is a vertex $a \in T_1$ that is adjacent to a vertex $b \in T_2$, then the set $\{1,a,b,3,4,5\}$ induces a $C_6$.\qed

\noindent \textbf{P\ref{obs:c5-ti-ti+2}}\textit{
    A vertex in $T_i$ is adjacent to at most one vertex in $T_{i+2}$ and vice versa. By symmetry, a vertex in $T_i$ is adjacent to at most one vertex in $T_{i+3}$ and vice versa.
}

\noindent {\it Proof}. If there is a vertex $a \in T_1$ that is adjacent to two vertices $b,c \in T_3$, then the set $\{a,b,2,3,c\}$ induces a $C_4$-$twin$.\qed

\noindent \textbf{P\ref{obs:c5-ti-ti+2-ti+3}}\textit{
    A vertex in $T_i$ is adjacent to at most one vertex in $T_{i+2} \cup T_{i+3}$.
}

\noindent {\it Proof}. If there is a vertex $a \in T_1$ that is adjacent to two vertices $b \in T_3$ and $c \in T_4$, then by Property P\ref{obs:c5-ti-ti+1}, $b$ and $c$ are not adjacent and the set $\{a,2,b,c\}$ induces a $claw$.\qed

\noindent \textbf{P\ref{obs:c5-no3cons-t}}\textit{
    There are no 3 consecutive non-empty sets $T_i$.
}

\noindent {\it Proof}. If there are vertices $a \in T_1$ and $b \in T_2$, then we may choose $c \in T_3$ such that by Property P\ref{obs:c5-ti-ti+2}, $c$ is one of the many vertices in $T_3$ that is non-adjacent to $a$. By Property P\ref{obs:c5-ti-ti+1}, $a$ and $b$ are not adjacent and $b$ and $c$ are not adjacent. Then the set $\{a,b,c,5\}$ induces a $4K_1$.\qed

\noindent \textbf{P\ref{obs:c5-xi-xi+1}}\textit{
    A vertex in $X_i$ is non-adjacent to at most one vertex in $X_{i+1}$ and vice versa. By symmetry, a vertex in $X_i$ is non-adjacent to at most one vertex in $X_{i+4}$ and vice versa.
}

\noindent {\it Proof}. If there is a vertex $a \in X_1$ that is non-adjacent to two vertices $b,c \in X_2$, then we claim that every other vertex $d \in X_1$ must be adjacent to exactly one of $\{b,c\}$. Otherwise if $d$ is non-adjacent to both, then the set $\{a,d,2,3,b,c\}$ induces a $bridge$ or, if $d$ is adjacent to both, then the set $\{a,1,d,2,b,c\}$ induces a $bridge$. Then two vertices $d,e \in X_1$ must be adjacent to the same vertex of $\{b,c\}$. WLOG, assume they are both adjacent to $b$ and not adjacent to $c$. Then the set $\{e,d,b,3,c,4\}$ induces a $bridge$.\qed

\noindent \textbf{P\ref{obs:c5-xi-xi+2}}\textit{
    A vertex in $X_i$ is adjacent to at most one vertex in $X_{i+2}$ and vice versa. By symmetry, a vertex in $X_i$ is adjacent to at most one vertex in $X_{i+3}$ and vice versa.
}

\noindent {\it Proof}. If there is a vertex $a \in X_1$ that is adjacent to two vertices $b,c \in X_3$, then the set $\{1,a,b,5,c\}$ induces a $C_4$-$twin$.\qed

\noindent \textbf{P\ref{obs:c5-no3cons-x}}\textit{
    There are no 3 consecutive non-empty sets $X_i$.
}

\noindent {\it Proof}. Consider a vertex $a \in X_1$. From Property P\ref{obs:c5-xi-xi+1}, we may choose $b \in X_2$ such that it is one of the many vertices in $X_2$ that are adjacent to $a$. From Property P\ref{obs:c5-xi-xi+1} and Property P\ref{obs:c5-xi-xi+2}, we may choose $c \in X_3$ such that $c$ is one of the many vertices in $X_3$ that is non-adjacent to $a$ and adjacent to $b$. Then the set $\{a,2,b,3,c,4\}$ induces a $bridge$.\qed

\noindent \textbf{P\ref{obs:c5-xi-xi+1-xi+3}}\textit{
    If $X_i$ and $X_{i+1}$ are non-empty, then $X_{i+3} \;\circled{0}\; X_i \cup X_{i+1}$
}

\noindent {\it Proof}. Consider a vertex $a \in X_1$. From Property P\ref{obs:c5-xi-xi+1}, we may choose $b \in X_2$ such that $b$ is one of the many vertices in $X_2$ that is adjacent to $a$. Consider $c \in X_4$. If $c$ is adjacent to $a$ but not $b$, then the set $\{a,b,4,c,3\}$ induces a $C_4$-$twin$. So if $c$ is adjacent to $a$, it must be adjacent to all of $a$'s neighbours in $X_2$. Since $X_2$ is big and by Property P\ref{obs:c5-xi-xi+1}, $a$ has at least two neighbours in $X_2$. Then so must $c$. This is a contradiction to Property P\ref{obs:c5-xi-xi+2}. Thus $X_1 \;\circled{0}\; X_4$ and, by symmetry, $X_2 \;\circled{0}\; X_4$.\qed

\noindent \textbf{P\ref{obs:c5-ti-xi}}\textit{
    $T_i \;\circled{1}\; X_i$ for all $i$. By symmetry, $T_i \;\circled{1}\; X_{i+4}$ for all $i$.
}

\noindent {\it Proof}. If there is a vertex $a \in T_1$ that is non-adjacent to a vertex $b \in X_1$, then the set $\{1,a,b,5\}$ induces a $claw$.\qed

\noindent \textbf{P\ref{obs:c5-ti-xi+1}}\textit{
    A vertex in $T_i$ is adjacent to at most one vertex in $X_{i+1}$. By symmetry, a vertex in $T_i$ is adjacent to at most one vertex in $X_{i+3}$.
}

\noindent {\it Proof}. If there is a vertex $a \in T_1$ that is adjacent to two vertices $b,c \in X_2$, then we claim that every other vertex $d \in T_1$ must be adjacent to exactly one of $\{b,c\}$. Otherwise, if $d$ is non-adjacent to both, then the set $\{d,1,a,2,b,c\}$ induces a $bridge$, or, if $d$ is adjacent to both, then the set $\{d,a,b,c,3,4\}$ induces a $bridge$. If both $b$ and $c$ have neighbours $d,e \in T_1 - a$ respectively, then there will be another vertex $f \in T_1$ that is adjacent to one of $\{b,c\}$. Then the set $\{b,c,d,e,f\}$ induces a $C_4$-$twin$. So WLOG, $b$ is adjacent to all of $T_1 - a$, while $c$ is non-adjacent to all of $T_1 - a$. So $b$ is adjacent to $d,e \in T_1$. Then the set $\{d,e,b,2,c,3\}$ induces a $bridge$.\qed

\noindent \textbf{P\ref{obs:c5-ti-xi+1}vv}\textit{
    A vertex in $X_{i+1}$ is adjacent to at most one vertex in $T_i$. By symmetry, a vertex in $X_{i+3}$ is adjacent to at most one vertex in $T_i$.
}

\noindent {\it Proof}. If there is a vertex $a \in X_2$ that is adjacent to two vertices $b,c \in T_1$, then we claim that every other vertex $d \in X_2$ must be adjacent to exactly one of $\{b,c\}$. Otherwise, if $d$ is non-adjacent to both, then the set $\{d,3,a,2,b,c\}$ induces a $bridge$, or, if $d$ is adjacent to both, then the set $\{b,c,a,d,3,4\}$ induces a $bridge$. Then by the assumption that $X_2$ is big, at least one of $\{b,c\}$ is adjacent to more than one vertex in $X_2$. Contradiction to Property P\ref{obs:c5-ti-xi+1}.\qed

\noindent \textbf{P\ref{obs:c5-ti-xi+2}}\textit{
    $T_i \;\circled{0}\; X_{i+2}$ for all $i$.
}

\noindent {\it Proof}. If there is a vertex $a \in T_1$ that is adjacent to a vertex $b \in X_3$, then the set $\{b,a,3,5\}$ induces a $claw$.\qed

\noindent \textbf{P\ref{obs:c5-ti-xi-xi+1}}\textit{
    One of $\{T_i, X_i, X_{i+1}\}$ must be empty. By symmetry, one of $\{T_{i+2}, X_i, X_{i+1}\}$ must be empty.
}

\noindent {\it Proof}. If there are vertices $a \in T_1$, $b \in X_1$, we may choose $c \in X_2$ such that it is one of the many vertices that, by Property P\ref{obs:c5-ti-xi+1}, is not adjacent to $a$ and, by Property P\ref{obs:c5-xi-xi+1}, is adjacent to $b$. By Property P\ref{obs:c5-ti-xi}, $a$ and $b$ are adjacent. Then the set $\{a,1,b,2,c,3\}$ induces a $bridge$.\qed

\noindent \textbf{P\ref{obs:c5-ti-ti+1-xi}}\textit{
    One of $\{T_i, X_i, T_{i+1}\}$ must be empty.
}

\noindent {\it Proof}. Consider vertices $a \in T_1$, $b \in X_1$ and $c \in T_2$. By Property P\ref{obs:c5-ti-xi}, $a$ and $b$ are adjacent and $b$ and $c$ are also adjacent. By Property P\ref{obs:c5-ti-ti+1}, $a$ and $c$ are not adjacent. Then the set $\{a,1,b,2,c,3\}$ induces a $bridge$.\qed

\noindent \textbf{P\ref{obs:c5-ti+1-xi-xi+1}}\textit{
    If $T_{i+1}$ is non-empty, $X_i \;\circled{1}\; X_{i+1}$.
}

\noindent {\it Proof}. Consider $a \in T_2$, $b \in X_1$, and $c \in X_2$. By Property P\ref{obs:c5-ti-xi}, $a$ and $b$ are adjacent and $a$ and $c$ are also adjacent. Then if $b$ is not adjacent to $c$, the set $\{1,b,a,c,4,5\}$ induces a $C_6$.\qed

\noindent \textbf{P\ref{obs:c5-ti+1-xi-xi+2}}\textit{
    If $T_{i+1}$ is non-empty, $X_i \;\circled{0}\; X_{i+2}$. By symmetry, if $T_{i+2}$ is non-empty, $X_i \;\circled{0}\; X_{i+2}$.
}

\noindent {\it Proof}. Consider vertices $a \in T_2$, $b \in X_1$, and $c \in X_3$. By Property P\ref{obs:c5-ti-xi}, $a$ and $b$ are adjacent. If $c$ is not adjacent to $a$, then if $c$ is adjacent to $b$, the set $\{b,a,c,1\}$ induces a $claw$. So if $c$ is not adjacent to $a$, then it is not adjacent to $b$.

Suppose $c$ is adjacent to $a$. Then by Property P\ref{obs:c5-ti-xi+1}, $a$ is not adjacent to any other vertex in $X_3$. Suppose $c$ is adjacent to $b$, since otherwise the observation holds. Since $X_3$ is big and by Property P\ref{obs:c5-ti-xi+1}, there are two vertices $d,e \in X_3$ that are not adjacent to $a$. By the previous argument, since $d$ and $e$ are not adjacent to $a$, they are also not adjacent to $b$. Then the set $\{a,b,c,3,d,e\}$ induces a $bridge$. Therefore, $c$ cannot be adjacent to $b$, so the observation holds.\qed

\noindent \textbf{P\ref{obs:c5-xi-ti-ti+2}}\textit{
    If $X_i$ is non-empty, $T_i \;\circled{0}\; T_{i+2}$
}

\noindent {\it Proof}. Consider vertices $a \in X_1$, $b \in T_1$, and $c \in T_3$. By Property P\ref{obs:c5-ti-xi}, $a$ and $b$ are adjacent. Suppose $c$ is not adjacent to $a$. If $b$ and $c$ are adjacent, then the set $\{3,a,b,c,2\}$ induces a $C_4$-$twin$. Now suppose $c$ is adjacent to $a$. By Property P\ref{obs:c5-ti-xi+1}, there is another vertex $d \in X_1$ that is not adjacent to $c$. If $b$ and $c$ are adjacent, then the set $\{b,d,3,c,2\}$ induces a $C_4$-$twin$. Thus, if a vertex $a \in X_1$ exists, any vertices $b \in T_1$ and $c \in T_3$ must be non-adjacent.\qed

\noindent \textbf{P\ref{obs:c5-xi-xi+2-ti+4}}\textit{
    If a vertex in $X_i$ is adjacent to vertices in $X_{i+2}$ and in $T_{i+4}$, then its neighbours in $X_{i+2} \cup T_{i+4}$ are a clique.
}

\noindent {\it Proof}. Consider vertices $a \in X_1$, $b \in X_3$ and $c \in T_5$. If $a$ is adjacent to both $b$ and $c$ but $b$ and $c$ are not adjacent, then the set $\{a,2,b,c\}$ induces a $claw$. So if a vertex in $X_1$ is adjacent to a vertex in $X_3$ and one in $T_5$, those vertices in $X_3$ and $T_5$ must also be adjacent.\qed

\noindent \textbf{P\ref{obs:c5-xi-xi+1-ti+1-clique}}\textit{
    If $X_i$, $X_{i+1}$, $T_{i+1}$ are all non-empty, then $X_i \cup X_{i+1} \cup T_{i+1}$ is a clique.
}

\noindent {\it Proof}. Consider the case when the sets $X_1$, $X_2$, and $T_2$ are non-empty. By Property P\ref{obs:c5-ti-xi}, $T_2 \cup X_2$ and $T_2 \cup X_1$ are cliques. By Property P\ref{obs:c5-ti+1-xi-xi+1}, $X_1 \cup X_2$ is also a clique. Thus, the union of the three sets forms a clique.\qed

\end{appendices}
\end{document}